\documentclass[12pt]{article}
\usepackage{amssymb}
\usepackage{tikz} 
\usepackage{mathrsfs}
\font\eufm=eufm10 at 12pt\font\eufms=eufm10\font\eufmss=eufm7\newfam\eufam
\textfont\eufam=\eufm\scriptfont\eufam=\eufms\scriptscriptfont\eufam=\eufmss
\def\build#1_#2^#3{\mathrel{\mathop{\kern 0pt#1}\limits_{#2}^{#3}}}\def\rde{\mathscr}
\def\Z{{\bf Z}}\def\R{{\rde R}}\def\P{{\rde P}}\def\Mv{{\rde M}{\rde V}}\def\V{{\rde V}}\def\E{{\rde E}}\def\F{{\rde F}}\def\G{{\rde G}}

\def\Hom{{\rm Hom}}\def\fl{\longrightarrow}\def\Ker{{\rm Ker}}
\def\A{{\rde A}}\def\B{{\rde B}}\def\C{{\rde C}}\def\D{{\rde D}}\def\Nil{{\rde N}\! il}\def\S{{\mathscr S}}
\def\rde{\mathscr}
\def\vc{{\scriptstyle\wedge }}
\def\cqfd{\hfill\vbox{\hrule\hbox{\vrule height6pt depth0pt\hskip 6pt \vrule height6pt}\hrule\relax}}
\def\noi{\noindent}\def\e{{\varepsilon}}
\def\hfl#1#2{\smash{\mathop{\hbox to 12 mm{\rightarrowfill}}\limits^{\scriptstyle#1}_{\scriptstyle#2}}}
\def\vfl#1#2{\llap{$\scriptstyle #1$}\left\downarrow\vbox to 6mm{}\right.\rlap{$\scriptstyle #2$}}
\def\bhfl#1#2{\smash{\mathop{\vbox{\hbox to 10 mm{\rightarrowfill}\nointerlineskip\hbox to 10 mm{\rightarrowfill}}}\limits^{\scriptstyle#1}_{\scriptstyle#2}}}

\def\diagram#1{\def\normalbaselines{\baselineskip=0pt\lineskip=10pt\lineskiplimit=1pt} \matrix{#1}}
\def\pv{\raise 2pt\hbox{$\bigwedge$}}\def\v{{}^\vee}
\begin{document}

\overfullrule=0pt

\vskip 64pt
\centerline{\bf Algebraic K-theory of generalized free products and functors Nil.}
\vskip 12pt
\centerline{Pierre Vogel\footnote{Universit\'e de Paris, Institut de Math\'ematiques de Jussieu-Paris Rive Gauche,
    B\^atiment Sophie Germain, Case 7012, 75205--Paris Cedex 13 France, Email: pierre.vogel@imj-prg.fr}}
\vskip 48pt
\noi{\bf Abstract.} In this paper, we extend Waldhausen's results on algebraic K-theory of generalized free products in a more general setting and we give some properties
of the Nil functors. As a consequence, we get new groups with trivial Whitehead groups.
\vskip 12pt
\noi{\bf Keywords:} Algebraic K-theory, functor Nil, Whitehead groups.

\noi{\bf Mathematics Subject Classification (2020):} 18E10, 19D35, 19D50.
\vskip 24pt
\noi{\bf Introduction.}
\vskip 12pt
Quillen's construction associates to any essentially small exact category $\A$ its algebraic K-theory which is an infinite loop space $K(\A)$ and this
correspondence is a functor from the category of essentially small exact categories to the category $\Omega sp_0$ of infinite loop spaces.

If $\A$ is the category $\P_A$ of finitely generated projective right modules over some ring $A$, one gets a functor $A\mapsto K(A)=K(\P_A)$ and this functor can be
enriched into a new functor $\underline K$ containing also the negative part of the algebraic K-theory. That is $\underline K$ is a functor from the category
of rings to the category $\Omega sp$ of $\Omega$-spectra and the natural transformation $K(A)\rightarrow \underline K(A)$ induces a homotopy equivalence from $K(A)$ to the
0-th term of $\underline K(A)$.

By a left-flat bimodule we mean a pair $(A,S)$ where $A$ is a ring and $S$ is an $A$-bimodule flat on the left. The left-flat bimodules form a category where a morphism
$(A,S)\fl(B,T)$ is a ring homomorphism $f:A\fl B$ together with a morphism of $A$-bimodules $\varphi:S\fl f^*(T)$.

For each left-flat bimodule $(A,S)$ one has an exact category $\Nil(A,S)$ where the objects are the pairs $(M,f)$ where $M$ is an object in $\P_A$ and
$f:M\rightarrow M\otimes_A S$ is a nilpotent morphism of right $A$-modules.

The correspondence $M\mapsto (M,0)$ induces a morphism $K(A)\rightarrow K(\Nil(A,S))$ and this morphism has a retraction coming from the functor $(M,f)\mapsto M$. Thus
there is a functor $Nil$ from the category of left-flat bimodules to $\Omega sp_0$ which is unique up to homotopy such that:
$$K(\Nil(A,S))\simeq K(A)\times Nil(A,S)$$
\vskip 12pt
\noi{\bf Theorem 1.} {\sl There is a functor $\underline N il$ from the category of left-flat bimodules to the category $\Omega sp$ of $\Omega$-spectra and a natural
transformation $Nil\rightarrow\underline Nil$ such that the following holds for every left-flat bimodule $(A,S)$:

$\bullet$ the map $Nil(A,S)\rightarrow\underline Nil(A,S)$ induces a homotopy equivalence from $Nil(A,S)$ to the 0-th term of $\underline Nil(A,S)$

$\bullet$ if $R$ is the tensor algebra of $S$, then there is a natural  homotopy equivalence in $\Omega sp$:
$$\underline K(R)\build\longrightarrow_{}^\sim \underline K(A)\times\Omega^{-1}(\underline Nil(A,S))$$

Moreover if $A$ is regular coherent on the right, every spectrum $\underline Nil(A,S)$ is contractible.}
\vskip 12pt
Following a terminology of Waldhausen, a ring homomorphism $\alpha:A\rightarrow B$ will be called pure if it is split injective as an $A$-bimodule homomorphism.
\vskip 12pt
\noi{\bf Theorem 2.} {\sl Let $\alpha:C\longrightarrow A$ and $\beta: C\longrightarrow B$ be pure ring homomorphisms. Let $R$ be the ring defined by the push-out diagram:
$$\diagram{C&\hfl{\alpha}{}&A\cr\vfl{\beta}{}&&\vfl{}{}\cr B&\hfl{}{}&R\cr}$$
and $\underline K'(R)$ be the $\Omega$-spectrum defined by the homotopy fibration in $\Omega sp$:
$$\underline K(C)\build\longrightarrow_{}^f\underline K(A)\times\underline K(B)\longrightarrow\underline K'(R)$$
where $f$ is the map $\underline K(\alpha)\times -\underline K(\beta)$.

Suppose $A$ and $B$ are $C$-flat on the left. Then there exist a left-flat bimodule $(C\times C,S)$ and a homotopy equivalence in $\Omega sp$:
$$\underline K(R)\build\longrightarrow_{}^\sim \underline K'(R)\times\Omega^{-1}(\underline Nil(C\times C,S))$$}
\vskip 12pt
\noi{\bf Theorem 3.} {\sl Let $C$ and $A$ be two rings and $\alpha$ and $\beta$ be two pure ring homomorphisms from $C$ to $A$. Let $R$ be the ring generated by $A$ and
an invertible element $t$ with the only relations:
$$\forall c\in C,\ \ \alpha(c)t=t\beta(c)$$
and $\underline K'(R)$ be the $\Omega$-spectrum defined by the homotopy fibration:
$$\underline K(C)\build\longrightarrow_{}^f\underline K(A)\longrightarrow\underline K'(R)$$
where $f$ is the map $\underline K(\alpha)-\underline K(\beta)$.

Suppose $A$ is $C$-flat on the left via both $\alpha$ and $\beta$. Then there exist a left-flat bimodule $(C\times C,S)$ and a homotopy equivalence in $\Omega sp$:
$$\underline K(R)\build\longrightarrow_{}^\sim \underline K'(R)\times\Omega^{-1}(\underline Nil(C\times C,S))$$}
\vskip 12pt
The connective part of these theorems are generalization of the Waldhausen's results in [Wa1], in the sense that the condition free on the left in [Wa1] may be replaced
by the condition flat on the left and (in the polynomial extension case) the condition on the right may be removed.
\vskip 24pt
Because of these results, the functor $\underline Nil$ detects in some sense the default of excision in algebraic K-theory. So it would be useful to know when the
spectrum $\underline Nil(A,S)$ is contractible, especially if $A$ is not regular coherent. Actually we have the following result:
\vskip 12pt
\noi{\bf Theorem 4.} {\sl Let $A$ and $B$ be two rings, $S$ be an $(A,B)$-bimodule and $T$ be a $(B,A)$-bimodule. Suppose $S$ and $T$ are flat on both sides.
Using projections $A\times B\fl A$ and $A\times B\fl B$, $S$ and $T$ may be considered as $A\times B$-bimodules. Then we have natural homotopy equivalences of spectra:
  $$\underline Nil(A\times B,S\oplus T)\build\fl_{}^\sim \underline Nil(A,S\build\otimes_B^{} T)$$
  $$\underline Nil(A\times B,S\oplus T)\build\fl_{}^\sim \underline Nil(B,T\build\otimes_A^{} S)$$}
\vskip 12pt
We have other results concerning the spectrum $\underline Nil(A,S)$ when the bimodule $S$ is a direct sum of bimodules:
$S=\build\oplus_{i\in I}^{} S_i$.

Let $W(I)$ be the set of words in the set $I$. This set is the unitary monoid freely generated by $I$. Let $CW(I)$ the set of cyclic words in $I$. The set $CW(I)$
is the quotient of $W(I)$ by the equivalence relation $uv\sim vu$ in $W(I)$. A word $u\in W(I)$ is said to be reduced if we have the following:
$$\forall v\in W(I),\ \ \forall p>1,\ \ \ u\not=v^p$$
The set of reduced words is denoted by $W_0(I)$ and its image in $CW(I)$ is denoted by $CW_0(I)$. A subset $X\subset W(I)$ is said to be admissible if the projection
$W(I)\fl CW(I)$ induces a bijection $X\build\fl_{}^\sim CW_0(I)$.

For every $u\in W(I)$, we have a bimodule $S_u$ defined by:
$$u=1\ \Longrightarrow\ \ S_u=A$$
$$u=vi,\ \hbox{with}\ \ i\in I,\ \ \Longrightarrow\ \ S_u=S_v\build\otimes_A^{}S_i$$
\vskip 12pt
Using these notations, we have this result:
\vskip 12pt
\noi{\bf Theorem 5:} {\sl Let $A$ be a ring and $S_i$, $i\in I$, be a family of $A$-bimodules. Suppose each bimodule $S_i$ is flat on both sides. Let $X$ be an
  admissible set in $W(I)$. Then we have a
  homotopy equivalence of spectra:
  $$\underline Nil(A,\build\oplus_i^{} S_i)\build\fl_{}^\sim\build\oplus_{u\in X}^{} \underline Nil(A,S_u)$$
  where $\oplus$ is the coproduct in the category of spectra.}
\vskip 12pt
Moreover, using theorems 4 and 5, we can deduce this last result:
\vskip 12pt
\noi{\bf Theorem 6:} {\sl Let $A$ be a ring, $S$ be an $A$-bimodule and $I$ be a set. For each $i\in I$, let $A_i$ be a right regular coherent ring, $E_i$ be an 
$(A,A_i)$-bimodule and $F_i$ be an $(A_i,A)$-bimodule. Suppose all these bimodules are flat on both sides. Then the inclusion:
$$S\subset S\ \oplus \build\oplus_i^{} E_i\build\otimes_{A_i}^{} F_i$$
induces a homotopy equivalence:
$$\underline Nil(A,S)\build\fl_{}^\sim  \underline Nil(A,S\ \oplus \build\oplus_i^{} E_i\build\otimes_{A_i}^{} F_i)$$}
\vskip 12pt
The paper is organized as follows:

In section 1, we construct many categories and functors in such a way they are defined in each of the three cases: the general polynomial extension case, the generalized
free product case and the generalized Laurent extension case. We prove also many algebraic properties of these categories and functors.

In section 2, we apply these properties to algebraic K-theory and prove theorems 1, 2 and 3.

The section 3 is devoted to the proof of theorems 4, 5 and 6 about Nil functors.

In the last section we apply all these theorems and get new results about Whitehead spectra. In particular we construct a class Cl$_1$ bigger that Waldhausen's class
Cl such that every group in Cl$_1$ has trivial Whitehead groups.
\vskip 24pt

\noi{\bf 1. The categories $\V$ and $\Mv$ and their algebraic properties.}
\vskip 12pt
In order to simplify the notations, the following writing conventions will often be used:
\vskip 12pt
$\bullet$ Convention 1: if $\Phi:\A\rightarrow\B$ is a functor, then for every morphism $\alpha$ in $\A$, its image under $\Phi$ will be still denoted by $\alpha$. So, if
$\alpha:X\rightarrow Y$ is a morphism, we have a morphism $\alpha:\Phi(X)\rightarrow\Phi(Y)$.
\vskip 12pt
$\bullet$ Convention 2: if $E$ is a right module over some ring $A$ and $F$ is a left module over the same ring, then the module $E\otimes_A F$ will be simply denoted be 
$EF$. In the same spirit, if $E$ is an $A$-bimodule, the tensor product $E\otimes_A\dots\otimes_A E$ of $n$ copies of $E$ will be denoted by $E^n$.
\vskip 12pt
\noi{\bf 1.1 Definition:} {\sl Let $A$ be a ring and $S$ be an $A$-bimodule. Let $M$ be a right $A$-module and $f:M\rightarrow MS$ be an $A$-linear map. So by iteration, we
  get for each integer $n>0$ a morphism $f^n:M\rightarrow MS^n$. We say that $f$ is nilpotent if every element in $M$ is killed by some power $f^n$ of $f$.}
\vskip 12pt
\noi{\bf 1.2 Lemma:} {\sl  Let $(A,S)$ be a left-flat bimodule. Let $M$ be a right $A$-module and $f:M\rightarrow MS$ be an $A$-linear map. Then $f$ is nilpotent if and
only if there is a filtration of $M$:
$$0=M_0\subset M_1\subset M_2\subset\dots  \subset M$$
by right $A$-submodules, such that:

$\bullet$ $M$ is the union of the $M_i$'s

$\bullet$ for every $i>0$, one has: $f(M_i)\subset M_{i-1}S$.}
\vskip 12pt
\noi{\bf Proof:} If such a filtration exists, then $f$ is clearly nilpotent.

Suppose $f$ is nilpotent. For every integer $n\geq0$, denote by $M_n$ the kernel of $f^n:M\fl MS^n$. By construction we have an increasing sequence
$$0=M_0\subset M_1\subset M_2\subset M_3\dots$$
and $M$ is the union of the $M_i$'s.

Since $S$ is flat on the left, we have, for every $n\geq0$, an isomorphism:
$$\Ker(f^n:MS\fl MS^{n+1})\simeq \Ker(f^n:M\fl MS^n)S$$
and then an equality: $M_{n+1}=f^{-1}(M_n S)$. The result follows.\cqfd
\vskip 12pt
\noi{\bf 1.3 The exact category $\Nil(A,S)$ and the space $Nil(A,S)$.}
\vskip 12pt
Let $(A,S)$ be a left-flat bimodule. The pairs $(M,f)$ where $M$ is a right $A$-module and $f:M\rightarrow MS$ is a nilpotent morphism are the objects of a category
denoted by $\Nil(A,S)\v$. Let $\A$ be the category of finitely generated projective right $A$-modules.
The full subcategory of $\Nil(A,S)\v$ generated by pairs $(M,f)$ with $M\in\A$ will be denoted by $\Nil(A,S)$.
\vskip 12pt
If $0\fl (M,f)\fl (M',f')\fl (M'',f'')\fl 0$ is a sequence
in $\Nil(A,S)$, we say that this sequence is exact if the following diagram is commutative with exact lines:
$$\diagram{0&\hfl{}{}&M&\hfl{}{}&M'&\hfl{}{}&M''&\hfl{}{}&0\cr &&\vfl{f}{}&&\vfl{f'}{}&&\vfl{f''}{}&&\cr 0&\hfl{}{}&MS&\hfl{}{}&M'S&\hfl{}{}&M''S&\hfl{}{}&0\cr}$$
With these exact sequences the category $\Nil(A,S)$ becomes an exact category in the sense of Quillen and its algebraic K-theory $K(\Nil(A,S))$ is a well defined infinite
loop space (see [Q]).

Actually $\Nil$ is a functor from the category of left-flat bimodules to the category of essentially small exact categories and exact functors.

We have two exact functors $M\mapsto (M,0)$ from $\A$ to $\Nil(A,S)$ and $(M,f)\mapsto M$
from $\Nil(A,S)$ to $\A$ inducing two maps:
$$K(A)\build\longrightarrow_{}^F K(\Nil(A,S))\build\longrightarrow_{}^G K(A)$$
where $G$ is a retraction of $F$. Denote by $Nil(A,S)$ the homotopy fiber of $G$. Then $Nil(A,S)$ is an infinite loop space and we have a decomposition:
$$K(\Nil(A,S))\simeq K(A)\times Nil(A,S)$$
\vskip 12pt
Throughout this paper, we'll consider many categories and functors and, in particular, many exact categories and their abelianizations, where an abelianization
of an exact category is defined as follows:
\vskip 12pt
\noi{\bf Definition:} {\sl Let $\E$ be an exact category. We say that a category $\E\v$ is an abelianization of $\E$ if the following holds:
 
   $\bullet$ $\E\v$ is an abelian category

  $\bullet$ $\E$ is a fully exact subcategory of $\E\v$ i.e. $\E$ is a full subcategory of $\E\v$ and, for every sequence $S=(0\fl X\fl Y\fl Z\fl0)$ in $\E$, $S$ is exact
in $\E$ if and only if $S$ is exact un $\E\v$
 
   $\bullet$ $\E$ is stable in $\E\v$ under extension.}
\vskip 12pt
Notice that the Gabriel-Quillen embedding theorem produces an abelianization for every essentially small exact category (see [TT], thm A.7.1 or [K], prop A.2).
\vskip 12pt
Following Waldhausen we consider three situations: the generalized free product case, the generalized Laurent extension case and the generalized polynomial extension case.

In case 1 (i.e. the generalized polynomial extension case), we have a ring $C$ and a $C$-bimodule $S$ which is flat on the left. In this case the ring $R$ is the
tensor algebra of $S$:
$$R=C\oplus S\oplus S^2\oplus S^3\oplus\dots$$

In case 2 (i.e. the generalized free product case) we have two pure morphisms of rings $\alpha: C\rightarrow A$ and $\beta: C\rightarrow B$ and we suppose that
$A$ and $B$ are $C$-flat on the left. We denote by $R$ the ring defined by the cocartesian diagram:
$$\diagram{C&\hfl{\alpha}{}&A\cr\vfl{\beta}{}&&\vfl{}{}\cr B&\hfl{}{}&R\cr}$$

In case 3 (i.e the generalized Laurent extension case) we have two rings $C$ and $A$ and two pure morphisms $\alpha$ and $\beta$ from $C$ to $A$. We suppose that
$A$ is $C$-flat on the left via both $\alpha$ and $\beta$ and we denote by $R$ the ring generated by $A$ and an invertible element $t$ with the only relations:
$$\forall c\in C,\ \ \alpha(c)t=t\beta(c)$$
So we have a morphism $\gamma:A\rightarrow R$ and $\gamma\circ\alpha$ and $\gamma\circ\beta$ are conjugate.

From now on we will consider $\alpha$ and $\gamma$ as inclusions. So in each case $A$, $B$ and $R$ are $C$-bimodules.
\vskip 12pt
We denote by $\A$, $\B$, $\C$ and $\R$ the categories of finitely generated projective right modules over the rings $A$, $B$, $C$ and $R$ respectively.
We set also: $\D=\C$ in case 1, $\D=\A\times\B$ in case 2 and $\D=\A$ in case 3. This category is the category of finitely generated projective right modules over the ring
$C$ or $A\times B$ or $A$.

The categories  $\A$, $\B$, $\C$, $\D$ and $\R$
are contained in the abelian categories $\A\v$, $\B\v$, $\C\v$, $\D\v$ and $\R\v$ of right-modules over the corresponding rings and these categories are abelianizations
of $\A$, $\B$, $\C$, $\D$ and $\R$ respectively. Notice that $\Nil(A,S)\v$ is also an abelianization of $\Nil(A,S)$.

We denote by $C_2$ the ring $C$ in case 1 and $C\times C$ in the other cases and also by $\C_2\v$ the category of right $C_2$-modules and by $\C_2$ the subcategory of
$\C_2\v$ generated by finitely generated projective modules. The category $\C_2\v$ is also an abelianization of $\C_2$.

We will define the $C_2$-bimodule $S$ and many categories and functors in order to give a common proof of theorems 1, 2 and 3 (at least for the connective part of it).
\vskip 12pt
\noi{\bf 1.4 The functors $s$ and $\sigma$ and the bimodule $S$.}
\vskip 12pt
Consider the case $C_2=C\times C$. Let $\pi_1$ and $\pi_2$ be the two projections
$C\times C\rightarrow C$. Consider a right $C$-module $M$ and an integer $i\in\{1,2\}$. The ring $C\times C$ acts on $M$ via $\pi_i$ and becomes a right
$C\times C$-module $M^i$. This functor $M\mapsto M^i$ from $\C\v$ to $\C_2\v$ has an adjoint functor (on both sides) from $\C_2\v$ to $\C\v$ denoted by $E\mapsto E_i$.
The two functors $M\mapsto M^1$ and $M\mapsto M^2$ induce an equivalence of categories from $\C\v\times\C\v$ to $\C_2\v$ and the functors $E\mapsto E_1$ and
$E\mapsto E_2$ induce an inverse of it.

The two functors $s_i: M\mapsto M^i$ from $\C$ to $\C_2$ (and also from $\C\v$ to $\C_2\v$) are exact and it is the same for they adjoint functors
$\sigma_i:E\mapsto E_i$ from $\C_2$ to $\C$ (and from $\C_2\v$ to $\C\v$). So $s_1$, $s_2$ and $s=s_1\oplus s_2$ are exact functors and their adjoint functors $\sigma_1$,
$\sigma_2$ and $\sigma=\sigma_1\oplus \sigma_2$ are also exact.

In case 1, $s$ and $\sigma$ are defined to be identities. Therefore $s$ and $\sigma$ are well defined in all cases: $s$ is an exact functor from $\C$ to $\C_2$ (and
also from $\C\v$ to $\C_2\v$) and $\sigma$ is an exact functor from $\C_2$ to $\C$ (and also from $\C_2\v$ to $\C\v$).

Moreover, for every module $E$ in $\C_2$ (or $\C_2\v$) the module $\sigma(E)$ is nothing else but the module $E$ equipped with the $C$-action induced by the identity or
the diagonal map from $C$ to $C_2$.
\vskip 12pt
We can do the same for left modules and we have functors $M\mapsto {}^i M$ and $E\mapsto {}_i E$.
In case of bimodules, we gets functors $M\mapsto {}^i M^j$ and $E\mapsto {}_i E_j$. Using these notations we have the following, for every right $C\times C$-module
$E$ and left $C\times C$-module $F$:
$$E\build\otimes_{C\times C}^{}F=EF=E_1\ {}_1 F\oplus E_2\ {}_2 F=\ \build\oplus_i^{} E_i\ {}_i F$$

In case 1, the $C_2$-bimodule $S$ is already defined.

Consider the case 2. Since $\alpha:C\rightarrow A$ and $\beta:C\rightarrow B$ are pure morphisms, $A$ and $B$ have two complements $A'$ and $B'$. These objects are
$C$-bimodules and we have two decompositions of $C$-bimodules: $A=\alpha(C)\oplus A'$ and $B=\beta(C)\oplus B'$. Then we define the bimodule $S$ by:
$${}_2S_1=A'\hskip 24pt {}_1 S_2=B'\hskip 24pt {}_1 S_1={}_2 S_2=0$$

Consider the case 3. Ring homomorphisms $\alpha$ and $\beta$ induce two left $C$-actions on $A$ and we get two $(C,A)$-bimodules ${}_\alpha A$ and ${}_\beta A$.
By doing the same on the right, we get four $C$-bimodules ${}_\alpha A_\alpha $, ${}_\alpha A_\beta$, ${}_\beta A_\alpha$ and ${}_\beta A_\beta$. Moreover, since
$\alpha$ and $\beta$ are pure morphisms, we have two decompositions of $C$-bimodules:
$${}_\alpha A_\alpha=\alpha(C)\oplus A'\hskip 24pt {}_\beta A_\beta=\beta(C)\oplus A''$$
Then we define the bimodule $S$ by:
$${}_2S_1=A'\hskip 24pt {}_1 S_2=A''\hskip 24pt {}_1 S_1={}_\beta A_\alpha\hskip 24pt {}_2 S_2={}_\alpha A_\beta$$

Then in the three cases $S$ is a well defined $C_2$-bimodule. It is easy to see that $S$ is flat on the left.
\vskip 12pt
\noi{\bf 1.5 The categories $\D$, ${\Mv}$ and $\V$ and the functors $T$, $F$ and $\widehat F$.}
\vskip 12pt

We have a functor $T:\D\v\longrightarrow \R\v$ defined as follows:

$\bullet$ in case 1: $T(E)=E\otimes_C R=ER$

$\bullet$ in case 2: $T(E_A,E_B)=E_A\otimes_A R\oplus E_B\otimes_B R=E_A\ R\oplus E_B\ R$

$\bullet$ in case 3: $T(E)=E\otimes_A R=ER$

It is easy to check that $T$ is an exact functor sending $\D$ to $\R$.
\vskip 12pt
Let $E$ be an object in $\D\v$, $M$ be an object in $\C\v$ and $\varphi:T(E)\fl MR$ be a morphism in $\R\v$. We say that $\varphi$ is admissible if the following holds:

$\bullet$ $\varphi(E)\subset M\oplus MS$ in case 1

$\bullet$ $\varphi(E_A)\subset MA$ and $\varphi(E_B)\subset MB$ in case 2 (with: $E=(E_A,E_B)$)

$\bullet$ $\varphi(E)\subset MA\oplus MtA\subset MR$ in case 3.

The set of admissible morphisms $\varphi:T(E)\fl MR$ will be denote by $\F(E,M)$.
\vskip 12pt

Following Waldhausen, we define a splitting diagram as a triple $X=(E,M,\varphi)$ with: $E\in\D\v$, $M\in\C\v$ and $\varphi\in{\rde F}(E,M)$.

The splitting diagram $(E,M,\varphi)$ is called a Mayer Vietoris presentation (resp. a splitting module) if $\varphi$ is surjective (resp. bijective). The splitting
modules, the Mayer Vietoris presentations and the splitting diagrams define three categories $\V\v\subset\overline{\Mv}\subset \S\v$. Moreover categories $\V\v$ and
$\S\v$ are abelian.

If we replace $\D\v$ and $\C\v$ by $\D$ and $\C$, we get three subcategories $\V\subset \Mv\subset \S$.

The correspondences $(E,M,\varphi)\mapsto E$ and $(E,M,\varphi)\mapsto M$ define two functors $\Phi_2:\S\fl \D$ and $\Phi_3:\S\fl \C$ (and also from
$\S\v$ to $\D\v$ and from $\S\v$ to $\C\v$). We have an extra functor $\Phi_1$ sending $(E,M,\varphi)$ to the kernel of $\varphi$.

Consider a sequence $S=(0\fl X\fl Y\fl Z\fl 0)$ in $\V$ or in $\S$. We say that this sequence is exact if it is sent to an exact sequence under $\Phi_2$ and
$\Phi_3$. If $S$ is a sequence in $\Mv$, we say that $S$ is exact if it is sent to an exact sequence under $\Phi_1$, $\Phi_2$ and $\Phi_3$.

With these exact sequences, $\V$, $\Mv$ and $\S$ become exact categories and the inclusions $\V\subset\Mv\subset\S$ are exact functors. Moreover
$\Phi_1:\Mv\fl\R$ is an exact functor. In some sense $\V$ is the kernel of the functor $\Phi_1:\Mv\fl\R$.

Since $\overline{\Mv}$ is not an abelian category, it will be useful to construct an abelian category $\Mv\v$ containing $\Mv$. The category $\Mv$ is equivalent to the
category of tuples $(U,E,M,\mu,\varphi)$ where $(U,E,M)$ is an object of $\R\times\D\times\C$ and $\mu:U\fl T(E)$ and $\varphi:T(E)\fl MR$ are morphisms in $\R$ such that
$\varphi$ is admissible and the sequence:
$$0\fl U\build\fl_{}^\mu T(E)\build\fl_{}^\varphi MR\fl 0$$
is exact. So the category $\Mv\v$ is defined as the category of tuples $(U,E,M,\mu,\varphi)$ where $(U,E,M)$ is an object of $\R\v\times\D\v\times\C\v$ and
$\mu:U\fl T(E)$ and $\varphi:T(E)\fl MR$ are morphisms in $\R\v$ such that $\varphi$ is admissible and $\varphi\mu=0$. It is easy to see that $\Mv\v$ is an abelian category
and the inclusions of exact categories $\V\subset\Mv$ and $\Mv\subset\S$ extend to functors $\V\v\fl\Mv\v$ and $\Mv\v\fl\S\v$.  Moreover the categories $\V\v$, $\Mv\v$ and
$\S\v$ are abelianizations of $\V$, $\Mv$ and $\S$ respectively.

\vskip 12pt
We have a functor $F$ from $\C_2\v$ to $\D\v$ defined as follows:

$\bullet$ in case 1, $F$ is the identity

$\bullet$ in case 2, $F(M)=(\sigma_1(M)A,\sigma_2(M)B)\in\D\v$

$\bullet$ in case 3, $F(M)=\sigma_1(M){}_\alpha A\oplus\sigma_2(M){}_\beta A\in\D\v$

\noi where ${}_\alpha A$ and ${}_\beta A$ are the module $A$ equipped with the $(C,A)$-bimodule structure induced by $\alpha$ and $\beta$ respectively.

This functor $F$ is exact and sends $\C_2$ to $\D$. It has a right adjoint functor $\widehat F$ from $\D\v$ to $\C_2\v$ and we have:

$\bullet$ in case 1, $\widehat F$ is the identity

$\bullet$ in case 2, $\widehat F(E_A,E_B)=s_1(E_A)\oplus s_2(E_B)$

$\bullet$ in case 3, $\widehat F(E)=s_1(E_\alpha)\oplus s_2(E_\beta)$

\noi where $E_\alpha$ and $E_\beta$ are the module $E$ equipped wihe the right $C$-module structure induced by $\alpha$ and $\beta$ respectively.
\vskip 12pt
In case: $C_2=C\times C$, we have another functor $M\mapsto \widetilde M$ from $\C_2\v$ (or $\C_2$) to itself defined by:
$$M=(M',M'')\ \ \Longrightarrow\ \ \widetilde M=(M'',M')$$
\vskip 12pt
\noi{\bf 1.6 Lemma:} {\sl Suppose $C_2=C\times C$. Then for every right $C\times C$-module $M$, we have natural isomorphisms:
$$s\sigma(M)\simeq M\oplus\widetilde M$$
$$\widehat F F(M)\simeq M\oplus \widetilde MS$$
Moreover the induced projection $s\sigma(M)\fl M$ and the induced injections $M\fl s\sigma(M)$ and $M\fl \widehat F F(M)$ are adjoint to identities.}
\vskip 12pt
\noi{\bf Proof:} Let $M=(M',M'')$ be a module in $\C_2=\C\times\C$.  We have:
$$s(\sigma(M))=s(M'\oplus M'')=s_1(M'\oplus M'')\oplus s_2(M'\oplus M'')\simeq M\oplus(s_1(M'')\oplus s_2(M'))\simeq M\oplus \widetilde M$$

In case 2, we have:
$$\widehat F F(M)=\widehat F(M'A,M''B)=(M'(C\oplus A'),M''(C\oplus B'))$$
$$\simeq (M',M'')\oplus (M'{}_2 S_1,M''{}_1S_2)\simeq M\oplus \widetilde M S$$

In the case 3 we have:
$$\widehat F F(M)=\widehat F(M'{}_\alpha A\oplus M''{}_\beta A)=s_1(M'{}_\alpha A_\alpha\oplus M''{}_\beta A_\alpha)\oplus s_2(M'{}_\alpha A_\beta\oplus M''{}_\beta
A_\beta)$$
$$=s_1(M'\oplus M'{}_2S_1\oplus M''{}_1 S_1)\oplus s_2(M'{}_2 S_2\oplus M''\oplus M''{}_1S_2)$$
$$\simeq M\oplus s_1(\widetilde M S_1)\oplus s_2(\widetilde M S_2)\simeq M\oplus \widetilde MS$$
and the result follows.\cqfd
\vskip 12pt
\noi{\bf 1.7 The module $M[S]$ and the transformations $e$, $\e$ and $\tau$.}
\vskip 12pt
For every $M\in\C_2\v$, we set:
$$M[S]=M\oplus MS\oplus MS^2\oplus\dots$$
and $M[S]$ is a right $C_2$-module. In case 1, $M[S]$ is isomorphic to $MR$.

We have a stabilization map $MS[S]\fl M[S]$ induced by the identities $MSS^i\fl MS^{i+1}$.
\vskip 12pt

\noi{\bf 1.8 Lemma:} {\sl There exist natural transformations:
  $$e:TF(P)\build\longrightarrow_{}^\sim \sigma(P)R$$
  $$\e:\sigma(\widehat F(E)[S])\build\longrightarrow_{}^\sim T(E)$$
  $$\tau:\sigma(s(M)[S])\longrightarrow MR$$
  for all  $P\in\C_2\v$, $E\in\D\v$ and $M\in\C\v$ such that:

  $\bullet$ $e$ is an isomorphism of $R$-modules

  $\bullet$ $\e$ is an isomorphism of $C$-modules

  $\bullet$ $\tau$ is an epimorphism of $C$-modules and the following diagram is exact (i.e. cartesian and cocartesian):
  $$\diagram{\sigma s(M)&\hfl{}{}&\sigma(s(M)[S])\cr\vfl{\tau_0}{}&&\vfl{}{\tau}\cr M&\hfl{}{}&MR\cr}$$
  where the horizontal maps are the canonical inclusions and $\tau_0$ is adjoint to the identity of $s(M)$.}
\vskip 12pt
\noi{\bf Proof:} In case 1, it is easy to see that $e$, $\e$, $\tau$ and $\tau_0$ can be chosen to be identities.

In the other cases consider a module $P\in\C_2\v$. So we have two $C$-modules $M=P_1$ and $N=P_2$. 

In case 2, we have:
$$T(F(P))=T(MA,NB)=MA\otimes_A R\oplus NB\otimes_B R\simeq MR\oplus NR\simeq(M\oplus N)R=\sigma(P)R$$
and we get the isomorphism $e$.

In case 3, we have:
$$T(F(P))=T(M{}_\alpha A\oplus N{}_\beta A)\simeq  M{}_\alpha R\oplus N{}_\beta R$$
but the multiplication on the left by $t$ induces a isomorphism of $(C,R)$-bimodules from ${}_\beta R$ to ${}_\alpha R$. Then we have:
$$T(F(P))\simeq  M{}_\alpha R\oplus N{}_\alpha R=\sigma(P){}_\alpha R=\sigma(P) R$$
and we get the isomorphism $e$.

In order to construct the morphism $\tau$, we need to give an explicit description of $R$ as a $C$-bimodule.

We set:
$${}_\alpha U_\alpha=A'\hskip 24pt {}_\beta U_\beta=B'\hskip 24pt {}_\alpha U_\beta={}_\beta U_\alpha=0$$
in case 2 and:
$${}_\alpha U_\alpha=A'\hskip 24pt {}_\beta U_\beta=tA''t^{-1}\hskip 24pt {}_\alpha U_\beta=At^{-1}\hskip 24pt {}_\beta U_\alpha=tA$$
in case 3.

For each sequence $\Sigma=(i_1,j_1,i_2,j_2,\dots,i_n,j_n)$ in the set $\{\alpha,\beta\}$, the ring structure on $R$ induce a well defined morphism of $C$-bimodules:
$$\Phi(\Sigma): {}_{i_1}U_{j_1} {}_{i_2}U_{j_2} \dots {}_{i_n}U_{j_n}\fl R$$
and the sum of these morphisms is an epimorphism. On the other hand it is easy to see that, for each $(i,j,k)$ in $\{\alpha,\beta\}$ the image of $\Phi(i,j,j,k)$ is
contained in $C\oplus {}_i U_k$. Therefore, if there is some $k<n$ such that $j_k=i_{k+1}$ in the sequence $\Sigma=(i_1,j_1,i_2,j_2,\dots,i_n,j_n)$, every element of the
image of $\Phi(\Sigma)$ is reducible. As a consequence the sum of the morphisms $\Phi(\Sigma)$, for each sequence $\Sigma=(i_1,j_1,i_2,j_2,\dots,i_n,j_n)$ such that
$j_k\not=i_{k+1}$ for every $k<n$ is still an epimorphism.

Actually this sum is an isomorphism and we have a description of $R$ as $C$-bimodule:
$$R=C\oplus \bigoplus\Bigl({}_{i_1}U_{j_1}{}_{i_2}U_{j_2}\dots {}_{i_n}U_{j_n}\Bigr)$$
the sum being taken over all non empty sequences $(i_1,j_1,i_2,j_2, \dots,i_n,j_n)$ in $\{\alpha,\beta\}$ such that $j_k\not=i_{k+1}$ for all $k<n$.

This fact was proven in [Wa1], p. 140 (or in [C]) for the case 2 and in [Wa1], p. 150 (with a suggestion of S. Cappell) for the case 3. 

Denote by $f$ (resp. $g$) the unique bijection from $\{1,2\}$ to $\{\alpha,\beta\}$ such that $f(1)=\beta$ (resp. $g(1)=\alpha$). Then, because of the definition of $S$,
we have in case 2:
$$\forall i,j\in\{1,2\},\ \ {}_i S_j={}_{f(i)} U_{g(j)}$$

In case 3, we check that the multiplication by $t$ or $1$ on the left and by $t^{-1}$ or $1$ on the right induce for each $i,j$ in $\{1,2\}$ an isomorphism of
$C$-bimodules ${}_i S_j\build\fl_{}^\sim {}_{f(i)} U_{g(j)}$.

Then in cases 2 and 3 we have an isomorphism of $C$-bimodules:
$$R\simeq C\oplus \bigoplus\Bigl({}_{i_1}S_{j_1}{}_{i_2}S_{j_2}\dots {}_{i_n}S_{j_n}\Bigr)$$
the sum being taken over all non empty sequences $(i_1,j_1,i_2,j_2, \dots,i_n,j_n)$ in $\{1,2\}$ such that $j_k=i_{k+1}$ for all $k<n$. Hence we get an isomorphism of
$C$-bimodules:
$$R\simeq C\oplus\build\oplus_{i,j}^{}({}_i S_j\oplus{}_i (S^2)_j\oplus{}_i(S^3)_j\oplus\dots)\simeq C\oplus\build\oplus_{i,j,n>0}^{}{}_i(S^n)_j$$

So we are able to define the morphism $\tau$. If $M$ is a right $C$-module we have:
$$MR\simeq M\oplus\build\oplus_{i,j,n>0}^{} M{}_i(S^n)_j\simeq M\oplus\build\oplus_{j,n>0}^{} (s(M)S^n)_j\simeq M\oplus \build\oplus_{n>0}^{}\sigma(s(M)S^n)$$
But we have:
$$\sigma(s(M)[S])=\sigma s(M)\oplus \build\oplus_{n>0}^{}\sigma(s(M)S^n)$$
Then we define the morphism $\tau$ to be the identity on the direct sum of the $\sigma(s(M)S^n)$ and the morphism $\tau_0:\sigma s(M)\rightarrow M$ induced by the
adjunction on the first term $\sigma s(M)$.

The last thing to do is to construct the isomorphism $\e$ in cases 2 and 3.

In case 2 with $E=(E_A,E_B)$, we have:
$$T(E)=E_A\otimes_A R\oplus E_B\otimes_B R\simeq E_A(C\oplus B'\oplus B'A'\oplus B'A'B'\oplus\dots)\oplus E_B(C\oplus A'\oplus A'B'\oplus A'B'A'\oplus\dots)$$
$$\simeq E_A\oplus E_B\oplus E_A\bigl(\build\oplus_{j,n>0}^{} {}_1(S^n)_j\bigr)\oplus E_B\bigl(\build\oplus_{j,n>0}^{} {}_2(S^n)_j\bigr)$$
$$\simeq\sigma(\widehat F(E))\oplus\bigl(\build\oplus_{n>0}^{} \sigma(\widehat F(E)S^n))\simeq \sigma(\widehat F(E)[S])$$
which give the isomorphism $\e$.

In the last case we have:
$$R\simeq C\oplus\build\bigoplus_{i,j,n>0}^{}{}_i(S^n)_j\simeq(C\oplus {}_2 S_1)(C\oplus\build\bigoplus_{j,n>0}^{}{}_1(S^n)_j)
\oplus{}_2 S_2(C\oplus\build\bigoplus_{j,n>0}^{}{}_2(S^n)_j)$$
$$\simeq A(C\oplus\build\bigoplus_{j,n>0}^{}{}_1(S^n)_j)\oplus At^{-1}(C\oplus\build\bigoplus_{j,n>0}^{}{}_2(S^n)_j)$$
$$\simeq A_\alpha(C\oplus\build\bigoplus_{j,n>0}^{}{}_1(S^n)_j)\oplus A_\beta(C\oplus\build\bigoplus_{j,n>0}^{}{}_2(S^n)_j$$
and we check that the isomorphism from $R$ to this last module is an isomorphism of $(A,C)$-bimodules. Since $E$ belongs to $\A\v$, we have;
$$T(E)\simeq E_\alpha(C\oplus\build\bigoplus_{j,n>0}^{}{}_1(S^n)_j)\oplus E_\beta(C\oplus\build\bigoplus_{j,n>0}^{}{}_2(S^n)_j)$$
$$\simeq\sigma(\widehat F(E))\oplus\build\bigoplus_{n>0}^{} \sigma(\widehat F(E)S^n)\simeq\sigma(\widehat F(E)[S])$$
which give the isomorphism $\e$ in this last case.\cqfd
\vskip 12pt
\noi{\bf 1.9 Remark:} We have an explicit description of the morphism $\tau$ and the isomorphism $\e$.

Consider two modules $M\in\C\v$ and $E\in\D\v$.

In case 1, set: $\overline u=u$ for each $u\in E$.

In case 2, the maps $E_A\fl E_AR$ and $E_B\fl E_BR$ induce a map $\sigma(E)\fl T(E)$ and in case 3, the maps $E_\alpha\fl ER$ and $E_\beta\fl EtR\fl ER$
induce also a map $\sigma(E)\fl T(E)$. Denote by $v\mapsto\overline v$ this map.

So we have a morphism $v\mapsto\overline v$ from $\sigma(E)$ to $T(E)$ in the three cases.

Denote also by $s\mapsto\overline s$ the isomorphism ${}_i S_j\build\fl_{}^\sim{}_{f(i)}U_{g(j)}\subset R$ in cases 2 or 3 and the identity $S\fl S$ in case 1.

With these notations, we have the following description of $\tau$ and $\e$:

In case 1, for each integer $n\geq0$, each $u\in M$, each $v\in E$ and each sequence $(s_1,s_2,\dots, s_n)$ in $S$ we have:
$$\tau(us_1s_2\dots s_n)=u\ \overline{s_1}\ \overline{s_2}\dots \overline{s_n}\in MR$$
$$\e(vs_1s_2\dots s_n)=\overline v\ \overline{s_1}\ \overline{s_2}\dots \overline{s_n}\in T(E)$$

In case 2 or 3, for each integer $n\geq0$, each sequence $(i_0,i_1,\dots,i_n)$ in $\{1,2\}$, each $u\in M$, each $v\in E_{i_0}$ and each sequence $(s_1,s_2,\dots s_n)$,
with $s_k\in{}_{i_{k-1}}S_{i_k}$, we have:
$$\tau(us_1s_2\dots s_n)=u\ \overline{s_1}\ \overline{s_2}\dots \overline{s_n}\in MR$$
$$\e(vs_1s_2\dots s_n)=\overline v\ \overline{s_1}\ \overline{s_2}\dots \overline{s_n}\in T(E)$$
\vskip 12pt
\noi{\bf 1.10 Lemma:} {\sl There is a natural transformation:
  $$\Lambda:\F(E,M)\fl\Hom(\widehat F(E),s(M)\oplus s(M)S)$$
  for each $(E,M)\in\times\D\v\times\C\v$ with the following properties:

  $\bullet$ $\Lambda$ is injective in case 2 or 3 and bijective in case 1

  $\bullet$ for each $(P,M)\in\C_2\v\times\C\v$, if $f$ is the map $P\fl \widehat F(F(P))$ induced by adjunction, the composite morphism:
  $$\F(F(P),M)\build\fl_{}^\Lambda \Hom(\widehat F(F(P)),s(M)\oplus s(M)S)\build\fl_{}^{f^*} \Hom(P,s(M)\oplus s(M)S)$$
is bijective

$\bullet$ for each $\varphi\in\F (E,M)$, we have a commutative diagram:
$$\diagram{\sigma(\widehat F(E)[S])&\hfl{g}{}&\sigma(s(M)[S])\cr\vfl{\e}{}&&\vfl{}{\tau}\cr T(E)&\hfl{\varphi}{}&MR\cr}$$
where $g$ is the composite morphism:
$$\sigma(\widehat F(E)[S])\build\fl_{}^{\Lambda(\varphi)}\sigma((s(M)\oplus s(M)S)[S])\build\fl_{}^h \sigma(s(M)[S])$$
and $h$ the morphism induced by the identity $s(M)[S]\fl s(M)[S]$ and the stabilization map $s(M)S[S]\fl s(M)[S]$.}

\vskip 12pt
\noi{\bf Proof:} Consider the case 1. We have an isomorphism from Hom$_\R(ER,MR)$ to Hom$_\C(E,MR)$ inducing an isomorphism ${\rde F}(E,M)\simeq \Hom(E,M\oplus MS)$
and we get the isomorphism $\Lambda:\F(E,M)\build\fl_{}^\sim\Hom(E,M\oplus MS)=\Hom(\widehat F(E),s(M)\oplus s(M)S)$.

Consider the case 2. We have $E=(E_A,E_B)\in\A\times\B$ and we get isomorphisms:
$$\Hom_\R(T(E),MR)\simeq \Hom_\R(E_A R\oplus E_B R,MR)\simeq \Hom_\A(E_A,MR)\oplus\Hom_\B(E_B,MR)$$
and then isomorphisms:
$${\rde F}(E,M)\simeq \Hom(E_A,MA)\oplus\Hom(E_B,MB)\simeq\Hom(E,(MA,MB))=\Hom(E,F(s(M)))$$

Consider now the last case. We have:
$$\Hom_\R(T(E),MR)\simeq\Hom_\R(ER,MR)\simeq\Hom_\A(E,MR)$$
and then:
$$\F(E,M)\simeq \Hom(E,MA\oplus MtA)\simeq\Hom(E,M{}_\alpha A\oplus M{}_\beta A)\simeq\Hom(E,F(s(M)))$$

Therefore, in case 2 and 3, we have an isomorphism:
$$\F(E,M)\build\fl_{}^\sim\Hom(E,F(s(M)))\leqno{(*)}$$

On the other hand, the morphism $\widehat F$ induces an injection:
$$\Hom(E,F(s(M)))\fl\Hom(\widehat F(E),\widehat F(F(s(M))))\simeq\Hom(\widehat F(E),s(M)\oplus s(M)S)$$
(see lemma 1.6) and we get the desired injection
$$\Lambda:\F(E,M)\fl\Hom(\widehat F(E),s(M)\oplus s(M)S)$$

Let $(P,M)$ be an object in $C_2\v\times\C\v$. In case 1, the morphism $P\fl \widehat F(F(P))$ is the identity and the composite map:
$$\F(F(P),M)\build\fl_{}^\Lambda \Hom(\widehat F(F(P)),s(M)\oplus s(M)S)\build\fl_{}^{f^*} \Hom(P,s(M)\oplus s(M)S)$$
is bijective.

Consider the other cases. Let $\varphi\in\F(F(P),M)$ be an admissible morphism and $\alpha:F(P)\fl F(s(M))$ be the corresponding morphism (via the isomorphism $(*)$).
Let $f:P\fl \widehat F F(P)$ be the morphism adjoint to the identity of $F(P)$. By
adjunction, the composite morphism $P\build\fl_{}^f\widehat F F(P))\build\fl_{}^\alpha\widehat F F(s(M))$ is the morphism obtained from $\alpha$ by adjunction and we
have a bijection $\Hom(F(P),F(s(M)))\simeq\Hom(P,\widehat F F(s(M))$. Hence we have a bijection:
$$\F(F(P),M)\simeq\Hom(P,\widehat F F(s(M))\simeq\Hom(P,s(M)\oplus s(M)S)$$
which is nothing else but the map $f^*\circ\Lambda$.
\vskip 12pt
Denote by $(D)$ the diagram of the lemma.

In case 1, $\e$ and $\tau$ are identities and we have: $g=\varphi$. Hence $(D)$ is commutative.

Consider the other cases. Via the bijection $\F(E,M)\simeq\Hom(E,F(s(M)))$, the morphism $\varphi\in\F(E,M)$ corresponds to a morphism $\widetilde\varphi:E\fl F(s(M))$
and we have a diagram:
$$\diagram{\sigma(\widehat F(E)[S])&\hfl{\widetilde\varphi}{}&\sigma(\widehat F F(s(M))[S])&\hfl{g_0}{}&\sigma(s(M)[S])\cr
\vfl{\e}{}&&\vfl{\e}{}&&\vfl{\tau}{}\cr T(E)&\hfl{\widetilde\varphi}{}&T(F(s(M)))&\hfl{\varphi_0}{}&MR\cr}$$

In this diagram, the square on the left is commutative by naturality and the square on the right $(D_0)$ is the diagram $(D)$ in the case: $E=F(s(M))$.
Moreover the total square is the diagram $(D)$. Hence, to prove the commutativity of $(D)$ it is enough to prove that $(D_0)$ is commutative.

In $(D_0)$, the morphism $g_0$ is induced by the isomorphism $\widehat F F(s(M))\build\fl_{}^\sim s(M)\oplus s(M)S$, the identity $s(M)[S]\fl s(M)[S]$ and the
stabilization map $s(M)S[S]\fl s(M)[S]$.

The morphism $\varphi_0$ is the composite map:
$$T(F(s(M)))\build\fl_{}^\sim (MA\oplus MB)R\simeq MR\oplus MR\build\fl_{}^+ MR$$
in case 2 and the composite map:
$$T(F(s(M)))\build\fl_{}^\sim (MA\oplus MtA)R\simeq MR\oplus MR\build\fl_{}^+ MR$$
in case 3. Hence $\varphi_0$ is the composite map:
$$T(F(s(M)))\build\fl_{\sim}^e\sigma s(M)R\build\fl_{}^a MR$$
where $a:\sigma s(M)\fl M$ is adjoint to the identity of $s(M)$.

Consider an element $u\in\widehat F F(s(M))$, an integer $n\geq0$, a sequence $(i_0,i_1,\dots,i_n)$ in $\{1,2\}$ and a sequence $(s_1,s_2,\dots,s_n)$ with
$s_k\in{}_{i_{k-1}}S_{i_k}$. Denote by $v$ the image of $u$ under the isomorphism $\widehat F F(s(M))\simeq s(M)\oplus s(M)S$.
If $v$ belongs to $s(M)$ or to $s(M)S_{i_0}$, we have, because of remark 1.9,  the following:
$$\varphi_0\e(us_1s_2\dots s_n)=v\ \overline{s_1}\ \overline{s_2}\dots \overline{s_n}=\tau g_0(us_1s_2\dots s_n)$$
and $(D_0)$ is therefore commutative.\cqfd
\vskip 12pt
\noi{\bf 1.11 The functors $\Phi:\Nil(C_2,S)\v\fl \V\v$ and $\Psi:\S\v\fl \Nil(C_2,S)\v$.}
\vskip 12pt
Let $H$ be a module in $\C_2\v$. Because of lemma 1.10, we have an isomorphism:
$$\zeta:{\rde F}(F(H),\sigma(H))\build\fl_{}^\sim \Hom(H,s\sigma(H)\oplus s\sigma(H)S)$$
Therefore, for every morphism $\theta:H\fl HS$, we have a unique morphism $\varphi_\theta$ in ${\rde F}(F(H),\sigma(H)$ such that $\zeta(\varphi_\theta)$ is the
composite morphism:
$$H\build\fl_{}^{1-\theta} H\oplus HS\build\fl_{}^i s\sigma H\oplus s\sigma H S$$
where $i:H\fl s\sigma H$ is adjoint to the identity. Thus $\Phi(H,\theta)=(F(H),\sigma(H),\varphi_\theta)$ is a well defined split diagram in $\S\v$. 
\vskip 12pt
\noi{\bf 1.12 Lemma:} {\sl The correspondence above induces two equivalences of categories $\Phi:\Nil(C_2,S)\build\fl_{}^\sim\V$ and 
$\Phi:\Nil(C_2,S)\v\build\fl_{}^\sim\V\v$.

Moreover there is a functor $\Psi:\S\v\fl\Nil(C_2,S)\v$ and a morphism of functors $\pi$ from $I\Phi\Psi$ to the identity, where $I:\V\v\fl\S\v$ is the 
inclusion, such that the following holds for every $X\in\S\v$ (with: $\Psi(X)=(H,\theta)$):

$\bullet$ we have a natural exact sequence in $\C\v$:
$$0\fl\Phi_1(X)\fl \sigma(H)\build\fl_{}^\pi \Phi_3(X)$$

$\bullet$ a splitting diagram $X$ belongs to $\overline{\Mv}$ (resp. to $\V\v$) if and only if the morphism $\pi:\Phi\Psi(X)\fl X$ is an epimorphism (resp. an 
isomorphism).}
\vskip 12pt
\noi{\bf Remark:} Actually, everything works without any flatness condition from subsection 1.4 to 1.11. But this condition is strongly needed for lemma 1.12,
essentially for constructing the functor $\Psi$.
\vskip 12pt
\noi{\bf Proof of lemma 1.12:} Let $(H,\theta)$ be an object of $\Nil(C_2,S)\v$. If $\theta=0$, the morphism $\varphi_0:T(F(H)\fl \sigma(H)R$ is nothing else but
the isomorphism $e$ defined in lemma 1.8. Therefore $\Phi(H,\theta_\theta)$ belongs to $\V\v$ in this case.

Suppose $\theta$ is nilpotent. Then, because of lemma 1.2, there is a filtration $0=H_0\subset H_1\subset H_2\subset\dots$ of $H$ such that $H$ is the union of the $H_i$'s
and, for all $i>0$, $\theta(H_i)$ is contained in $H_{i-1}S$. 

On the other hand it is easy to see that the functor $\Phi$ is exact. So we get a filtration:
$$0=\Phi(H_0,\theta)\subset\Phi(H_1,\theta)\subset\Phi(H_2,\theta)\subset\dots \subset\Phi(H,\theta)$$
Because each $\Phi(H_i/H_{i-1},\theta)=\Phi(H_i/H_{i-1},0)$ belongs to $\V\v$, each $\Phi(H_i,\theta)$ is also in $\V\v$ and then $\Phi(H,\theta)$ is a splitting module.

Therefore $\Phi$ is a functor from $\Nil(C_2,S)\v$ to $\V\v$ and also from $\Nil(C_2,S)$ to $\V$.

Let $X=(E,M,\varphi)$ be a splitting diagram with $E\in\D\v$ and $M\in\C\v$. Denote by $P$ the module $\widehat F(E)\in\C_2\v$. By composing the morphism
$\varphi\e:\sigma(P[S])\fl MR$ with the identity $P[S]=\sigma(P[S])$, we get a $\Z$-linear map $f:P[S]\fl MR$.
Denote by $H$ the $\Z$-submodule $f^{-1}(M)$ and by $i$ the inclusion map $H\fl P[S]$.

Because of lemma 1.10, $\varphi$ is determined by a morphism $\psi=\Lambda(\varphi)$ from $P=\widehat F(E)$ to $s(M)\oplus s(M)S$. Then we have two morphisms
$\lambda:P\fl s(M)$ and $\gamma:P\fl s(M)S$ such that: $\psi=\lambda-\gamma$.
Moreover the composite map $P[S]=\sigma(P[S])\build\fl_{}^{\varphi\e} MR$ is equal to : $\tau(\lambda-\gamma)$.

We have:
$$u\in H\ \ \Longleftrightarrow\ \ \varphi\e(u)\in M\ \ \Longleftrightarrow\ \ \tau(\lambda-\gamma)(u)\in M$$
$$\Longleftrightarrow\ \ \forall k\geq0,\ \ \tau(\lambda(u_{k+1})-\gamma u_k))=0$$

But the morphism $\tau:\sigma(s(M)S^k)\fl MR$ is injective for all $k>0$ (see lemma 1.8). Then we have:
$$u\in H\ \ \Longleftrightarrow\ \ \forall k\geq0,\ \ \lambda(u_{k+1})=\gamma(u_k)$$

Since $\lambda$ and $\gamma$ are morphisms in $\C_2\v$, $H$ is a $C_2$-module and we have an exact sequence in $\C_2\v$:
$$0\fl H\build\fl_{}^i P[S]\build\fl_{}^\delta N[S]$$
where $N$ is the module $s(M)S$ and $\delta$ is the morphism sending $u=u_0+u_1+u_2+\dots$ (with $u_k\in PS^k$ for every $k\geq0$) to:
$$\delta(u_0+u_1+u_2+\dots)=\build\sum_{k\geq0}^{}(\lambda(u_{k+1})-\gamma(u_k))$$

If $U$ is a module in $\C_2\v$, we have a morphism $\theta_U:U[S]\fl U[S]S$ sending $u_0+u_1+u_2+\dots\in U[S]$ (with $u_k\in PS^k$ for every $k\geq0$) to:
$$\theta_U(u_0+u_1+u_2+\dots)=u_1+u_2+u_3+\dots$$
Notice that $\theta_U^k$ sends $u_0+u_1+u_2+\dots\in U[S]$ to $u_k+u_{k+1}+\dots$ and $\theta_U$ is nilpotent. Hence $(U[S],\theta_U)$ belongs to $\Nil(C_2,S)\v$.

Consider the following diagram:
$$\diagram{P[S]&\hfl{\delta}{}&N[S]\cr\vfl{\theta_P}{}&&\vfl{\theta_N}{}\cr P[S]S&\hfl{\delta}{}&N[S]S\cr}\leqno{(D)}$$
Let $n\geq 0$ be an integer, $u$ be an element in $P$ and $s_1,s_2,\dots,s_n$ be elements in $S$. We have the following:

$$n=0\ \ \Longrightarrow\ \ \delta\theta_P(u)=\theta_N\delta(u)=0$$
$$n=1\ \ \Longrightarrow\ \ \delta\theta_P(us_1)=\theta_N\delta(us_1)=-\gamma(u)s_1$$
$$n>1\ \ \Longrightarrow\ \ \delta\theta_P(us_1\dots s_n)=\theta_N\delta(us_1\dots s_n)=\lambda(u)s_1\dots s_n-\gamma(u)s_1\dots s_n$$
and the diagram $(D)$ is commutative.

Since $S$ is flat on the left, $\Nil(C_2,S)\v$ is an abelian category and there is a unique nilpotent morphism $\theta:H\fl HS$ such that the following sequence is exact
in $\Nil(C_2,S)\v$:
$$0\fl (H,\theta)\build\fl_{}^i (P[S],\theta_P)\build\fl_{}^\delta (N[S],\theta_N)$$

Therefore $\Psi(X)=(H,\theta)$ is a well defined object in $\Nil(C_2,S)\v$ and we get the desired functor $\Psi:\S\v\fl\Nil(C_2,S)\v$.
\vskip 12pt
Consider the splitting diagram $\Phi\Psi(X)=(F(H),\sigma(H),\varphi_\theta)$. The morphism $\varphi_\theta$ corresponds to the composite morphism
$H\build\fl_{}^{1-\theta}H\oplus HS\fl s\sigma(H)\oplus s\sigma(H)S$. We have to construct a morphism $\pi:\Phi\Psi(X)\fl X$ in $\S\v$. This morphism is
given by two morphisms $\pi_0:F(H)\fl E$ and $\pi_1:\sigma(H)\fl M$.

For each $u\in H$ we set: $\pi'_0(u)=u_0$, with $i(u)=u_0+u_1+\dots$ and $u_k\in PS^k$ for all $k$. So we have two morphisms $\pi'_0:H\fl\widehat F(E)$ and 
$\lambda\pi'_0:H\fl s(M)$ and, by adjunction, two morphisms $\pi_0:F(H)\fl E$ and $\pi_1:\sigma(H)\fl M$.

For $u\in H$, with: $i(u)=u_0+u_1+\dots$, we have:
$$\psi\pi'_0(u)=\psi(u_0)=\lambda(u_0)-\gamma(u_0)=\lambda(u_0)-\lambda(u_1)=\lambda\pi'_0(u))-\lambda\pi'_0\theta(u)=\lambda\pi'_0(1-\theta)(u)$$
and the following diagram is commutative:
$$\diagram{H&\hfl{\psi_0}{}&s\sigma(H)\oplus s\sigma(H)S\cr\vfl{\pi'_0}{}&&\vfl{\pi_1}{}\cr P&\hfl{\psi}{}&s(M)\oplus s(M)S\cr}$$
Therefore the two morphisms $\pi_0:F(H)\fl E$ and $\pi_1:\sigma(H)\fl M$ induce a well defined morphism $\pi:\Phi\Psi(X)\fl X$ and we get a commutative diagram:
$$\diagram{T(F(H))&\hfl{\varphi_\theta}{\sim}&\sigma(H) R\cr \vfl{\pi_0}{}&&\vfl{\pi_1}{}\cr T(E)&\hfl{\varphi}{}&MR\cr}\leqno{(D_X)}$$

We have an exact sequence in $\R\v$:
$$0\fl \Phi_1(X)\build\fl_{}^{\mu} T(E)\build\fl_{}^\varphi MR$$
and then an exact sequence in $\C\v$:
$$0\fl \Phi_1(X)\build\fl_{}^{\e^{-1}\mu} \sigma(P[S])\build\fl_{}^{\varphi\e} MR$$
Therefore we have a commutative diagram in $\C\v$ with exact lines:
$$\diagram{0&\hfl{}{}&\Phi_1(X)&\hfl{}{}&\sigma(H)&\hfl{\pi_1}{}&M\cr
&&\vfl{=}{}&&\vfl{\e i}{}&&\vfl{j}{}\cr
0&\hfl{}{}&\Phi_1(X)&\hfl{\mu}{}&T(E)&\hfl{\varphi}{}&MR\cr}$$
where $j:M\fl MR$ is the inclusion, and the top line of this diagram is the desired exact sequence.
\vskip 12pt
We have the following equivalences:
$$X\in\overline{\Mv}\ \ \Longleftrightarrow\ \ \hbox{the morphism}\ \ \varphi:T(E)\fl MR\ \ \hbox{is surjective}$$
$$\Longleftrightarrow\ \ \hbox{the image of}\ \ \varphi\ \ \hbox{contains}\ \ M$$
$$\Longleftrightarrow\ \ \hbox{the image of}\ \ \sigma(P[S])\fl MR\ \ \hbox{contains}\ \ M$$
$$\Longleftrightarrow\ \ \sigma(H)\build\fl_{}^{\pi'_0}\sigma(P)\build\fl_{}^{\lambda'}M\ \ \hbox{is surjective}$$
where $\lambda'$ is adjoint to $\lambda$. Therefore $X$ belongs to $\overline{\Mv}$ if and only if the morphism $\pi_1:\sigma(H)\fl M$ is surjective.

Suppose $X$ belongs to $\overline{\Mv}$. Let $E'$ be the image of $\pi_0:F(H)\fl E$. Because of the diagram above, $T(E')$ contains the image of $\mu:\Phi_1(X)\fl T(E)$
and the image of $\e i:\sigma(H)\fl T(E)$. But $T(E')$ is a $R$-submodule of $T(E)$ and $\varphi(T(E'))$ is a $R$-submodule of $MR$ containing $M$. Therefore
$\varphi:T(E')\fl MR$ is surjective and $T(E')$ contains the kernel $\Phi_1(X)$ of $\varphi$. Hence we have: $T(E')=T(E)$ and then: $E'=E$.

Consequently, if $X$ belongs to $\overline{\Mv}$, $\pi_1:\sigma(H)\fl M$ and $\pi_0:F(H)\fl E$ are surjective and $\pi:\Phi\Psi(X)\fl X$ is an epimorphism (in $\S\v$).
Conversely, if $\pi:\Phi\Psi(X)\fl X$ is an epimorphism, $\pi_1:\sigma(H)\fl M$ is surjective and $X$ belongs to $\overline{\Mv}$.
\vskip 12pt
Suppose $X$ is in $\V\v$. Then $\pi_0$ and $\pi_1$ are epimorphisms. Because of the exact sequence:
$$0\fl\Phi_1(X)\fl\sigma(H)\build\fl_{}^{\pi_1} M$$
the morphism $\pi_1:\sigma(H)\fl M$ is an isomorphism. Therefore in the diagram $(D_X)$, $\pi_0$ and
$\varphi$ are isomorphisms and $\pi:\Phi\Psi(X)\fl X$ is an isomorphism too. As a consequence, the functor $\Phi$ induces two equivalences of categories
$\Nil(C_2,S)\fl \V$ and $\Nil(C_2,S)\v\fl \V\v$.

Conversely, if $\pi$ is an isomorphism, $\pi_0$ and $\pi_1$ are isomorphisms and $\varphi$ is an isomorphism too. Therefore $X$ belongs to $\V\v$.\cqfd
\vskip 12pt
\noi{\bf 1.13 Lemma:} {\sl Let $(A,S)$ be a left-flat bimodule and $X=(M,\theta)$ be an object in $\Nil(A,S)\v$. Let $V$ be a finitely generated
  projective right $A$-module and $f: V\fl M$ be a morphism. Then there exist an object $Y=(M',\theta')\in\Nil(A,S)$, a split injective morphism $f':V\fl M'$ and a morphism
  $g:Y\fl X$ making the following diagram commutative:
  $$\diagram{V&\hfl{f'}{}&M'\cr \vfl{=}{}&&\vfl{g}{}\cr V&\hfl{f}{}&M\cr}$$}
\vskip 12pt
\noi{\bf Proof:} Denote by $\A\v$ the category of right $A$-modules and $\A\subset\A\v$ the category of finitely generated projective modules in $\A\v$. Because of lemme
1.2 there is a filtration:
$$0=M_0\subset M_1\subset M_2\subset\dots$$
of $M$ by $A$-modules such that $M$ is the union of the $M_i$'s and $\theta(M_i)\subset M_{i-1}S$ for every $i>0$. Since $V$ is finitely generated, there is an integer
$n>0$ such that $M_n$ contains the image of $f:V\fl M$.

Then we'll construct modules $F_i\in\A$, morphisms $h_i:F_i\fl M_i$ for $i=0,1,\dots,n$ and morphisms $\theta':F_i\fl F_{i-1} S$ for $i=1,2,\dots,n$ such that:

$\bullet$ $F_n=V$, $F_0=0$ and $h_n=f$

$\bullet$ for each $i=1,2,\dots,n$, the following diagram is commutative:
$$\diagram{F_i&\hfl{h_i}{}&M_i\cr \vfl{\theta'}{}&&\vfl{\theta}{}\cr F_{i-1}S&\hfl{h_{i-1}}{}&M_{i-1}S\cr}\leqno{(D_i)}$$

Let $p$ be an integer with $0\leq p\leq n$. Denote by $H(p)$ the following property:

$\bullet$ There exist modules $F_i\in\A$ and morphisms $h_i:F_i\fl E_i$ for $i=p,\dots, n$ and also morphisms $\theta':F_i\fl F_{i-1}S$ for $i=p+1,\dots,n$ such that
$F_n=V$, $h_n=f$, the diagram $(D_i)$ is commutative for $i=p+1,\dots,n$ and $F_0=0$ (if $p=0$).

This property is clearly true if $p=n$. Suppose $H(p)$ is true with $p>0$ and consider the composite morphism $\lambda:F_p\fl M_p\build\fl_{}^\theta M_{p-1}S$. If $p=1$,
this morphism is trivial and we set: $F_0=0$. Therefore the property $H(p-1)=H(0)$ is true.

Consider the case $p>1$. Since $F_p$ is finitely generated, $M_{p-1}$ contains a finitely generated submodule $M'$ such that: $\lambda(F_p)\subset M'S$. Let $F_{p-1}$ be
a module in $\A$ and $\mu:F_{p-1}\fl M'$ be an epimorphism. Since $F_p$ is projective the morphism $F_p\fl M'S$ factorizes through $F_{p-1}S$ and we have a commutative
diagram:
$$\diagram{F_p&\hfl{\theta'}{}&F_{p-1}S\cr\vfl{=}{}&&\vfl{\mu}{}\cr F_p&\hfl{\lambda}{}&M'S\cr}$$
So we define the morphism $h_{p-1}$ as the composite map: $F_{p-1}\build\fl_{}^\mu M'\subset M_{p-1}$ and we have the property $H(p-1)$.

By induction we obtain the property $H(0)$ and all the data are constructed.

Then we set: $M'=F_0\oplus F_1\oplus\dots\oplus F_n$. The morphisms $h_i$ induce a morphism $g:M'\fl M$ and the morphisms $\theta':F_i\fl F_{i-1}S$ induce a morphism
$\theta':M'\fl M'S$. The lemma is now easy to check.\cqfd
\vskip 12pt
Denote by $\R'$ the full subcategory of $\R$ generated by modules on the form $U=VR$ with $V\in\C$. This category is exact and cofinal in $\R$, that is, for each module
$M\in\R$, there is a module $M'\in\R$ such that $M\oplus M'$ belongs to $\R'$.
\vskip 12pt
\noi{\bf 1.14 Lemma:} {\sl Let $X$ be a split diagram in $\S\v$, $V$ be a module in $\R'$ and $f:V\fl \Phi_1(X)$ be a morphism in $\R\v$. Then there exist an object
  $Y\in\Mv$, a morphism $g:Y\fl X$ in $\S\v$ and an isomorphism $\e:V\build\fl_{}^\sim \Phi_1(Y)$ such that the following diagram is commutative:
  $$\diagram{V&\hfl{\e}{\sim}&\Phi_1(Y)\cr\vfl{=}{}&&\vfl{g}{}\cr V&\hfl{f}{}&\Phi_1(X)\cr}$$

  Moreover, if $X$ belongs to $\Mv$, the morphism $g$ can be chosen to be an epimorphism in $\S\v$.}
\vskip 12pt
\noi{\bf Proof:} The split diagram $X$ is a triple $(E,M,\varphi)$ with $E\in\D\v$ and $M\in\C\v$. Since $V$ belongs to $\R'$, there is a module $W\in\C$ such that:
$V=WR$ and we get a morphism $f':W\fl \Phi_1(X)$ in $\C\v$.

Denote by $K=(H,\theta)$ the object $\Psi(X)\in\Nil(C_2,S)\v$. Because of lemma 1.12, we have an exact sequence:
$$0\fl \Phi_1(X)\build\fl_{}^\mu\sigma(H)\build\fl_{}^\pi M$$
and the morphism $\mu f':W\fl\sigma(H)$ has an adjoint $\lambda:s(W)\fl H$. Because of lemma 1.13, there are an object $K'=(H',\theta')\in\Nil(C_2,S)$, a morphism
$h:K'\fl K$ and a split injective morphism $\lambda':s(W)\fl H'$ making the following diagram commutative:
$$\diagram{s(W)&\hfl{\lambda'}{}&H'\cr\vfl{=}{}&&\vfl{h}{}&\cr s(W)&\hfl{\lambda}{}&H\cr}$$

The morphism $\widetilde\lambda': W\fl\sigma(H')$ adjoint to $\lambda'$ is still split injective and we get a commutative diagram with exact lines:
$$\diagram{0&\hfl{}{}&W&\hfl{\widetilde\lambda'}{}&\sigma(H')&\hfl{\pi'}{}&M'&\hfl{}{}&0\cr &&\vfl{f'}{}&&\vfl{h}{}&&\vfl{}{}&&\cr
0&\hfl{}{}&\Phi_1(X)&\hfl{\mu}{}&\sigma(H)&\hfl{\pi}{}&M&&\cr}$$
where $M'$ is the cokernel of $\widetilde\lambda'$. 

We have now an object $Y'=\Phi(K')=(F(H'),\sigma(H'),\varphi_{\theta"})\in\V$, an object $Y=(F(H'),M',\pi'\varphi_{\theta'})\in\Mv$ and an epimorphism $u:Y'\fl Y$. But
the morphism $\Phi(K')\fl\Phi(K)\build\fl_{}^\pi X$ vanishes on the kernel of $u$ and factorizes by a morphism $g:Y\fl X$. So we have morphisms
$Y'\fl Y\fl X$ inducing a commutative diagram with exact lines:
$$\diagram{0&\hfl{}{}&0&\hfl{}{}&T(F(H'))&\hfl{}{}&\sigma(H')R&\hfl{}{}&0\cr&&\vfl{}{}&&\vfl{}{}&&\vfl{\pi'}{}&&\cr
0&\hfl{}{}&WR&\hfl{}{}&T(F(H'))&\hfl{}{}&M'R&\hfl{}{}&0\cr&&\vfl{f}{}&&\vfl{}{}&&\vfl{}{}&&\cr
0&\hfl{}{}&\Phi_1(X)&\hfl{}{}&T(E)&\hfl{}{}&MR&&\cr}$$

Hence, we get an isomorphism $\e:V=WR\fl \Phi_1(Y)$ and a commutative diagram:
$$\diagram{V&\hfl{\e}{\sim}&\Phi_1(Y)\cr\vfl{=}{}&&\vfl{g}{}\cr V&\hfl{f}{}&\Phi_1(X)\cr}$$

Suppose $X$ belongs to $\Mv$. Because of lemma 1.12, the morphism $\pi:\Phi\Psi(X)\fl X$ is an epimorphism in $\S\v$ inducing two epimorphisms $F(H)\fl E$ and
$\sigma(H)\fl M$.

Since $E$ is finitely generated, $H$ contains a finitely generated $C_2$-submodule $H_0\subset H$ such that the composite morphism $F(H_0)\fl F(H)\fl E$ is an
epimorphism. Therefore there exist a module $P\in\C$ and a morphism $u:s(P)\fl H$ such that the image of $u$ contains the submodule $H_0$. Because of lemma 1.13,
there are an object $K'=(H',\theta')\in\Nil(C_2,S)$, a morphism
$h:K'\fl K$ and a split injective morphism $\lambda'\oplus u':s(W)\oplus s(P)\fl H'$ making the following diagram commutative:
$$\diagram{s(W)\oplus s(P)&\hfl{\lambda'\oplus u'}{}&H'\cr\vfl{=}{}&&\vfl{h}{}&\cr s(W)\oplus s(P)&\hfl{\lambda\oplus u}{}&H\cr}$$
Therefore the composite morphism $F(H')\build\fl_{}^h F(H)\fl E$ is an epimorphism.

If we continue the construction above by using this morphism $h:K'\fl K$, we get a morphism $g:Y\fl X$ such that the morphism $\Phi_2(Y)\fl \Phi_2(X)$ is isomorphic to
the morphism $F(H')\fl E$ which is an epimorphism. Hence $\Phi_2(Y)\fl \Phi_2(X)$ is also an epimorphism. On the other hand, we have a commmutative diagram:
$$\diagram{T(\Phi_2(Y))&\hfl{}{}&\Phi_3(Y)R\cr\vfl{}{}&&\vfl{}{}\cr T(E)&\hfl{\varphi}{}&MR\cr}$$
where $\varphi:T(E)\fl MR$ is surjective. Then $\Phi_3(Y)R\fl MR$ is surjective and $\Phi_3(Y)\fl M$ is surjective too. The result follows.\cqfd
\vskip 24pt
\noi{\bf 2. Algebraic K-theory of categories $\V$ and $\Mv$.}
\vskip 12pt
\noi{\bf 2.1 About Waldhausen K-theory.}
\vskip 12pt
A Waldhausen category is a category with a zero object and two subcategories: the category of cofibrations and the category of (weak-)equivalences. These categories have
to satisfy certain conditions (see [Wa2]). Waldhausen associates to any essentially small Waldhausen category $\C$ an infinite loop space $K(\C)$ and $K$ (called the
Waldhausen K-theory functor) is a functor from the category of essentially small Waldhausen categories to the category of infinite loop spaces.

An exact category $\E$ may be considered as a Waldhausen category, where a cofibration is an admissible monomorphism of $\E$ (i.e. a morphism $f$ appearing is an exact
sequence $0\fl X\build\fl_{}^f Y\fl Z\fl 0$ in $\E$) and an equivalence is an isomorphism in $\E$. Moreover, if $\E$ is essentially small, we have a natural homotopy
equivalence from the Quillen K-theory of $\E$ to the Waldhausen K-theory of $\E$.
\vskip 12pt
To every exact category $\E$ we can associate the following category $\E_*$:

The objects of $\E_*$ (called the $\E$-complexes) are the complexes:
$$C=\Bigl(\dots\build\fl_{}^d C_n\build\fl_{}^d C_{n-1}\build\fl_{}^d C_{n-2}\build\fl_{}^d \dots\Bigr)$$
where each $C_n$ is an object of $\E$ and each $d$ is a morphism of $\E$ such that the sum $\build\oplus_n^{} C_n$ exists in $\E$ and each morphism $d^2$ is zero.

The morphisms in this category are morphisms respecting degrees and differentials. A sequence $0\fl X\fl Y\fl Z\fl 0$ in $\E_*$
is said to be exact if it induces an exact sequence in $\E$ on each degree. With these exact sequences, $\E_*$ becomes an exact category.

If $\E$ is the category of right modules (resp. the category of finitely generated projective right modules) over a ring $A$, the $\E$-complexes are called $A$-complexes
(resp. finite $A$-complexes).

Suppose $\E$ is an exact subcategory of an abelian category $\E\v$. Then $\E_*$ is a Waldhausen category, where cofibrations are admissible monomorphisms and
equivalences are morphisms inducing an isomorphism in homology (where homologies are computed in $\E\v$). Moreover $\E_*$ is saturated and has a cylinder functor
satisfying the cylinder axiom (in the sense of Waldhausen [Wa2]). We have the following result ([Wa2], [We]):
\vskip 12pt
\noi{\bf 2.2 Gillet-Waldhausen theorem:} {\sl Let $\E$ be an essentially small exact category contained in an abelian category $\E\v$. Suppose $\E$ is stable in $\E\v$ by
  kernel of epimorphisms. Then the inclusion $\E\subset\E_*$ of Waldhausen categories induces a homotopy equivalence in K-theory.}
\vskip 12pt
Let $\E$ be an essentially small exact category and $C$ and $C'$ be two $\E$-complexes. For each integer $n$, we set:
$$\Hom(C,C')_n=\prod_p\Hom_{\E}(C_p,C'_{n+p})$$
and $\Hom(C,C')$ is a graded $\Z$-module. We have on $\Hom(C,C')$ a natural differential $d$ of degree $-1$ defined by:
$$\forall f\in\Hom(C,C')_n,\ \ d(f)=d\circ f-(-1)^n f\circ d$$
An element of $\Hom(C,C')_n$ is called a linear map of degree $n$, a cycle in $\Hom(C,C')_n$ is called a morphism of degree $n$ and a boundary of $\Hom(C,C')_n$ is called
a homotopy of degree $n$.  The morphisms of degree $0$ are the morphisms in the category $\E_*$.

In this category, we have also a notion of $n$-cone:

Consider a morphism $f:X\fl Y$ in $\E_*$ and an integer $n\in\Z$. We set: $C=X\oplus Y$. So we have four linear maps: $i:Y\fl C$, $p:C\fl X$, $r:C\fl Y$ and $s:X\fl C$.
The map $i$ is an injection, $p$ is a projection, $r$ is a retraction of $i$ and $s$ is a section of $p$. There is a unique way to modify degrees and differentials on $C$
such that the following holds:
$$\partial^\circ i=n\hskip 24pt \partial^\circ r=-n\hskip 24pt \partial^\circ p=-1-n\hskip 24pt \partial^\circ s=1+n$$
$$d(i)=0\hskip 24pt d(p)=0\hskip 24pt d(r)=-(-1)^n fp\hskip 24pt d(s)=if$$
With these new degrees and differentials, $C$ is an $\E$-complex called the $n$-cone of $f$. If $n=0$, $C$ is the classical mapping cone, the map $i:Y\fl C$ is a
cofibration in $\E_*$ and we have an exact sequence in $\E_*$:
$$0\fl X\fl T(f)\fl C\fl 0$$
where $T(f)$ is the cylinder of $f$.

If $n=-1$, the map $p:C\fl X$ is also a morphism in $\E_*$.
\vskip 12pt
Let $F:\A\fl\B$ be an exact functor between two Waldhausen categories. We say that $F$ has the approximation property if the following holds:

$\bullet$ (App1) a morphism in $\A$ is an equivalence if and only if its image under $F$ is an equivalence

$\bullet$ (App2) for every $(X,Y)\in\A\times\B$ and every morphism $f:F(X)\fl Y$, there exist a morphism $\alpha:X\fl X'$ in $\A$
and a commutative diagram in $\B$:
$$\diagram{F(X)&\hfl{f}{}&Y\cr\vfl{\alpha}{}&&\vfl{=}{}\cr F(X')&\hfl{f'}{\sim}& Y\cr}$$
where $f'$ is an equivalence. We have the following theorem of Waldhausen ([Wa2]):
\vskip 12pt
\noi{\bf 2.3 Waldhausen approximation theorem:} {\sl Let $F:\A\fl\B$ be an exact functor between two essentially small saturated Waldhausen categories. Suppose $\A$ has
  a cylinder functor satisfying the cylinder axiom and $F$ has the approximation property. Then $F$ induces a homotopy equivalence in K-theory.}

\vskip 12pt
\noi{\bf 2.4 Lemma:} {\sl The functor $\Phi_2\times\Phi_3:\Mv\fl \D\times \C$ induces a homotopy equivalence in K-theory:
  $$K(\Mv)\build\fl_{}^\sim K(\D\times\C)=K(\D)\times K(\C)$$}
\vskip 12pt
\noi{\bf Proof:} The proof will be done by using Waldhausen K-theory.

The three exact categories $\Mv$, $\D$ and $\C$ are contained in the abelian categories $\Mv\v$, $\D\v$ and $\C\v$ respectively. Moreover each of these exact categories
are closed by kernel of epimorphisms in the corresponding abelian categories. Therefore, because of Gillet-Waldhausen theorem, it's enough to prove that the
functor $\Phi_2\times\Phi_3:\Mv_*\fl \D_*\times \C_*$ induces a homotopy equivalence in K-theory.

By replacing the category of equivalences of $\Mv_*$ by the category of morphisms $f:X\fl Y$ inducing a homology equivalence $\Phi_2(X)\build\fl_{}^\sim\Phi_2(Y)$ we get a new
Waldhausen category denoted by $\Mv'_*$.

Denote also by $\Mv^0_*$ the Waldhausen subcategory of $\Mv_*$ of objects $X$ with acyclic $\Phi_2(X)$. Because of the fibration theorem (see [Wa2]), the sequence:
$$\Mv^0_*\fl\Mv_*\fl\Mv'_*$$
induces a fibration in K-theory.

Denote by $\D^0_*$ the Waldhausen subcategory of $\D_*$ of acyclic complexes in $\D_*$. Since each morphism in $\D^0_*$ is an equivalence, $\D^0_*$ has trivial
K-theory. We have a commutative diagram:
$$\diagram{\Mv^0_*&\hfl{}{}&\Mv_*&\hfl{}{}&\Mv'_*\cr \vfl{\Phi_2\times\Phi_3}{}&&\vfl{\Phi_2\times\Phi_3}{}&&\vfl{\Phi_2}{}\cr
\D^0_*\times\C_*&\hfl{}{}&\D_*\times\C_*&\hfl{}{}&\D_*\cr}$$
where each line induces a fibration in K-theory. Therefore it will be enough to prove that functors $\Phi_2:\Mv'_*\fl\D_*$ and $\Phi_2\times\Phi_3:\Mv^0_*\fl
\D^0_*\times\C_*$ induce homotopy equivalences in K-theory or, equivalently, that $\Phi_2:\Mv'_*\fl\D_*$ and $\Phi_3:\Mv^0_*\fl\C_*$ induce homotopy equivalences in
K-theory. 

Because of the approximation theorem of Waldhausen, in order to prove that $\Phi_2:\Mv'_*\fl\D_*$ and $\Phi_3:\Mv^0_*\fl\C_*$ induces a homotopy equivalence in
K-theory, it's enough to show that these two functors have the approximation property.

The property (App1) it easy to check. Then the last thing to do is to show that $\Phi_2:\Mv'_*\fl\D_*$ and $\Phi_3:\Mv^0_*\fl\C_*$ have the property (App2).
\vskip 12pt
Consider an object $X\in\Mv'_*$, an object $F\in\D_*$ and a morphism $f:\Phi_2(X)\fl F$ in $\D_*$. The object $X$ is a triple $X=(E,M,\varphi)$ where $\E$ is a
$\D$-complex, $M$ is a $\C$-complex and $\varphi$ is an element in $\F(E,M)$ inducing a surjective morphism $T(E)\fl MR$. So $f$ is a morphism from $E$ to $F$ in $\D_*$. 

Denote by $E_1$ the $-1$-cone of $f$. So we have linear maps $i:F\fl E_1$, $p:E_1\fl E$, $r:E_1\fl F$ and $s: E\fl E_1$. The map $p$ is an epimorphism in $\D_*$, $i$ is
a morphism of degree $-1$ and we have:
$$d(s)=if\hskip 24pt d(r)=fp$$

Consider the triple $(E_1,M,\varphi p)$. Since $p$ is an epimorphism, this triple is an object $X_1\in\Mv_*$ and we have a morphism $g:X_1\fl X$. Denote by $Y$ the
$0$-cone of $g$. We have a cofibration $j:X\fl Y$ and four linear maps: $j:E\fl \Phi_2(Y)$, $q:\Phi_2(Y)\fl E_1$, $\rho:\Phi_2(Y)\fl E$ and $\sigma:E_1\fl\Phi_2(Y)$.
Moreover we have;
$$d(\sigma)=jp\hskip 24pt d(\rho)=-pq$$
Consider the linear map $g=f\rho+rq$. The differential $d(g)$ vanishes and $g$ is a morphism from $\Phi_2(Y)$ to $F$. It is easy to see that $g$ is surjective and its
kernel is isomorphic to the $0$-cone of the identity of $E$. Moreover we have: $gj=f$. Therefore $g:\Phi_2(Y)\fl F$ is a homology equivalence and the following diagram is
commutative:
$$\diagram{\Phi_2(X)&\hfl{f}{}&F\cr\vfl{j}{}&&\vfl{=}{}\cr \Phi_2(Y)&\hfl{g}{\sim}&F\cr}$$
Then the functor $\Phi_2:\Mv'_*\fl\D_*$ has the property (App2) and $\Phi_2:\Mv'_*\fl\D_*$ induces a homotopy equivalence in K-theory.
\vskip 12pt
Consider an object $X\in\Mv^0_*$, an object $N\in\C_*$ and a morphism $f:\Phi_3(X)\fl N$ in $\C_*$. The object $X$ is a triple $(E,M,\varphi)$ where $\E$ is an acyclic
$\D$-complex, $M$ is a $\C$-complex and $\varphi$ is an element in $\F(E,M)$ inducing a surjective morphism $T(E)\fl MR$. So $f$ is a morphism from $M$ to $N$.

Let $U$ be the $-1$-cone of the identity of $N$. Then $U$ is acyclic and we have an epimorphism $p:U\fl N$. Consider the composite morphism:
$$\varphi':T(Fs(U))\build\fl_{}^p T(Fs(N))\build\fl_\sim^e\sigma s(N)R\fl NR$$
This morphism is surjective and the triple $(E\oplus Fs(U),N,f\varphi\oplus\varphi')$ is an object $Y$ in $\Mv^0_*$. Moreover we have a morphism $g:X\fl Y$ inducing the
inclusion $E\subset E\oplus  Fs(U)$ and the morphism $f:M\fl N$, making the following diagram commutative:
$$\diagram{\Phi_3(X)&\hfl{f}{}&N\cr \vfl{g}{}&&\vfl{=}{}\cr\Phi_3(Y)&\hfl{=}{}&N\cr}$$
Hence $\Phi_3:\Mv^0_*\fl \C_*$ has the property (App2) and this functor induces a homotopy equivalence in K-theory.\cqfd
\vskip 12pt

Denote by $\Mv'$ the full subcategory of $\Mv$ consisting of objects $X\in\Mv$ such that $\Phi_1(X)$ belongs to $\R'$. Then an object $X\in\Mv$ belongs to $\Mv'$ if and
only if there is an isomorphism $\Phi_1(X)\simeq VR$ for some $V\in\C$. This category is exact and cofinal in $\Mv$. Moreover the functor $\Phi_1$ sends $\Mv'$ to $\R'$.

We define the following categories:

$\bullet$ The category $\E_0$ of modules in $\R'$ and isomorphisms

$\bullet$ The category $\E_1$ of objects in $\Mv'$ and morphisms inducing isomorphisms under $\Phi_1$

$\bullet$ The category $\E_2$ of objects in $\Mv'$ and epimorphisms inducing isomorphisms under $\Phi_1$

$\bullet$ The category $\E_3$ of objects in $\Mv'\times\V$, where a morphism in $\E_3$ from $(X,V)$ to $(Y,W)$ is an morphism $X\oplus V\fl Y\oplus W$ in $\E_2$ sending
$X$ to $Y$.

We have four functors $f_1:\E_1\fl \E_0$, $f_2:\E_2\fl \E_1$, $g_1:\E_3\fl\E_1$ and $g_2:\E_3\fl \E_2$ where $f_1$ is induced by $\Phi_1$, $f_2$ is the inclusion and
$g_1$ (resp. $g_2$) is the correspondence $(X,V)\mapsto X$ (resp. $(X,V)\mapsto X\oplus V$).

In order to prove a connective version of theorems 1, 2 and 3, we'll need to prove that the sequence $\V\subset \Mv'\fl \R'$ induces a fibration in K-theory and, for that,
it will be useful to prove that $f_1f_2$ is a homotopy equivalence.
\vskip 12pt
\noi{\bf 2.5 Lemma:} {\sl The functor $f_1$ is a homotopy equivalence.}
\vskip 12pt
\noi{\bf Proof:} It is enough to prove that the fiber category $f_1/U$ is contractible for each object $U\in\E_0=\R'$.

Consider an object $U$ in $\R'$ and denote by $\F$ the fiber category $f_1/U$. By applying lemma 1.14 with $X=0$ and $V=U$, we get an object $Y\in\Mv$ and an isomorphism
$U\simeq \Phi_1(Y)$. Hence the category $\F$ is nonempty

Consider two objects $(X_1,\e_1)$ and $(X_2,\e_2)$ in $\F$. For $i=1,2$, $\e_i$ is an isomorphism from $\Phi_1(X_i)$ to $U$.

By applying lemma 1.14 with $X=X_1\oplus X_2$ and $f=\e_1^{-1}\oplus \e_2^{-1}:U\fl \Phi_1(X)$, we get an object $Y\in\Mv$, a morphism $g:Y\fl X_1\oplus X_2$ and an
isomorphism $\e:U\build\fl_{}^\sim\Phi_1(Y)$ such that $\e_1^{-1}\oplus \e_2^{-1}$ is the composite morphism $U\build\fl_{}^\e\Phi_1(Y)\build\fl_{}^g\Phi_1(X)$.
Therefore $(Y,\e^{-1})$ is an object in $\F$ and we have two morphisms $(Y,\e^{-1})\fl (X_1,\e_1)$ and $(Y,\e^{-1})\fl (X_2,\e_2)$.

Consider two objects $(X_1,\e_1)$ and $(X_2,\e_2)$ in $\F$ and two morphisms $h_1$ and $h_2$ from $(X_1,\e_1)$ to $(X_2,\e_2)$. Denote by $X$ the kernel of $h_1-h_2$
(in $\S\v$).

For $i=1,2$, we have an exact sequence:
$$0\fl \Phi_1(X_i)\build\fl_{}^{\mu_i} T\Phi_2(X_i)\build\fl_{}^{\varphi_i}\Phi_3(X_i)R\fl 0$$
and then an exact sequence:
$$0\fl U\build\fl_{}^{\mu_i\e_i^{-1}} T\Phi_2(X_i)\build\fl_{}^{\varphi_i}\Phi_3(X_i)R\fl 0$$
So we get for each $i=1,2$ a commutative diagram with exact lines:
$$\diagram{0&\hfl{}{}&U&\hfl{\mu_1\e_1^{-1}}{}&T\Phi_2(X_1)&\hfl{\varphi_1}{}&\Phi_3(X_1)R&\hfl{}{}&0\cr
  &&\vfl{=}{}&&\vfl{h_i}{}&&\vfl{h_i}{}&&\cr
  0&\hfl{}{}&U&\hfl{\mu_2\e_2^{-1}}{}&T\Phi_2(X_2)&\hfl{\varphi_2}{}&\Phi_3(X_2)R&\hfl{}{}&0\cr}$$
and then an exact sequence:
$$0\fl U\fl T\Phi_2(X)\fl \Phi_3(X)R$$
inducing an isomorphism $\e:U\fl\Phi_1(X)$. By applying lemma 1.14 with $\e$, we get an object $(Y,u^{-1})\in\F$ and a morphism $k:(Y,u^{-1})\fl (X_1,\e_1)$ such that:
$h_1 k=h_2 k$.

Because of these properties, each fiber category $\F$ is cofiltered and then contractible. Hence the functor $f_1$ is a homotopy equivalence.\cqfd
\vskip 12pt
\noi{\bf 2.6 Lemma:} {\sl The functor $g_1$ is a homotopy equivalence.}
\vskip 12pt
\noi{\bf Proof:} Let $X$ be an object of $\Mv'$ and $\F$ be the fiber category $g_1/X$. Applying lemma 1.14 to the morphism $0\fl \Phi_1(X)$, we get an object $E\in\Mv$,
a morphism $\alpha:E\fl X$ such that: $\Phi_1(E)=0$ and $\Phi_2(E)\fl \Phi_2(X)$ and $\Phi_3(E)\fl \Phi_3(X)$ are epimorphisms. Then $E$ belongs to $\V$ and the morphism
$\alpha:E\fl X$ induces epimorphisms on $\Phi_2$ and $\Phi_3$. Therefore, for each morphism $f:Y\fl X$ in $\E_1$, the morphism $Y\oplus E\build\fl_{}^{f\oplus\alpha}X$
belongs to $\E_2$.

An object in $\F$ is a triple $(Y,V,f)$ where $(Y,V)$ belongs to $\Mv'\times\V$ and $f:Y\fl X$ is a morphism in $\E_1$. A morphism $\varphi:(Y,V,f)\fl(Y',V',f')$ is a
morphism $\varphi:Y\oplus V\fl Y'\oplus V'$ in $\E_2$ sending $Y$ to $Y'$ such that: $f=f'\varphi$.

The category $\F$ is nonempty because it contains the object $(X,0,\hbox{Id})$.

We have three functors $G_0,G_1,G_2$ from $\F$ to $\F$ sending each $(Y,V,f)\in\F$ to $G_0(Y,V,f)=(Y,V,f)$, $G_1(Y,V,f)=(Y,V\oplus E,f)$ and
$G_2(Y,V,f)=(X,0,\hbox{Id})$ respectively.

The inclusion $0\subset E$ induces a morphism $G_0\fl G_1$. The morphism $f\oplus0\oplus\alpha:Y\oplus V\oplus E\fl X$ induces a morphism
$(Y,V\oplus E,f)\fl(X,0,\hbox{Id})$ and we get a morphism $G_1\fl G_2$. Therefore the identity of $\F$ is homotopic to $G_1$ and then to $G_2$ which is constant. Hence
$\F$ is contractible and, since each fiber category of $g_1$ is contractible, $g_1$ is a homotopy equivalence.\cqfd
\vskip 12pt
\noi{\bf 2.7 Lemma:} {\sl The functor $g_2$ is a homotopy equivalence.}
\vskip 12pt
\noi{\bf Proof:} Let $X$ be an object in $\E_1$ and $\F$ be the fiber category $g_2/X$. An object in $\F$ is a triple $(Y,V,f)$ where $(Y,V)$ belongs to $\E_3$ and
$f:Y\oplus V\fl X$ is a morphism in $\E_2$. A morphism $\varphi:(Y,V,f)\fl (Y',V',f')$ in $\F$ is a morphism $\varphi:Y\oplus V\fl Y'\oplus V'$ in $\E_2$ such that
$\varphi(Y)\subset Y'$ and $f=f'\varphi$. Therefore it is easy to see that, for every $(Y,V,f)\in\F$, we have a unique morphism in $\F$ from $(Y,V,f)$ to $(X,0,\hbox{Id})$.
Hence $\F$ has a final object and is contractible. Since each fiber category of $g_2$ is contractible, $g_2$ is a homotopy equivalence.\cqfd
\vskip 12pt
\noi{\bf 2.8 Lemma:} {\sl The functor $f_2$ is a homotopy equivalence.}
\vskip 12pt
\noi{\bf Proof:} The inclusion $X\subset X\oplus V$  for all $(X,V)\in\E_3$ induces a morphism from $g_1$ to $f_2g_2$ and $g_1$ is homotopic to $f_2g_2$. But $g_1$ and
$g_2$ are homotopy equivalences. Therefore $f_2$ is a homotopy equivalence too.\cqfd
\vskip 12pt
\noi{\bf 2.9 Lemma:} {\sl The following diagram of exact categories:
$$\V\subset\Mv'\build\fl_{}^{\Phi_1}\R'$$
  induces a fibration in K-theory.}
\vskip 12pt
\noi{\bf Proof:} In [Wa1] (lemma 10.2 p. 206), Waldhausen proved that $Q\V\fl Q\Mv'\fl Q\F_R$ is a homotopic fibration in a situation similar to ours. In our situation
we'll prove, essentially in the same way, that the sequence $Q\V\fl Q\Mv'\fl Q\R'$ is a homotopic fibration.

Following Waldhausen's notations, if $\F$ is an exact subcategory of an exact category $\E$, we have a bicategory $Q^{ep}(\E,\F)$ where the horizontal maps form the
Quillen's category $Q(\E)$ and the vertical morphisms form the category of epimorphisms in $\E$ with kernel in $\F$. In particular we have an equivalence between
$Q(\E)$ and $Q^{ep}(\E,0)$ and $Q(\E)$ may be considered as a bicategory.

We have a commutative diagram of bicategories:
$$\diagram{Q(\V)&\fl{}{}&Q^{ep}(\V,\V)\cr\vfl{}{}&&\vfl{}{}\cr Q(\Mv')&\fl{}{}&Q^{ep}(\Mv',\V)\cr}$$
which is homotopically cartesian. Moreover $Q^{ep}(\V,\V)$ is contractible. Therefore the diagram of bicategories:
$$Q(\V)\fl Q(\Mv')\fl Q^{ep}(\Mv',\V)$$
is a homotopic fibration.

On the other hand the morphism $f=f_1f_2:\E_2\fl \R'$ induces a morphism $f_*:Q^{ep}(\Mv',\V)\fl Q^{ep}(\R',0)$ and we want to prove that $f_*$ is a homotopy equivalence.
Actually, the proof of lemma 10.2 in [Wa1] works exactly the same in our situation except maybe for the sublemma (p. 209). In this sublemma, we have a filtered object
$$M_1\subset M_2\subset\dots M_{n-1}\subset M_n$$
in $\Mv'$ where each quotient $M_i/M_{i-1}$ is in $\Mv'$. Since $\Phi_1$ is exact and each module in $\R'$ is projective, the morphism
$\Phi_1(M_n)\fl\Phi_1(M_n/M_{n-1})$ is surjective and has a section $s$ from $U=\Phi_1(M_n/M_{n-1})$ to $\Phi_1(M_n)$.

Because of lemma 1.14, there are an object $N\in\Mv$, a morphism $g:N\fl M_n$ inducing epimorphisms on $\Phi_2$ and $\Phi_3$, and an isomorphism
$U\build\fl_{}^\sim\Phi_1(N)$ making the following diagram commutative:
$$\diagram{U&\hfl{\sim}{}&\Phi_1(N)\cr\vfl{=}{}&&\vfl{g}{}\cr U&\hfl{s}{}&\Phi_1(M_n)\cr}$$

Therefore the morphism $N\fl M_n$ induces epimorphisms on $\Phi_2$ and $\Phi_3$ and the composite morphism $N\fl M_n\fl M_n/M_{n-1}$ is an epimorphism with kernel in
$\V$. Hence the sublemma can be proven in our situation and, since $f$ is a homotopy equivalence, the proof of lemma 10.2 applies completely here. The lemma follows.
\cqfd
\vskip 12pt
\noi{\bf 2.10 Lemma:} {\sl Let $\R''$ be the full subcategory of $\R$ generated by the image of $\Phi_1:\Mv\fl \R$. Then the diagram:
  $$\V\subset\Mv\build\fl_{}^{\Phi_1}\R''$$
  induces a fibration in K-theory.}
\vskip 12pt
\noi{\bf Proof:} Let $X=(E,M,\varphi)$ be an object in $\Mv$. Since $E$ is projective in $\D$, there is a module $E_1\in\D$ such that $E\oplus E_1$ is free in case 1 or 3
and on the form $E\oplus E_1=(F_A,F_B)\in\A\times\B$ with $F_A$ and $F_B$ free. Therefore $T(E\oplus E_1)$ is free in all cases. Let $X_1$ be the object
$(E_1,0,0)\in\Mv$. We have an exact sequence in $\R$:
$$0\fl \Phi_1(X\oplus X_1)\fl T(E\oplus E_1)\fl MR\fl 0$$
and $\Phi_1(X\oplus X_1)$ is stably in $\R'$. Hence there is another object $X_2\in\Mv$ such that $\Phi_1(X\oplus X_1\oplus X_2)$ belongs to $\R'$ and $\Mv'$ is cofinal in
$\Mv$.

Consider the following commutative diagram:
$$\diagram{\Mv'&\hfl{\Phi_1}{}&\R'\cr\vfl{}{}&&\vfl{}{}\cr\Mv&\hfl{\Phi_1}{}&\R''\cr}\leqno{(D)}$$
where the vertical maps are the canonical cofinal inclusions.

Let $K$ be the fiber product of $K_0(\Mv)$ and $K_0(\R')$ over $K_0(\R'')$. The commutativity of the diagram $K_0(D)$ induce a map $\lambda:K_0(\Mv')\fl K$.
Since the map $K_0(\Mv')\fl K_0(\Mv)$ is injective (by cofinality) and factorizes through $K$, the map $\lambda$ is injective.

Let $w=(u,v)$ be an element in $K$. Then we have: $u\in K_0(\Mv)$, $v\in K_0(\R')$ and $\Phi_1(u)$ and $v$ are the same in $K_0(\R'')$.

For every object $X$ in some exact category $\A$, the class of $X$ in the Grothendieck group $K_0(\A)$ will be denoted by $[X]$.

So there are two objects $X$, $Y$ in $\Mv$ such that: $u=[X]-[Y]$. Since $\Mv'$ is cofinal in $\Mv$, there is an object $Y_1$ in $\Mv$ such that $Y\oplus Y_1$ is in
$\Mv'$. Let us set: $w'=w+\lambda[Y\oplus Y_1]$ and $X'=X\oplus Y_1$. We have: $w'=(u',v')$ with $u'=[X']$ and $[\Phi_1(X')]$ belongs to $K_0(\R')$.
Then $\Phi_1(X')$ is stably isomorphic to a module in $\R'$. Up to adding to $X$ (and then to $X'$) an object in $\Mv'$ on the form $(E,0,0)$, we may as well suppose that
$\Phi_1(X')$ belongs to $\R'$ and then that $X'$ belongs to $\Mv'$.

Therefore we have: $w'-\lambda[X']=(0,v'')$ where $v''$ is an element in $K_0(\R')$ killed in $K_0(\R'')$. But the morphism $K_0(\R')\fl K_0(\R'')$ is injective. So we
have:
$$w'-\lambda[X']=0\ \ \Longrightarrow\ \ w'=\lambda(X')$$
and $w'$ and then $w$ are in the image of $\lambda$. Therefore $\lambda$ is surjective and then bijective.

Consequently the diagram $K_0(D)$ is exact (cartesian and cocartesian) and, by cofinality, the diagram $K(D)$ is homotopically cartesian. The lemma follows.\cqfd
\vskip 12pt
As a consequence we get the following result, which is, in some sense, a connective version of theorems 1, 2 and 3:
\vskip 12pt
\noi{\bf 2.11 Proposition:} {\sl Let $X$ be the homotopy fiber  of the map $K(\C_2)\fl K(\D)\times F(\C)$ induced by the functor $F\times\sigma:\C_2\fl \D\times\C$.
Then there is a natural homotopy equivalence $\Omega K(R)\build\fl_{}^\sim Nil(C_2,S)\times X$.}
\vskip 12pt
\noi{\bf Proof:} Because of lemma 1.12 the functor $\Phi:\Nil(C_2,S)\fl\V$ is an equivalence of categories and the lemma 2.10 implies that the following diagram:
$$\Nil(C_2,S)\build\fl_{}^\Phi \Mv\build\fl_{}^{\Phi_1}\R''$$
induces a fibration in K-theory. Hence $\Omega K(\R'')\simeq\Omega K(R)$ is homotopically equivalent to the homotopy fiber of the map induced by $\Phi:\Nil(C_2,S)\fl\Mv$
in K-theory.

Because of lemma 2.4, we have a commutative diagram
$$\diagram{\Nil(C_2,S)&\hfl{\Phi}{}&\Mv\cr\vfl{=}{}&&\vfl{}{\Phi_2\times\Phi_3}\cr\Nil(C_2,S)&\hfl{\Phi'}{}&\D\times\C\cr}$$
where the functor $\Phi_2\times\Phi_3$ induces a homotopy equivalence in K-theory. Therefore $\Omega K(R)$ is homotopically equivalent to the homotopy fiber of the map
$\Phi':K(\Nil(C_2,S))\fl K(\D\times\C)$.

The functor $\Phi'$ sends an object $(H,\theta)\in\Nil(C_2,S)$ to the pair $(F(H),\sigma(H))\in \D\times\C$ and $\Phi'$ factorizes by the forgetful map $\Nil(C_2,S)\fl
\C_2$. Therefore $\Omega K(R)$ is homotopically equivalent to the homotopy fiber of $\Phi'':K(C_2)\times Nil(C_2,S)\fl K(\D)\times K(\C)$ where $\Phi''$ is trivial on
$Nil(C_2,S)$ and induced by the functor $F\times\sigma$ on $K(\C_2)$. The lemma follows.\cqfd
\vskip 12pt
\noi{\bf 2.12 The spectra $\underline K$ and $\underline Nil$.}
\vskip 12pt
The K-theory of Quillen is a functor K from the category of rings to the category of infinite loop spaces. There are different methods to construct a so called negative
K-theory: that is a functor $K'$ from
the category of rings to the category $\Omega sp$ of $\Omega$-spectra (see [B] and [KV]) and a natural homotopy equivalence from $K(A)$ to the 0-th term of $K'(A)$
such that the following sequence is exact for every ring $A$ and every integer $i\in\Z$:
$$0\fl K_i(A)\fl K_i(A[t])\oplus K_i(A[t^{-1}])\fl K_i(A[t,t^{-1}])\fl K_{i-1}(A)\fl 0$$
where $K_i(A)$ is the $i$-th homotopy group of $K'(A)$. This exact sequence was proven by Bass [B] for $i=1$ and generalized by Quillen for $i>0$ [Q]. The morphisms of
this exact sequence are induced by the inclusions $A\subset A[t]\subset A[t,t^{-1}]$, $A\subset A[t^{-1}]\subset A[t,t^{-1}]$ except for the map
$\partial:K_i(A[t,t^{-1}])\fl K_{i-1}(A)$. But this map $\partial$ has a section induced by the multiplication by $t\in K_1(\Z[t,t^{-1}])$. Therefore the exact sequence
above is natural in $A$.

Inspired by the Karoubi-Villamayor method [KV], we'll define our version $\underline K(A)$ of negative K-theory as follows:

Denote by $E$ the set of infinite square matrices with entries in $\Z$ having only finitely many nonzero entries in each row and each column. This ring has a two-sided
ideal $M(\Z)$ of matrices having only finitely many nonzero entries. So we set: $\Sigma=E/M(\Z)$.

For every ring $A$, $EA=E\otimes_\Z A$ and $\Sigma A=\Sigma\otimes_\Z A$ are rings and the morphism $f:EA\fl \Sigma A$ is a surjective ring homomorphism. It is easy to see
that the kernel of $f$ is isomorphism, as a pseudo ring, to $M(A)$ and that $EA$ is a flasque ring. Therefore (see [KV]) we have a natural  homotopy equivalence:
$$K(A)\build\fl_{}^\sim \Omega K(\Sigma A)$$
and the sequence $K(\Sigma^n A)$ is an $\Omega$-spectrum. This spectrum will be denoted by $\underline K(A)$ and $\underline K$ is a negative K-theory. For each integer
$i\in\Z$ we set also: $K_i(A)=\pi_i(\underline K(A))$.
\vskip 12pt
\noi{\bf 2.13 Lemma:} {\sl Let $(A,S)$ be a left-flat bimodule. Then for each integer $n\geq0$, $(\Sigma^n A,\Sigma^n S)$ is a left-flat bimodule and the sequence
  $Nil(\Sigma^n A,\Sigma^n S)$ is an $\Omega$-spectrum denoted by $\underline Nil(A,S)$. Moreover we have a natural homotopy equivalence from $Nil(A,S)$ to the 0-th term
  of $\underline Nil(A,S)$ and, for each integer $i\in\Z$, we have an exact sequence which is natural on the left-flat bimodule $(A,S)$ :
  $$0\fl Nil_i(A,S)\fl Nil_i(A[t],S[t])\oplus Nil(A[t^{-1}],S[t^{-1}])\fl\hskip 100pt$$
  $$\hskip 190pt Nil_i(A[t,t^{-1}],S[t,t^{-1}])\fl Nil_{i-1}(A,S)\fl 0$$
  where $Nil_i(?)$ is the $i$-th homotopy group of $\underline Nil(?)$.}
\vskip 12pt
\noi{\bf Proof:} Let $(A,S)$ be a left-flat bimodule. For every ring $B$, $BS=B\otimes_Z S$ is a $BA$-bimodule. Since $S$ is flat on the left $S$, as a left $A$-module,
is isomorphic to a filtered colimit of free left $A$-modules $E_i$ and $BS$ is also  isomorphic to a filtered colimit of free $BA$-modules $BE_i$. Therefore $BS$ is
flat on the left and $(BA,BS)$ is a left-flat bimodule. In particular each $(\Sigma^n A,\Sigma^n S)$ is a left-flat bimodule.
Moreover we have a natural isomorphism of rings: $\Sigma(A[S])\simeq (\Sigma A)[\Sigma S]$.

Consider the case 1 (with $C=A$). Because of proposition 2.11, we have a natural homotopy equivalence:
$$\Omega K(A[S])\build\fl_{}^\sim Nil(A,S)\times X$$
where $X$ is the homotopy fiber of the map: $K(\C_2)\build\fl_{}^{F\times\sigma} K(\D\times\C)$.
But in case 1, $\C_2$, $\D$ and $\C$ are equal to the category $\A$ of finitely generated projective right $A$-modules and functors $F$ and $\sigma$ are the identity.
Therefore $X$ is nothing else but the loop space of $K(A)$. So we get a natural homotopy equivalence:
$$\Omega K(A[S])\build\fl_{}^\sim Nil(A,S)\times \Omega K(A)$$
By naturality, we get a homotopy equivalence from $Nil(A,S)$ to the homotopy fiber of the map $\Omega K(A[S])\fl\Omega K(A)$ induced by the canonical ring homomorphism
$A[S]\fl A$.

We have a commutative diagram:
$$\diagram{K(A[S])&\hfl{}{}&K(A)\cr\vfl{\sim}{}&&\vfl{\sim}{}\cr\Omega K(\Sigma(A[S]))&\hfl{}{}&\Omega K(\Sigma A)\cr}$$
where the horizontal maps are induced by the canonical morphism $A[S]\fl A$ and the vertical maps are homotopy equivalences.

But $\Sigma(A[S])$ is isomorphic to $(\Sigma A)[\Sigma S]$. So we get a commutative diagram:
$$\diagram{\Omega K(A[S])&\hfl{}{}&\Omega K(A)\cr\vfl{\sim}{}&&\vfl{\sim}{}\cr\Omega^2K((\Sigma A)[\Sigma S]))&\hfl{}{}&\Omega^2K(\Sigma A)\cr}$$
inducing a homotopy equivalence $Nil(A,S)\build\fl_{}^\sim \Omega Nil(\Sigma A,\Sigma S)$. Then the sequence of spaces $Nil(\Sigma^n A,\Sigma^n,S)$ is a well defined
$\Omega$-spectrum $\underline Nil(A,S)$ and we have a natural decomposition:
$$\Omega\underline K(A[S])\build\fl_{}^\sim \underline Nil(A,S)\times \Omega\underline K(A)$$

Let us set: $Nil_i(A,S)=\pi_i(\underline Nil(A,S))$. Then we have, for every left-flat bimodule $(A,S)$ and every integer $i\in\Z$ an isomorphism
$K_i(A[S])\simeq Nil_{i-1}(A,S)\oplus K_i(A)$.

For the ring $A$ we have, for every integer $i\in\Z$, the following exact sequence $S_i(A)$:
$$0\fl K_i(A)\fl K_i(A[t])\oplus K_i(A[t^{-1}])\fl K_i(A[t,t^{-1}])\fl K_{i-1}(A)\fl 0$$
For $B=\Z[t]$ or $B=\Z[t^{-1}]$ or $B=\Z[t,t^{-1}]$, the ring $B(A[S])$ is isomorphic to $(BA)[BS]$. Then the sequence $S_{i+1}(A[S])$ decomposes into $S_{i+1}(A)$ and
the following sequence:
$$0\fl Nil_i(A,S)\fl Nil_i(A[t],S[t])\oplus Nil_i(A[t^{-1}],S[t^{-1}])\fl\hskip 100pt $$
$$\hskip 180 pt Nil_i(A[t,t^{-1}],S[t,t^{-1}])\fl Nil_{i-1}(A,S)\fl 0$$
which, as a consequence, is exact.\cqfd
\vskip 12pt
\noi{\bf 2.14 Proofs of theorems 1, 2 and 3.}
\vskip 12pt
In the proof of lemma 2.13, we have constructed the $\Omega$-spectrum $\underline Nil(?)$ and the first two properties of theorem 1 have already been proven.

Suppose $A$ is regular coherent on the right, i.e. every finitely presented right $A$-module has a finite resolution by finitely generated projective modules. The category
$\A$ of finitely generated projective
right $A$-modules is contained in the category $\A'$ of finitely presented right $A$-modules. The category $\Nil(A,S)$ is also contained in the category $\Nil'(A,S)$ of
pairs $(H,\theta)$ with $H\in\A'$ and $\theta:H\fl HS$ nilpotent. Moreover $\A$ is stable in $\A'$ under extension and kernel of admissible epimorphism and the inclusion
$\Nil(A,S)\subset \Nil'(A,S)$ have the same properties. Therefore, by the resolution theorem [Q], the inclusions $\A\subset \A'$ and $\Nil(A,S)\subset \Nil'(A,S)$ induce
homotopy equivalences in K-theory.

On the other hand, we have an inclusion $\A'\fl \Nil'(A,S)$ sending $H$ to the pair $(H,0)$ and every object $(H,\theta)$ in $\Nil'(A,S)$ has a finite filtration with
subquotients in $\A'$. Moreover $\A'$ and $\Nil'(A,S)$ are abelian categories and $\A'$ is closed in $\Nil'(A,S)$ under subobjects and quotients. Therefore, by the
devissage theorem [Q], the inclusion $\A'\subset \Nil'(A,S)$ induces a homotopy equivalence in K-theory. As a consequence the inclusion $\A\subset
\Nil(A,S)$ induce a homotopy equivalence in K-theory and the space $Nil(A,S)$ is contractible.

Hence for every left-flat bimodule $(A,S)$, with $A$ regular coherent, the space $Nil(A,S)$ is contractible and $Nil_i(A,S)$ is trivial for every $i\geq0$.

Consider the following property $E(n)$ where $n$ is any integer:

$\bullet$ For every left-flat bimodule $(A,S)$ with $A$ regular coherent on the right, the module $Nil_i(A,S)$ is trivial for every $i\geq n$.

The property $E(0)$ is then satisfied. Suppose $E(n)$ is true and take a left-flat bimodule $(A,S)$ with $A$ regular coherent on the right. Then $A[t]$, $A[t^{-1}]$ and
$A[t,t^{-1}]$ are also regular coherent on the right and, in the exact sequence:
$$0\fl Nil_n(A,S)\fl Nil_n(A[t],S[t])\oplus Nil_n(A[t^{-1}],S[t^{-1}])\fl\hskip 100pt $$
$$\hskip 180 pt Nil_n(A[t,t^{-1}],S[t,t^{-1}])\fl Nil_{n-1}(A,S)\fl 0$$
all the modules are zero except $Nil_{n-1}(A,S)$. Hence the module $Nil_{n-1}(A,S)$ is also trivial and the property $E(n-1)$ is true.

By induction $E(n)$ is satisfied for all $n$. Hence $\underline Nil(A,S)$ is contractible for every left-flat bimodule $(A,S)$ with $A$ regular coherent on the right and
theorem 1 is proven.
\vskip 12pt
Consider the case 2. The ring $R$ is defined by the cocartesian diagram:
$$\diagram{C&\hfl{\alpha}{}&A\cr\vfl{\beta}{}&&\vfl{}{}\cr B&\hfl{}{}&R\cr}\leqno{(D)}$$
By tensoring by $\Sigma^n$, we get a cocartesian diagram:
$$\diagram{\Sigma^nC&\hfl{\alpha}{}&\Sigma^nA\cr\vfl{\beta}{}&&\vfl{}{}\cr \Sigma^nB&\hfl{}{}&\Sigma^nR\cr}\leqno{(\Sigma^n D)}$$

We remark that the morphism $\alpha:\Sigma^n C\fl \Sigma^n A$ (resp. $\beta:\Sigma^n C\fl \Sigma^n B$) is pure with complement $\Sigma^n A'$ (resp. $\Sigma^n B'$).
Therefore in this new situation, we get new rings: $\Sigma^n C$, $\Sigma^n A$, $\Sigma^n B$, $\Sigma^n R$, $\Sigma^n C\times\Sigma^n C$ and a new bimodule $\Sigma^n S$.

Because of proposition 2.11, we have a homotopy equivalence between $\Omega K(\Sigma^nR)$ and $Nil(\Sigma^n C\times\Sigma^n C,\Sigma^n S)\times X_n$ where $X_n$ is the
homotopy fiber of $f:K(\Sigma^n C\times\Sigma^n C)\fl K(\Sigma^n A)\times K(\Sigma^n B)\times K(\Sigma^n C)$.

The morphism $f$ is induced by the functor $\C\times\C\fl \A\times\B\times\C$ sending $(M,M')\in\C\times\C$ to $(MA,M'B,M\oplus M')$. Denote by $f_n(\alpha)$ (resp.
$f_n(\beta)$) the map $K(\Sigma^n C)\fl K(\Sigma^n A)$ (resp. $K(\Sigma^n C)\fl K(\Sigma^n B)$) induced by $\alpha$ (resp. $\beta$). Then $X_n$ is homotopy equivalent
to the homotopy fiber of $f_n(\alpha)-f_n(\beta): K(\Sigma^n C)\fl K(\Sigma^n A)\times K(\Sigma^n B)$ and we have a homotopy cartesian diagram of spaces:
$$\diagram{\Omega H(\Sigma^n C)&\hfl{\alpha}{}&\Omega K(\Sigma^n A))\cr\vfl{\beta}{}&&\vfl{}{}\cr \Omega K(\Sigma^n B)&\hfl{}{}& X_n\cr}$$

By naturality of the homotopy equivalence: $K(?)\simeq\Omega K(\Sigma ?)$, we get a homotopy equivalence $X_n\fl \Omega X_{n+1}$ and the sequence of $X_n$ defines an
$\Omega$-spectrum $X$ together with a homotopy cartesian diagram of spectra:
$$\diagram{\underline K(C)&\hfl{\alpha}{}&\underline K(A)\cr\vfl{\beta}{}&&\vfl{}{}\cr \underline K(B)&\hfl{}{}&\Omega^{-1} X\cr}$$
and that finishes the proof of theorem 2.
\vskip 12pt
Consider now the case 3. We proceed as before and we get a homotopy equivalence between $\Omega K(\Sigma^nR)$ and $Nil(\Sigma^n C\times\Sigma^n C,\Sigma^n S)\times X_n$
where $X_n$ is the homotopy fiber of the map $f:K(\Sigma^n C\times\Sigma^n C)\fl K(\Sigma^n A)\times K(\Sigma^n C)$ induced by the functor
$F\times\sigma:\C\times\C\fl \A\times\C$ sending $(M,M')\in\C\times\C$ to $(M{}_\alpha A\oplus M'{}_\beta A,M\oplus M')$.

Denote by $f_n(\alpha)$ (resp. $f_n(\beta)$) the map $K(\Sigma^n C)\fl K(\Sigma^n A)$ induced by $\alpha$ (resp. $\beta$). Then $X_n$ is homotopy equivalent
to the homotopy fiber of $f_n(\alpha)-f_n(\beta)$.

As before, we get a homotopy equivalence $X_n\fl \Omega X_{n+1}$ and the sequence of $X_n$ defines an
$\Omega$-spectrum $X$ which is the homotopy fiber of the map $f(\alpha)-f(\beta):\underline K(C)\fl \underline K(A)$. Therefore we have a homotopy fibration of spectra:
$$\underline K(C)\build\fl_{}^{f(\alpha)-f(\beta)} \underline K(A)\fl \Omega^{-1}X$$
and that finishes the proof of theorem 3.\cqfd
\vskip 24pt
\noi{\bf 3. Properties of the functor $\underline Nil$.}
\vskip 12pt
Consider two rings $A$ and $B$. Every right $A\times B$-module $M$ is determined by two right modules $M_a$ and $M_b$, where $M_a$ is an $A$-module and $M_b$ is a
$B$-module. By setting: $R_a=A$ and $R_b=B$, we see that $M_i$ is a right $R_i$-module for each $i\in\{a,b\}$.

If $E$ is an $A\times B$-bimodule, $E$ is determined by four bimodules ${}_a E_a$, ${}_a E_b$, ${}_b E_a$ and ${}_b E_b$, and for each $i,j$ in $\{a,b\}$, ${}_i E_j$ is
a $(R_i,R_j)$-bimodule.

Suppose $f:M\fl ME$ is a morphism of right $A\times B$-modules. Then $f$ is determined by four morphisms ${}_i f_j:M_j\fl M_i\ {}_i E_j$ and ${}_i f_j$ is a morphism
of right $R_j$-modules.
\vskip 12pt
\noi{\bf 3.1 Lemma:} {\sl Let $A$ and $B$ be two rings and $E$ be an $A\times B$-bimodule. Suppose $E$ is flat on the left. Then the correspondences
  $$(M,f)\mapsto (M_b,{}_b f_b)$$
  $$(M,f)\mapsto (M_a,{}_a f_a+\build\sum_{k\geq0}^{}{}_a f_b({}_b f_b)^k {}_b f_a)$$
  induce two well defined functors:
  $$\Phi_1:\Nil(A\times B,E)\fl\Nil(B,{}_b E_b)$$
  $$\Phi_2:\Nil(A\times B,E)\fl\Nil(A,{}_aE_a\oplus\build\oplus_{k\geq0}^{}{}_a E_b({}_b E_b)^k{}_b E_a)$$
  
    Moreover, if $E_b$ is flat on the right, these functors induce a homotopy equivalence of spectra:
    $$\underline Nil(A\times B,E)\build\fl_{}^\sim \underline Nil(B,{}_b E_b)\times\underline Nil(A,{}_a E_a\build\oplus_{k\geq0}^{}{}_a E_b({}_b E_b)^k{}_b E_a)$$}
\vskip 12pt
This lemma will be proven in 3.7.
\vskip 12pt
Using this result, we are able to prove theorem 4. Consider two rings $A$ and $B$, an $(A,B)$-bimodule $S$ and a $(B,A)$-bimodule $T$. Suppose $S$ and $T$ are flat on
both sides. Define the $A\times B$-bimodule $E$ by:
$${}_a E_b=S\hskip 24pt {}_b E_a=T\hskip 24pt {}_a E_a={}_b E_b=0$$
Actually, this bimodule is the bimodule $S\oplus T$, where $S$ and $T$ are considered as $A\times B$-bimodules (via the projections $A\times B\fl A$ and $A\times B\fl B$).
Moreover $E=S\oplus T$ is flat on both sides.

Applying lemma 3.1, we get two functors:
$$\Phi_1:\Nil(A\times B,E)\fl\Nil(B,0)$$
$$\Phi_2:\Nil(A\times B,E)\fl\Nil(A,ST)$$
and a homotopy equivalence of spectra:
$$\Phi(S,T):\underline Nil(A\times B,E)\build\fl_{}^\sim \underline Nil(A,ST)$$
By exchanging the roles of $A$ and $B$ , and $S$ and $T$ we get also a homotopy equivalence of spectra:
$$\Phi(T,S)=\underline Nil(A\times B,E)\build\fl_{}^\sim \underline Nil(B,TS)\eqno{\cqfd}$$
\vskip 12pt
\noi{\bf 3.2 Proof of theorem 5.}
\vskip 12pt
Let $A$ be a ring, $I$ be a set and $S_i$, $i\in I$ be a family of $A$-bimodules flat on both sides. The direct sum of the $S_i$'s will be denoted by $S$.

Let $(M,f)$ be an object of the category $\Nil(A,S)$. The morphism $f:M\fl MS$ decomposes into a finite sum: $f=\build\sum_{i\in I}^{} f_i$, where $f_i$ is a morphism
from $M$ to $MS_i$.

If $u=i_1 i_2 \dots i_p$ is a word in $W(I)$, we set:
$$f_u=f_{i_1}f_{i_2}\dots f_{i_p}\ \ \hbox{and}\ \ \ S_u=S_{i_1}S_{i_2}\dots S_{i_p}$$
If $J$ is a subset of $W(I)$, we set also:
$$f_J=\build\sum_{u\in J}^{} f_u\ \ \hbox{and}\ \ \ S_J=\build\oplus_{u\in J}^{} S_u$$
We check that $f_u$ is a morphism from $M$ to $M S_u$ and $f_J$ is a morphism from $M$ to $MS_J$.

Let $u$ be a non empty word in $W(I)$. Then the correspondence $(M,f)\mapsto (M,f_u)$ induces an exact functor $\varphi_u:\Nil(A,S)\fl\Nil(A,S_u)$. If $J$ is a subset of
$W(I)$ that does not contains the empty word, the correspondence $(M,f)\mapsto (M,f_J)$ induces also an exact functor $\varphi_J:\Nil(A,S)\fl\Nil(A,S_J)$. These functors
are compatible with tensor product with any power of $\Sigma$ and induce morphisms of spectra on K-theory.

Consider two words $u$ and $v$ in $W(I)$. We have four $A$-bimodules: $S_u$, $S_v$, $S_{uv}$ and $S_{vu}$. Using notations in 1.4, for each $A$-bimodule $T$ and each
$(i,j)\in\{1,2\}^2$ we get an $A\times A$-bimodule ${}^i T^j$ where $A\times A$ acts on the left of $T$ via the $i$-th projection and on the right via the $j$-th
projection. We have a commutative diagram of exact categories:
$$\diagram{\Nil(A,S)&\hfl{F}{}&\Nil(A\times A,{}^1 S_u^2\oplus {}^2 S_v^1)&\hfl{\psi(u,v)}{}&\Nil(A,S_{uv})\cr\vfl{=}{}&&\vfl{G}{\sim}&&\cr
  \Nil(A,S)&\hfl{F'}{}&\Nil(A\times A,{}^2 S_u^1\oplus {}^1 S_v^2)&\hfl{\psi(v,u)}{}&\Nil(A,S_{vu})\cr}$$
where $F$ and $F'$ are defined by:
$$F(M,f)=(M^1\oplus M^2,f_u+f_v)\hskip 48pt F'(M,f)=(M^1\oplus M^2,f_v+f_u)$$
for each $(M,f)\in\Nil(A,S)$ and the morphisms $G$, $\psi(u,v)$ and $\psi(v,u)$ are defined by sending each $(M,M',f+f')$ (with $M,M'\in\A$, $f:M\fl M' S_v$ and
$f':M'\fl M S_u$) to
$$G(M,M',f+f')=(M',M,f+f')\hskip 24pt \psi(u,v)(M,M',f+f')=(M,f'f)$$
$$\psi(v,u)(M,M',f+f')=(M',ff')$$

Because of theorem 4, this diagram induces a homotopy commutative diagram of $\Omega$-spectra:
$$\diagram{\underline Nil(A,S)&\hfl{F}{}&\underline Nil(A\times A,{}^1 S_u^2\oplus {}^2 S_v^1)&\hfl{\psi(u,v)}{\sim}&\underline Nil(A,S_{uv})\cr
  \vfl{=}{}&&\vfl{G}{\sim}&&\vfl{H}{\sim}\cr
  \underline Nil(A,S)&\hfl{F'}{}&\underline Nil(A\times A,{}^2 S_u^1\oplus {}^1 S_v^2)&\hfl{\psi(v,u)}{\sim}&\underline Nil(A,S_{vu})\cr}$$
where $\psi(u,v)$ and $\psi(v,u)$ are homotopy equivalences and $H$ is the map $\psi(u,v)^{-1}G\psi(v,u)$ (up to homotopy).

On the other hand we have: $\psi(u,v)F=\varphi_{uv}$ and $\psi(v,u)F'=\varphi_{vu}$. Then we have a homotopy commutative diagram:
$$\diagram{\underline Nil(A,S)&\hfl{\varphi_{uv}}{}&\underline Nil(A,S_{uv})\cr\vfl{=}{}&&\vfl{H}{\sim}\cr
  \underline Nil(A,S)&\hfl{\varphi_{vu}}{}&\underline Nil(A,S_{vu})\cr}$$
and the homotopy class of the map $\varphi_u$ depends only on the class of $u$ in $CW(I)$. Therefore, to prove the theorem, its is enough to
prove it for some admissible set $X$.
\vskip 12pt
\noi{\bf 3.3 Lemma:} {\sl Suppose $I$ has exactly two elements $i$ and $j$. Let $J\subset W(I)$ be the set of elements $i^kj$, for $k\geq0$. Then the functors:
  $\varphi_i:\Nil(A,S)\fl \Nil(A,S_i)$ and $\varphi_J:\Nil(A,S)\fl \Nil(A,S_J)$ induce a homotopy equivalence of spectra:
  $$\underline Nil(A,S)\build\fl_{}^\sim \underline Nil(A,S_i)\times\underline Nil(A,S_J)$$}
\vskip 12pt
\noi{\bf Proof:} We apply lemma 3.1 in the following case:
$$B=A\hskip 24pt {}_a E_a=S_i\hskip 24pt {}_b E_b=0\hskip 24pt {}_b E_a=A\hskip 24pt {}_a E_b=S_j$$
and we get a homotopy equivalence of spectra:
$$\underline Nil(A\times A,E)\build\fl_{}^\sim \underline Nil(A,S_i\oplus S_j)=\underline Nil(A,S)$$
induced by the functor $(P,Q,f)\mapsto (P,{}_a f_a+{}_a f_b\ {}_b f_a)$.

By exchanging the role of $a$ and $b$, we have also a homotopy equivalence of spectra:
$$\underline Nil(A\times A,E)\build\fl_{}^\sim \underline Nil(A,S_i)\times\underline Nil(A,\build\oplus_{k\geq0}^{} S_i^k S_j)=\underline Nil(A,S_i)\times
\underline Nil(A,S_J)$$
and this map is induced by the two functors:
$$(P,Q,f)\mapsto (P,{}_a f_a)$$
$$(P,Q,f)\mapsto (Q,\build\sum_{k\geq0}^{}{}_b f_a({}_a f_a)^k{}_a f_b)$$
Because of theorem 4, this last functor is, after applying the functor $\underline Nil$, equivalent to the functor:
$$(P,Q,f)\mapsto (P,\build\sum_{k\geq0}^{}({}_a f_a)^k{}_a f_b\ {}_b f_a)$$
Therefore we get a homotopy equivalence of spectra:
$$\underline Nil(A,S)\build\fl_{}^\sim \underline Nil(A,S_i)\times\underline Nil(A,S_J)$$
induced by the two functors:
$$(P,f)\mapsto (P,f_i)$$
$$(P,f)\mapsto (P,\build\sum_{k\geq0}^{}f_i^k f_j)$$
and the lemma follows.\cqfd
\vskip 12pt
From now one, the coproduct in the category of spectra will be denoted by $\oplus$ and the trivial spectrum will be denoted by $0$. Actually, is $E_j$ is a family of
spectra, the spectrum $\build\oplus_j^{} E_j$ is nothing else but the filtered colimit of finite products of the $E_j$'s.

Let $X$ be an admissible set in $W(I)$. Then the map: $X\subset W(I)\build\fl_{}^\pi CW(I)$ induces a bijection $X\build\fl_{}^\sim CW_0(I)$.

For every object $(M,f)\in\Nil(A,S)$ there are only finitely many non zero morphisms $f_u$ and the product of the $\varphi_u$, for $u\in X$, induces a map:
$$\underline Nil(A,S)\fl\build\prod_{u\in X}^{} \underline Nil(A,S_u)$$
with values in $\build\oplus_{u\in X}^{} \underline Nil(A,S_u)$. Hence we have a morphism of spectra:
$$F:\underline Nil(A,S)\fl\build\oplus_{u\in X}^{} \underline Nil(A,S_u)$$
and the last thing to do is to prove that $F$ is a homotopy equivalence.

On the other hand, we have a homotopy commutative diagram:
$$\diagram{\build\lim_\longrightarrow^{}\underline Nil(A,\build\oplus_{j\in J}^{} S_j)&\hfl{}{}&\build\lim_\longrightarrow^{}
  \Bigl(\build\oplus_{u\in X\cap W(J)}^{}\underline Nil(A,S_u)\Bigr)\cr\vfl{\sim}{}&&\vfl{\sim}{}\cr\underline Nil(A,S)&\hfl{}{}&\build\oplus_{u\in X}^{}
  \underline Nil(A,S_u)\cr}$$
where the limit is taken over all finite subset $J$ of $I$. Moreover vertical arrows of this diagram are homotopy equivalences and, in order to prove the theorem, it is
enough to consider the case where $I$ is finite.

If $I$ has at most $1$ elements there is nothing to prove. So we may suppose that $I$ is a finite set with at least $2$ elements.
\vskip 12pt
If $J$ is a subset of $W(I)$ and $j$ an element of $J$, we denote by $Z(J,j)$ the set of words in $W(I)$ on the form $j^pj'$, where $j'$ is any element in $J$ distinct
from $j$ and $p$ is any non negative integer. 
\vskip 12pt
\noi{\bf 3.4 Lemma:} {\sl Suppose $I$ is finite with at least $2$ elements. Then there exists a sequence $(J_n,j_n)$, for $n\geq 0$, with 
  the following properties:

  $\bullet$ for every integer $n\geq0$, $J_n$ is a subset of $W(I)$, $j_n$ is an element of $J_n$ and $J_{n+1}$ is the set $Z(J_n,j_n)$

  $\bullet$ the set $Y=\{j_0,j_1,j_2,j_3,\dots\}\subset W(I)$ is admissible

  $\bullet$ for every integer $p>0$, there is an integer $m\geq0$ such that:
  $$\forall n>m,\ \forall u\in J_n, |u|>p$$
where $|u|$ is the length of $u$.}
\vskip 12pt
\noi{\bf Proof:} In order to define the sequence $(J_n,j_n)$, it's enough to define $J_0$ and to choose each $j_n$ in $J_n$. So we set: $J_0=I$ and, for each $n$, $j_n$
is chosen to be an element of $J_n$ of minimal length in $W(I)$.

Let $n\geq0$ be an integer. Set: $J=J_n$, $j=j_n$ and denote by $J'$ the complement of $j$ in $J$. The inclusion $J_{n+1}=Z(J,j)\subset W(I)$ factorizes through $W(J)$ and
$J_{n+1}$ can be considered as a subset of $W(J)$. Every word $u$ in $W(J)$ is written uniquely in the following form:
$$u=j^{n_0}j_1j^{n_1}j_2j^{n_2}\dots j_pj^{n_p}$$
with: $p\geq0$, $n_*\geq0$, $j_*\in J'$.

Suppose $u$ is reduced.

If $p=0$, then $u$ is a power of $j$ and, because $u$ is reduced, we have: $u=j$.

If $p>0$, then $u$ is, up to conjugation, on the form:
$$u=j^{n_0}j_1j^{n_1}j_2\dots j^{n_{p-1}}j_p$$
and $u$ is conjugate to an element in $W_0(Z(J,j))=W_0(J_{n+1})$. Moreover $u$ is reduced in $W(J)$ if and only if $u$ is reduced in $W(I)$.

Therefore every element in $CW_0(J)=CW_0(J_n)$ belongs to the image of the map $CW(J_{n+1})\fl CW(J_n)$ except the element $j=j_n$ and the inclusion
$J_{n+1}\subset J_n$ induces a bijection $\{j_n\}\coprod CW_0(J_{n+1})\build\fl_{}^\sim CW_0(J_n)$. As a consequence, we have, for each $n\geq0$, a bijection:
$$\{j_0,j_1,\dots,j_{n-1}\}\coprod CW_0(J_n)\build\fl_{}^\sim CW_0(I)$$

For each integer $n\geq0$, denote by $p_n$ the minimal length of the words in $J_n$. Since $J_{n+1}$ is contained in $W(J_n)$, we have: $p_{n+1}\geq p_n$ and the sequence
$(p_n)$ is increasing. We'll prove that this sequence is unbounded.

Let $n\geq0$ be an integer. For every $u\in J_n$, we have: $|u|\geq p_n$. For every integer $m$, denote by $H_m$ the set of elements $u\in J_m$ with $|u|=p_n$
and by $CH_m$ its image in $CW(I)$. Since $I$ is finite and every element of $H_m$ is reduced in $W(J_m)$, $CH_m$ is a finite set contained in $CW_0(I)$.

Denote by $q$ the cardinal of $CH_n$. Since $j_n$ is an element in $J_n$ of minimal length $j_n$ belongs to $H_n$. Then we have: $q>0$ and, because of the bijection
$\{j_n\}\coprod CW_0(J_{n+1})\build\fl_{}^\sim CW_0(J_n)$, the cardinal of $CH_{n+1}$ is $q-1$. For the same reason, we have the following:
$$\forall i\in\{0,1,\dots,q\}, \ \hbox{card}(CH_{n+i})=q-i$$
and then:
$$\hbox{card}(CH_{n+q})=0\ \ \hbox{and}\ \ p_{n+q}>p_n$$

Therefore the sequence $(p_n)$ is unbounded and the last property of the lemma is proven.

Since the map $\{j_0,j_1,\dots,j_{n-1}\}\coprod CW_0(J_n)\build\fl_{}^\sim CW_0(I)$ is bijective, the map $f:Y\fl CW_0(I)$ is injective. Let $u$ be an element
in $CW_0(I)$ with length $p$. Since the sequence $p_n$ is unbounded, there is an integer $n$ such that $p_n>p$ and $u$ doesn't belong to the image of
$CW_0(J_n)\fl CW_0(I)$. Hence $u$ belongs to the image of $\{j_0,j_1,\dots,j_{n-1}\}\fl CW_0(I)$ and $f$ is surjective. Therefore $Y$ is admissible and the lemma is
proven.\cqfd
\vskip 12pt
For every $n>0$ we denote by $Y_n$ the set $\{j_0,j_1,\dots,j_{n-1}\}$.
\vskip 12pt
\noi{\bf 3.5 Lemma:} {\sl For all integer $n\geq 0$, the maps $\varphi_u:\underline Nil(A,S)\fl\underline Nil(A,S_u)$, for $u\in Y_n$ and
  the map $\varphi_{J_n}:\underline Nil(A,S)\fl\underline Nil(A,S_{J_n})$ induce a homotopy equivalence of spectra:
  $$G_n:\underline Nil(A,S)\build\fl_{}^\sim\build\oplus_{u\in Y_n}^{}\underline Nil(A,S_u)\oplus\underline Nil(A,S_{J_n})$$}
\vskip 12pt
\noi{\bf Proof:} We'll prove the lemma by induction on $n$. The map $G_0$ is the identity (and then a homotopy equivalence).

Suppose $n\geq0$ and $G_n$ is a homotopy equivalence. Consider the bimodule $S'=S'_i\oplus S'_j$ with:
$$S'_i=S_{j_n}$$
$$S'_j=S_{J_n\setminus\{j_n\}}$$
By applying lemma 3.3 with this bimodule, we get a homotopy equivalence of spectra:
$$\underline Nil(A,S')=\underline Nil(A,S_{J_n})\build\fl_{}^\sim \underline Nil(A,S_{j_n})\oplus\underline Nil(A,\widehat S)$$
where $\widehat S$ is the following bimodule:
$$\widehat S=\build\oplus_{u,k}^{} \Bigl(S_{j_n}\Bigr)^kS_{u}$$
the sum being taken over all integer $k\geq0$ and all $u\not=j_n$ in $J_n$. So we have isomorphisms:
$$\widehat S\simeq\ \build\oplus_{u,k}^{} S_{j_n^k u}\simeq \build\oplus_{v\in J_{n+1}}^{} S_v=S_{J_{n+1}}$$
and then a homotopy equivalence:
$$\underline Nil(A,S_{J_n})\build\fl_{}^\sim \underline Nil(A,S_{j_n})\oplus \underline Nil(A,S_{J_{n+1}})$$
Therefore we have homotopy equivalences:
$$\underline Nil(A,S)\build\fl_{}^\sim\build\oplus_{u\in Y_n}^{}\underline Nil(A,S_u)\oplus\underline Nil(A,S_{J_n})$$
$$\build\fl_{}^\sim\build\oplus_{u\in Y_n}^{}\underline Nil(A,S_u)\oplus \underline Nil(A,S_{j_n})\oplus \underline Nil(A,S_{J_{n+1}})$$
$$\build\fl_{}^\sim\build\oplus_{u\in Y_{n+1}}^{}\underline Nil(A,S_u)\oplus \underline Nil(A,S_{J_{n+1}})$$
and $G_{n+1}$ is a homotopy equivalence. The lemma follows by induction.\cqfd
\vskip 12pt
We are now able to finish the proof of theorem 5.

As we have said before it is enough to consider the case where $I$ is finite with at least two elements.

Consider the map
$\Phi:\underline Nil(A,S)\fl\build\oplus_{u\in Y}^{}\underline Nil(A,S_u)$ induced by the functors $\varphi_u:\Nil(A,S)\fl\Nil(A,S_u)$. We have to prove that $G$ is a
homotopy equivalence and, for doing that, it will be enough to prove that $G$ induces an isomorphism:
$$\Phi_i:Nil_i(A,S)\fl\build\oplus_{u\in Y}^{}Nil_i(A,S_u)$$
for every integer $i\in\Z$.

Let $i\in \Z$ be an integer and $y$ be an element in $\build\oplus_{u\in Y}^{}Nil_i(A,S_u)$. There is an integer $n\geq0$ such that $y$ is in the direct sum of
$Nil_i(A,S_u)$, for $u\in Y_n$. Because of lemma 3.5, there is an element $x\in Nil_i(A,S)$ such that:
$$G_n(x)=y\oplus 0\in \build\oplus_{u\in Y_n}^{}Nil_i(A,S_u)\oplus Nil_i(A,S_{J_n})$$
Hence $x$ is sent to $y$ by $\Phi_i$ and $\Phi_i$ is surjective.

Suppose $i\geq 0$. Let $x\in Nil_i(A,S)$ be an element killed by $\Phi_i$. This element can be lifted in an element $y\in K_i(\Nil(A,S))=\pi_{i+1}(B(Q\Nil(A,S)))$ and
there is a finite subcategory $\E$ of the Quillen category $Q\Nil(A,S)$ such that $y$ can be lifted in an element $z\in\pi_{i+1}(B\E)$. This category $\E$ involves
only finitely many objects $(M,f)\in\Nil(A,S)$. Denote by $\F$ the set of these morphisms $f$.

Since each $f\in\F$ is nilpotent, there is an integer $p$ such that $f_u$ is trivial for each $f\in\F$ and each element $u\in W(I)$ of length $\geq p$. Because of the last
property of lemma 3.4, the set of integers $n$ such that there is some $u\in J_n$ and some $f\in\F$ with $f_u\not=0$ is finite. Hence, for $n$ big enough, the composite
functor
$$\E\fl Q\Nil(A,S)\build\fl_{}^{\varphi_{J_n}}Q\Nil(A,S_{J_n})$$ factorizes through the category $Q \P_A$ and the composite map
$$\Omega B\E\fl \Omega BQ\Nil(A,S_{J_n})\fl Nil(A,S_{J_n})$$
is trivial. Hence $z$ is killed in $Nil_i(A,S_{J_n})$ and $x\in Nil_i(A,S)$ is killed by $\varphi_{J_n}:Nil_i(A,S)\fl Nil(A,S_{J_n})$.

But $x$ is killed by $\Phi_i$ and $x$ is also killed by the map $Nil_i(A,S)\fl \build\oplus_{u\in Y_n}^{}Nil_i(A,S_u)$. Therefore $x$ is killed  by $G_n$ which is a
homotopy equivalence and $x$ is zero. Hence the morphism $\Phi_i$ is bijective for every integer $i\geq0$.

If $i$ is a negative integer, we may replace $A$ and $S$ by $\Sigma^p A$ and $\Sigma^p S$ for some $p>-i$ and the bijectivity of $\Phi_{i+p}$ for
$(\Sigma^p A,\Sigma^p S)$ implies that $\Phi_i$ (for $(A,S)$) is bijective. Hence $\Phi$ is a homotopy equivalence of spectra. Moreover this is true for the set $Y$ and
then for every admissible set in $W(I)$. Then we get the desired result and the theorem follows.\cqfd
\vskip 12pt
\noi{\bf 3.6 Proof of theorem 6.}
\vskip 12pt

Denote by $S'$ the following bimodule:
$$S'=S\oplus \build\oplus_i^{} E_i\build\otimes_{A_i}^{} F_i$$
Let $J$ be the disjoint union of $I$ and $\{0\}$. We set:
$$\forall i\in I,\ \ S_i=E_i\build\otimes_{A_i}^{} F_i=E_i\ F_i$$
$$S_0=S$$
and we have:
$$S'=\build\oplus_{j\in J}^{} S_j$$
Because of theorem 5, there exist a family of $A$-bimodules $U_k$, $k\in K$ such that the following holds:

$\bullet$ each $U_k$ is flat on both sides

$\bullet$ $\underline Nil(A,S')\simeq\underline Nil(A,S)\ \oplus\build\oplus_k^{} \underline Nil(A,U_k)$

$\bullet$ each $U_k$ has the form: $U_k=S_{j_1}S_{j_2}\dots S_{j_n}$, with $j_1,j_2,\dots,j_n$ in $J$ and $j_q\in I$ for some $q$.

Then, because of theorem 4, there exist a family of $A$-bimodules $V_k$, $k\in K$ and elements $i_k\in I$ such that

$\bullet$ each $V_k$ is flat on both sides

$\bullet$ $\underline Nil(A,S')\simeq\underline Nil(A,S)\ \oplus\build\oplus_k^{} \underline Nil(A,S_{i_k}\ V_k)$

\noi and we have, for each $k$:
$$\underline Nil(A,S_{i_k}\ V_k)=\underline Nil(A, E_{i_k} F_{i_k} V_k)\simeq \underline Nil(A_{i_k},F_{i_k} V_k\ E_{i_k})$$
But $A_{i_k}$ is regular coherent and each spectrum $\underline Nil(A_{i_k},F_{i_k} V_k\ E_{i_k})$ is contractible. The result follows.\cqfd
\vskip 12pt
\noi{\bf 3.7 Proof of lemma 3.1.} 
\vskip 12pt
With the notations of lemma 3.1, we denote by $F=B\oplus {}_b E_b\oplus \bigl({}_b E_b\bigr)^2\oplus\dots$ the tensor algebra of ${}_b E_b$ and by $\widehat E$ the
bimodule $\widehat E={}_a E_a\oplus {}_a E_b\ F\ {}_b E_a$.

For every object $(M,f)\in\Nil(A\times B,E)$, we set:
$$g_1={}_a f_a,\ g_2={}_a f_b,\ g_3={}_b f_a,\ g_4={}_b f_b$$
Since $f$ is nilpotent, there is an integer $n>0$ such that $f^n=0$. Then, for every integers $i_1,i_2,\dots,i_n$ in $\{1,2,3,4\}$ the morphism
$g_{i_1}g_{i_2}\dots g_{i_n}$ is a part of $f^n$ and we have:
$g_{i_1}g_{i_2}\dots g_{i_n}=0$. Hence the two morphisms ${}_b f_b$ and $\widehat f={}_a f_a+\build\sum_{k\geq0}^{} {}_a f_b\bigl({}_b f_b\bigr)^k{}_b f_a$ are
nilpotent and we have a functor
$$\Phi(A,B,E):\Nil(A\times B,E)\fl\Nil(B,{}_b E_b)\times\Nil(A,\widehat E)$$
sending each object $(M,f)\in\Nil(A\times B,E)$ to:
$$\Phi(A,B,E)(M,f)=\Bigl((M_b,{}_b f_b),(M_a,\widehat f)\Bigr)$$
Moreover this functor is exact.

It is easy to see that lemma 3.1 is equivalent to the fact that $\Phi(\Sigma^n A,\Sigma^n,B,\Sigma^n E)$ induces, for all $n\geq 0$, a homotopy equivalence in
K-theory. Therefore it is enough to prove that $\Phi(A,B,E)$ induces a homotopy equivalence in K-theory for every left-flat bimodule $(A\times B,E)$ such that $E_b$ is
flat on the right and that will be done by using K-theory of Waldhausen categories.

In our situation, we have three exact categories: the categories $\A$, $\B$ and $\C$ of finitely generated projective right modules over the rings $A$, $B$ and
$C=A\times B$ respectively. These categories are contained in the corresponding abelian categories $\A\v$, $\B\v$ and $\C\v$ of right modules over the corresponding rings.

If $(A,S)$ is a left-flat bimodule, we have also an exact category $\Nil(A,S)$ and a Waldhausen category $\Nil(A,S)_*$. Moreover $\Nil(A,S)$ is contained in the
abelian category $\Nil(A,S)\v$ (see part 1.3) and, as a subcategory of $\Nil(A,S)\v$, $\Nil(A,S)$ is stable under taking kernel of epimorphisms. Hence the Gillet-Waldhausen
theorem applies to the categories $\Nil(?,?)$ and we have a commutative diagram:
$$\diagram{\Nil(A\times B,E)&\hfl{\Phi(A,B,E)}{}&\Nil(B,{}_b E_b)\times\Nil(A,\widehat E)\cr\vfl{}{}&&\vfl{}{}\cr
\Nil(A\times B,E)_*&\hfl{\Phi(A,B,E)}{}&\Nil(B,{}_b E_b)_*\times\Nil(A,\widehat E)_*\cr}$$
of Waldhausen categories where the vertical functors induce homotopy equivalences in K-theory.

Therefore, in order to prove the lemma, it is enough to prove that the functor:
$$\Phi(A,B,E)_*:\Nil(A\times B,E)_*\fl\Nil(B,{}_b E_b)_*\times\Nil(A,\widehat E)_*$$
induces a homotopy equivalence in K-theory.

If $(A,S)$ is a left-flat bimodule, it is easy to see that an object of $\Nil(A,S)_*$ is
nothing else but a pair $(M,f)$, where $M$ is in $\A_*$ and $f:M\fl MS$ is a nilpotent morphism in $\A\v_*$. Then the functor $\Phi(A,B,E)_*$ is given by:
$$\Phi(A,B,E)_*(M,f)=\Bigl((M_b,{}_b f_b),(M_a,\widehat f)\Bigr)$$
for every $M\in\C_*=\A_*\times\B_*$.

Consider the Waldhausen category $\E=\Nil(A\times B,E)_*$. We may define a new subcategory of equivalences by saying that $\varphi:(M,f)\fl (M',f')$ is an equivalence
if the induced morphism $M_a\fl M'_a$ is an isomorphism in homology. With this new equivalences, we get a new Waldhausen category $\E'$. We have also a Waldhausen
subcategory $\E_0$ of $\E$ generated by the objects $(M,f)\in\E$ such that $M_a$ is acyclic. Denote also by $\Nil_0$ the Waldhausen subcategory of $\Nil(A,\widehat E)_*$
generated by pairs $(M,f)$ with $M$ acyclic. Hence we have a commutative diagram of essentially small Waldhausen categories:
$$\diagram{\E_0&\hfl{}{}&\E&\hfl{}{}&\E'\cr\vfl{\Phi_0}{}&&\vfl{\Phi(A,B,E)}{}&&\vfl{\Phi'}{}\cr
  \Nil(B,{}_b E_b)_*\times\Nil_0&\hfl{}{}&\Nil(B,{}_b E_b)_*\times\Nil(A,\widehat E)_*&\hfl{}{}&\Nil(A,\widehat E)_*\cr}$$

Since every morphism in $\Nil_0$ is an equivalence, the category $\Nil_0$ has trivial K-theory and, because of the fibration theorem, the two lines of the diagram induce
fibrations in K-theory. Hence in order to prove the lemma, it's enough to prove that $\Phi_0$ and $\Phi'$ induce homotopy equivalences on K-theory and that's equivalent to
show that $\Phi'$ and the functor:
$$\Phi'_0:\E_0\build\fl_{}^{\Phi_0}\Nil(B,{}_b E_b)_*\times\Nil_0\build\fl_{}^{pr_1}\Nil(B,{}_b E_b)_*$$
induce homotopy equivalences on K-theory. Moreover these two functors have the approximation properties (App1).
\vskip 12pt
\noi{\bf 3.8 Lemma:} {\sl The functor $\Phi'_0$ induces a homotopy equivalence in K-theory.}
\vskip 12pt
\noi{\bf Proof:} Because of the approximation theorem 2.3, we just have to prove that $\Phi'_0$ has the property (App2).

Let $X=(M,f)\in\E_0$ and $Y=(Q,g)\in\Nil(B,{}_b E_b)_*$ be two objects in the corresponding categories. A morphism $\varphi:\Phi'_0(X)\fl Y$ is represented by a morphism
$\varphi:M_b\fl Q$ in $\B_*$ making the following diagram commutative:
$$\diagram{M_b&\hfl{\varphi}{}&Q\cr\vfl{{}_b f_b}{}&&\vfl{g}{}\cr M_b\ {}_b E_b&\hfl{\varphi}{}&Q\ {}_b E_b\cr}$$
Denote by $C$ the cylinder of $\varphi$ and by $C'$ its mapping cone (or its 0-cone as defined in 2.1). The morphisms ${}_b f_b$ and $g$ induce two nilpotent morphisms
$\lambda:C\fl C\ {}_b E_b$ and $\lambda': C'\fl C'\ {}_b E_b$. By naturality, we have a commutative diagram in $\C\v_*$:
$$\diagram{M_b&\hfl{i}{}&C&\hfl{p}{}&Q\cr\vfl{{}_b f_b}{}&&\vfl{\lambda}{}&&\vfl{g}{}\cr M_b\ {}_b E_b&\hfl{i}{}&C\ {}_b E_b&\hfl{p}{}&Q\ {}_b E_b\cr}$$
where $i: M_b\fl C$ is a cofibration and $p:C\fl Q$ a homotopy equivalence. Moreover we have: $\varphi=pi$.

Define the object $M'\in \C_*$ by: $M'_a=M_a$ and $M'_b=C$. In order to define the desired object $X'\in\E_0$, we have to construct the morphism $f':M'\fl M'E$.

We set: ${}_a f'_a={}_a f_a$, ${}_b f'_b=\lambda$ and ${}_b f'_a=i\ {}_b f_a$. Since $i:M_b\fl C$ is a cofibration and $M_a\ {}_a E_b$ is acyclic, there is no
obstruction to extend ${}_a f_b: M_b\fl M_a\ {}_a E_b$ to a morphism: ${}_a f'_b: C\fl M_a\ {}_a E_b$. Hence we get a morphism $f':M'\fl M'E$.

Let $\overline C$ be the finite complex in $\C_*$ defined by: $\overline C_a=0$ and $\overline C_b=C'$ and $\overline\lambda:\overline C\fl\overline C E$ be the morphism
in $\C\v_*$ associated with $\lambda'$ under the canonical bijection: $\Hom_{\C\v_*}(\overline C,\overline C E)\simeq \Hom_{\B\v_*}(C',C'{}_b E_b)$. Since $\lambda'$ is
nilpotent, $\overline\lambda$ is also nilpotent.

We have a commutative diagram in $\C\v_*$ with exact rows:
$$\diagram{0&\hfl{}{}&M&\hfl{}{}&M'&\hfl{}{}&\overline C&\hfl{}{}&0\cr &&\vfl{f}{}&&\vfl{f'}{}&&\vfl{\overline\lambda}{}&&\cr
  0&\hfl{}{}&ME&\hfl{}{}&M'E&\hfl{}{}&\overline C E&\hfl{}{}&0\cr}$$
and, since $f$ and $\overline\lambda$ are nilpotent, $f'$ is nilpotent also. More precisely, if we have: $f^p=0$ and $\overline\lambda^q=0$, then we have $f'^{p+q}=0$.

Since $M'_a$ is acyclic, $X'=(M',f')$ is an object in $\E_0$ and the morphism $i$ induces a morphism $\alpha:X\fl X'$. Moreover we have the following commutative diagram:
$$\diagram{\Phi'_0(X)&\hfl{\varphi}{}&Y\cr\vfl{\alpha}{}&&\vfl{=}{}\cr\Phi'_0(X')&\hfl{p}{}&Y\cr}$$
and the property (App2) is satisfied. The lemma follows.\cqfd
\vskip 12pt
So the last thing to do is to prove that the functor $\Phi':\E'\fl\Nil(A,\widehat E)_*)$ induces a homotopy equivalence in K-theory. As before, it's enough to prove that
$\Phi'$ has the property (App2). In order
to do that we'll need two technical results:
\vskip 12pt
\noi{\bf 3.9 Lemma:} {\sl Let $A$ and $B$ be two rings and $S$ be an $(A,B)$-bimodule. Suppose $S$ is flat on the left. Let $X$ be a finite $B$-complex, $Y$ be an
$A$-complex and $f:X\fl YS$ be a morphism of $B$-complexes. Then there exist a finite $A$-complex $Y'$, a morphism $g:Y'\fl Y$ and a commutative diagram:
$$\diagram{X&\hfl{}{}&Y'S\cr\vfl{=}{}&&\vfl{g}{}\cr X&\hfl{f}{}&YS}$$}
\vskip 12pt
\noi{\bf 3.10 Lemma:} {\sl Let $A$ be a ring, $X$ and $Y$ two $A$-complexes and $f:X\fl Y$ be a morphism. Suppose the module $\build\oplus_n^{} X_n$ is finitely
  presented and each $Y_n$ is flat. Then $f$ factorizes through a finite $A$-complex.}
\vskip 12pt
These two lemmas will be proven at the end of the section.
\vskip 12pt
As a consequence of these two lemmas we have the following result:
\vskip 12pt
\noi{\bf 3.11 Lemma:} {\sl Let $A$, $B$ and $C$ be three rings, $S$ be a $(C,B)$-bimodule and T be a $(B,A)$-bimodule. Suppose $S$ is flat on the right and $T$ is flat
on the left. Let $X$, $Y$ and $Z$ be finite complexes over the rings $A$, $B$ and $C$ respectively and $f:X\fl YT$ and $g:Y\fl ZS$ be two morphisms. We suppose that the
composite morphism:
$$X\build\fl_{}^f YT\build\fl_{}^g ZST$$
is zero. Then the morphism $g:Y\fl ZS$ factorizes through a finite complex $Y'$ by a morphism $\lambda:Y\fl Y'$ in such a way that the  composite morphism:
$$X\build\fl_{}^f YT\build\fl_{}^\lambda Y'T$$
is zero.}
\vskip 12pt
\noi{\bf Proof:} Let $K$ be the kernel of $g$. Since $T$ is flat on the left we have two exact sequences:
$$0\fl K\build\fl_{}^i Y\build\fl_{}^g ZS$$
$$0\fl KT\build\fl_{}^i YT\build\fl_{}^g ZST$$
and the morphism $f$ factorizes through $KT$. Because of lemma 3.9 there exist a finite $B$-complex $K'$, a morphism $j:K'\fl K$ and a commutative diagram:
$$\diagram{X&\hfl{f'}{}&K'T\cr\vfl{=}{}&&\vfl{ij}{}\cr X&\hfl{f}{}&YT\cr}$$

Denote by $Y_1$ the cokernel of $ij:K'\fl Y$. Since $gi$ is the zero morphism, there is a morphism $g_1:Y_1\fl ZS$ making the following diagram commutative:
$$\diagram{Y&\hfl{g}{}&ZS\cr\vfl{}{}&&\vfl{=}{}\cr Y_1&\hfl{g_1}{}&ZS\cr}$$
Since $S$ is flat on the right, $ZS$ is flat and, because of lemma 3.10, $g_1$ factorizes through a finite complex $Y'$. The lemma follows.\cqfd 
\vskip 12pt
\noi{\bf 3.12 Lemma:} {\sl The functor $\Phi':\E'\fl\Nil(A,\widehat E)_*)$ has the property (App2).}
\vskip 12pt
\noi{\bf Proof:} Let's set:
$$H={}_a E_a\hskip 24pt K={}_b E_b\hskip 24pt S={}_a E_b\hskip 24pt T={}_b E_a\hskip 24pt \widehat K=F=B\oplus K\oplus K^2\oplus \dots$$
Let $X$ and $Y$ be two objects in $\Nil(A\times B,E)_*$ and $\Nil(A,\widehat E)_*$ respectively and $\varphi$ be a morphism from $\Phi'(X)$ to $Y$. By setting:
$(U,f)=X$, $(M,g)=Y$, $P=U_a$, $Q=U_b$, we see that $P$ and $M$ are finite $A$-complexes, $Q$ is a finite $B$-complex, $\varphi$ is a morphism from $P$ to $M$ and
$g$ is written as a finite sum:
$$g=\lambda+\build\sum_{k\geq0}^{}\mu_k$$
where $\lambda$ is a morphism from $M$ to $MH$ and each $\mu_k$ is a morphism from $M$ to $MSK^k T$. Moreover the following diagrams are commutative:
$$\diagram{P&\hfl{\varphi}{}&M\cr\vfl{{}_a f_a}{}&&\vfl{\lambda}{}\cr PH&\hfl{\varphi}{}&MH\cr}\hskip 64pt
\diagram{P&\hfl{\varphi}{}&M\cr\vfl{{}_a f_b\bigl({}_b f_b\bigr)^k{}_b f_a}{}&&\vfl{\mu_k}{}\cr PSK^k T&\hfl{\varphi}{}&MSK^k T\cr}$$

We want to construct an object $X'\in\Nil(A\times B,E)_*$ and a morphism $\alpha:X\fl X'$ such that $\Phi'(\alpha):\Phi'(X)\fl \Phi'(X')$ is isomorphic to 
$\varphi:\Phi'(X)\fl Y$. So we want to have: $X'=(U',f')$ and $U'_a=M$ and, to determine $X'$, we need to define $Q'=U'_b$, the morphism $f'$ and the morphism 
$\alpha:Q\fl Q'$. Actually, $Q'$ will be constructed as a direct sum: $Q'=\build\oplus_{i\geq0}^{} Q_i$ such that ${}_b f'_b(Q_{i+1})$ is contained in $Q_iK$ for all 
$i\geq 0$ and is zero for $i=-1$ and ${}_a f'_b$ vanishes on each $Q_i$, for $i>0$.

So we need to define finite $B$-complexes $Q_i$, a morphism $e:Q_0\fl MS$ and morphisms $\alpha_i:Q\fl Q_i$, $\theta_i:Q_{i+1}\fl Q_iK$ and $\beta_i: M\fl Q_i T$ (with
$Q_i=0$ for $i$ big enough) and the morphism $f'$ is defined by:
$${}_a f'_a=\lambda\hskip 24pt {}_b f'_b=\build\sum_{i\geq0}^{}\theta_i\hskip 24pt {}_b f'_a=\build\sum_{i\geq0}^{}\beta_i\hskip 24pt {}_a f'_b=e\ pr_0$$
where $pr_0:Q'\fl Q_0$ is the projection. The morphism $\alpha$ is equal to $\varphi:P\fl M$ on $P$ and to $\build\sum_{i\geq0}^{} \alpha_i:Q\fl Q'$ on $Q$.

But we have two conditions: the fact that $\alpha$ is a morphism and the equality: $\Phi'(X')=Y$. These conditions are equivalent to:

$\varphi\ {}_a f_b=e\alpha_0$

$\alpha_i\ {}_b f_a=\beta_i\varphi$

$\theta_i\alpha_{i+1}=\alpha_i\ {}_b f_b$

$\mu_i=e\theta_0\theta_1\dots\theta_{i-1}\beta_i$

\noi for all $i\geq0$.

For technical reasons, we introduce the morphism $e_i=e\theta_0\theta_1\dots\theta_{i-1}$ from $Q_i$ to $MSK^i$. Then we have to construct, for
each $i\geq0$, the complex $Q_i$ and morphisms $e_i:Q_i\fl MSK^i$, $\alpha_i:Q\fl Q_i$, $\beta_i: M\fl Q_i T$ and $\theta_i:Q_{i+1}\fl Q_iK$, with the following
properties:

$A(i)$: $\varphi\ {}_a f_b\bigl({}_b f_b\bigr)^i=e_i\alpha_i$

$B(i)$: $\alpha_i\ {}_b f_a=\beta_i\varphi$

$C(i)$: $\theta_i\alpha_{i+1}=\alpha_i\ {}_b f_b$

$D(i)$: $\mu_i=e_i\beta_i$

$E(i)$: $e_{i+1}=e_i\theta_i$

\noi for all $i\geq0$.

Since the sum of the $\mu_k$'s is finite, there is an integer $n>0$ such that: $\mu_i=0$ for all $i>n$.

So we'll construct $(Q_i, e_i, \alpha_i, \beta_i, \theta_i)$ by induction. Let $i\geq0$ be an integer and suppose that $(Q_j, e_j, \alpha_j, \beta_j, \theta_j)$ is
defined for all $j>i$ such that the properties $A(j)$, $B(j)$, $C(j)$, $D(j)$, $E(j)$, are satisfied for all $j>i$. We begin this induction with $i=n$ by
setting: $Q_j=0$ for all $j>n$.

We have to construct $Q_i$ and morphisms $e_i$, $\alpha_i$, $\beta_i$ and $\theta_i$.

Consider the following morphisms:
$$\varphi\ {}_a f_b\bigl({}_b f_b\bigr)^i:Q\fl MSK^i$$
$$\mu_i:M\fl MSK^i T$$
$$e_{i+1}:Q_{i+1}\fl MSK^{i+1}$$
These morphisms induce a morphism $h:Q\oplus M\oplus Q_{i+1}\fl MSK^i(B\oplus T\oplus K)$ and, because of lemma 3.9, there are a finite complex $Q_i$, a morphism
$e_i:Q_i\fl MSK^i$ and a commutative diagram:
$$\diagram{Q\oplus M\oplus Q_{i+1}&\hfl{h'}{}&Q_i(B\oplus T\oplus K)\cr\vfl{=}{}&&\vfl{e_i}{}\cr Q\oplus M\oplus Q_{i+1}&\hfl{h}{}&MSK^i(B\oplus T\oplus K)\cr}$$
The morphism $h'$ induces morphisms:
$$\alpha_i:Q\fl Q_i$$
$$\beta_i:M\fl Q_i T$$
$$\theta_i:Q_{i+1}\fl Q_i K$$
and properties $A(i)$, $D(i)$, $(E(i)$ are satisfied. Denote by $u$ and $v$ the defaults of properties $B(i)$ and $C(i)$:
$$u=\beta_i\varphi-\alpha_i\ {}_b f_a$$
$$v=\alpha_i\ {}_b f_b-\theta_i\alpha_{i+1}$$
Because of properties $A(i)$ and $D(i)$, we have:
$$e_i u=e_i\beta_i\varphi-e_i\alpha_i\ {}_b f_a=\mu_i\varphi-\varphi\ {}_a f_b\bigl({}_b f_b\bigr)^i{}_b f_a=0$$
and, because of properties $A(i)$, $E(i)$ and $A(i+1)$, we have:
$$e_i v=e_i\alpha_i\ {}_b f_b-e_i\theta_i\alpha_{i+1}=\varphi\ {}_a f_b\bigl({}_b f_b\bigr)^{i+1}-e_{i+1}\alpha_{i+1}=0$$
Since a tensor product of bimodules which are flat on the right is flat on the right, we can apply the lemma 3.11 to
morphisms $u\oplus v: P\fl Q_i(T\oplus K)$ and $e_i:Q_i\fl MSK^i$. Thus $e_i$ factorizes through a finite complex $Q'_i$ by a morphism $\e:Q_i\fl Q'_i$ such
that: $\e(u\oplus v)=0$.

Hence, up to replacing $Q_i$ by $Q'_i$, we may as well suppose that $A(i)$, $B(i)$, $C(i)$, $D(i)$, $E(i)$ are satisfied.
Therefore $Q_i$, $e_i$, $\alpha_i$, $\beta_i$, $\theta_i$ are defined and $A(i)$, $B(i)$, $C(i)$, $D(i)$, $E
(i)$ are satisfied for all $i\geq0$.

Then the finite complex $Q'=\build\oplus_{i\geq0}^{} Q_i$ is constructed and the morphism $f'$ is defined by:
$${}_a f'_a=\lambda$$
$${}_b f'_b=\build\oplus_{i\geq0}^{}\theta_i$$
$${}_a f'_b=e_0pr_0$$
$${}_b f'_a=\build\oplus_{i\geq0}^{}\beta_i$$
Hence the desired object $X'$ is constructed and $\Phi'$ has the approximation property (App2). The lemma follows and then follow lemma 3.1 and theorems 4, 5 and 6.\cqfd
\vskip 12pt
\noi{\bf 3.13 Proof of lemma 3.9.}
\vskip 12pt
The situation is the following: $(A,S)$ is a left-flat bimodule, $X$ is a finite $A$-complex, $Y$ is an $A$-complex and $f:X\fl YS$ is a morphism of $A$-complexes. We want
to construct a finite $A$-complex $Y'$ and a morphism $Y'\fl Y$ such that $f:X\fl YS$ factorizes through $Y'S$. 

For each integer $n$, denote by $X(n)$ the $n$-skeleton of $X$. Let $n$ be an integer. Suppose the $n$-skeleton $Y'(n)$ of $Y'$ is constructed in such a way that we have
a morphism $g_n:Y'(n)\fl Y$ and a commutative diagram:
$$\diagram{X(n)&\hfl{h_n}{}&Y'(n)S\cr\vfl{=}{}&&\vfl{g_n}{}\cr X(n)&\hfl{f}{}&YS\cr}\leqno{(D_n)}$$
If $n$ is small enough, $X(n)$ is null and $Y'(n)$ can be chosen to be zero.

Denote by $Z_n$ the kernel of the morphism $d:Y'_n\fl Y'_{n-1}$ and by $U_{n+1}$ the module defined by the cartesian square:
$$\diagram{U_{n+1}&\hfl{\alpha}{}&Y_{n+1}\cr \vfl{\beta}{}&&\vfl{d}{}\cr Z_n&\hfl{g_n}{}&Y_b\cr}$$

The composite morphism $X_{n+1}\build\fl_{}^d X_n\build\fl_{}^{h_n}Y'_n S$ takes values in $Z_n S$ and induces, together with the morphism $f:X_{n+1}\fl Y_{n+1}S$,
a well defined morphism $\lambda:X_{n+1}\fl U_{n+1}$. So we get a commutative diagram:
$$\diagram{X_{n+1}&\hfl{\lambda}{}&U_{n+1}S&\hfl{\alpha}{}&Y_{n+1}S\cr\vfl{d}{}&&\vfl{\gamma}{}&&\vfl{d}{}\cr X_n&\hfl{h_n}{}&Y'_n S&\hfl{g_n}{}&Y_n S\cr}$$
where $\gamma$ is the composite morphism $U_{n+1}\build\fl_{}^\beta Z_n\subset Y'_n$.

Since $X_{n+1}$ is finitely generated, there is a finitely generated submodule $M$ in $U_{n+1}$ such that $\lambda(X_{n+1})$ is contained in $MS$. Let $Y'_{n+1}$ be
a finitely generated projective $A$-module and $\mu:Y'_{n+1}\fl M$ be an epimorphism. Since $X_{n+1}$ is projective, the morphism $\lambda:X_{n+1}\fl Y'S$ can be lifted
in a morphism $X_{n+1}\fl Y'_{n+1}S$ and we get a commutative diagram:
$$\diagram{X_{n+1}&\hfl{h}{}&Y'_{n+1}S&\hfl{g}{}&Y_{n+1}\cr\vfl{d}{}&&\vfl{d}{}&&\vfl{d}{}\cr X_n&\hfl{}{}&Y'_nS&\hfl{}{}&Y_nS\cr}$$
where $d:Y'_{n+1}\fl Y'_n$ is the morphism $\gamma\mu$ and $g:Y'_{n+1}\fl Y_{n+1}$ is the morphism $\alpha\mu$.

Thus we have constructed the complex $Y'(n+1)$ and the commutative diagram $(D_{n+1})$. By induction we have $Y'(n)$ and the commutative diagram $(D_n)$ for every integer
$n$ and, for $n$ big enough, $Y'(n)$ and the diagram $(D_n)$ is a solution of the problem.\cqfd 
\vskip 12pt 
\noi{\bf 3.14 Proof of lemma 3.10.} 
\vskip 12pt
In the lemma, $f:X\fl Y$ is a morphism between two $A$-complexes, the direct sum of the $X_n's$ if finitely presented and each $Y_n$ is flat.

For every $A$-complex $E$, denote by $E^n$ its $n$-coskeleton i.e. the quotient of $E$ by its $(n-1)$--skeleton. 

Let $n$ be an integer. Suppose that the $n$-coskeleton $F^n$ of a finite complex $F$ is constructed in such a way that the morphism $f:X^n\fl Y^n$ induced by $f:X\fl Y$
factorizes through $F^n$ via two morphisms $\alpha:X^n\fl F^n$ and $\beta:F^n\fl Y^n$. If $n$ is big enough $X^n$ is trivial and we may set: $F^n=0$.

We have a commutative diagram:
$$\diagram{X_{n+1}&\hfl{\alpha}{}&F_{n+1}&\hfl{\beta}{}&Y_{n+1}\cr\vfl{d}{}&&\vfl{d}{}&&\vfl{d}{}\cr X_n&\hfl{\alpha}{}&F_n&\hfl{\beta}{}&Y_n\cr}$$

Let $E$ be the $A$-module defined by the cocartesian square:
$$\diagram{X_n&\hfl{}{}&F'_n\cr\vfl{d}{}&&\vfl{}{}\cr X_{n-1}&\hfl{}{}&E\cr}$$
where $F'_n$ is the cokernel of the morphism $d:F_{n+1}\fl F_n$. Since $X_n$, $X_{n-1}$ and $F'_n$ are finitely presented, the $A$-module $E$ is also finitely presented.

We have a commutative diagram:
$$\diagram{X_{n+1}&\hfl{\alpha}{}&F_{n+1}&\hfl{\beta}{}&Y_{n+1}\cr\vfl{d}{}&&\vfl{d}{}&&\vfl{d}{}\cr X_n&\hfl{\alpha}{}&F_n&\hfl{\beta}{}&Y_n\cr
  \vfl{d}{}&&\vfl{\delta}{}&&\vfl{d}{}\cr X_{n-1}&\hfl{}{}&E&\hfl{}{}&Y_{n-1}\cr}$$
where the composite morphism $F_{n+1}\build\fl_{}^d F_n\fl E$ is trivial.

But $E$ is finitely presented and $Y_{n-1}$ is flat. Therefore the morphism $E\fl Y_{n-1}$ factorizes through a finitely generated projective $A$-module $F_{n-1}$
and, together with the composite morphism $F_n\build\fl_{}^\delta E\fl F_{n-1}$, we get the desired finite complex $F^{n-1}$ and a commutative diagram:
$$\diagram{X^{n-1}&\hfl{\alpha}{}&F^{n-1}&\hfl{\beta}{}&Y^{n-1}\cr\vfl{}{}&&\vfl{}{}&&\vfl{}{}\cr X^n&\hfl{\alpha}{}&F^n&\hfl{\beta}{}&Y^n\cr}$$

So we construct the complexes $F^n$ inductively and, for $n$ small enough, the morphism $f:X\fl Y$ factorizes through the finite complex $F=F^n$.\cqfd
\vskip 12pt
\noi{\bf Remark:} It is not clear that the right flatness condition is necessary in theorems 4, 5 and 6. Actually this condition is only used in order to prove that the
functor $\Phi':\E'\fl \Nil(A,\widehat E)_*$ (in the proof of lemma 3.1) is a homotopy equivalence. The proof given here needs the right flatness condition (in the lemma 3
.10) but another proof without this condition is still possible. 
\vskip 24pt
\noi{\bf 4 Whitehead spectra.}
\vskip 12pt
If $E$ is an $\Omega$-spectrum and $X$ is a space, we denote by $H(X,E)$ the $\Omega$-spectrum associated to the smash product $X\vc E$. For every $i\in\Z$, we have:
$$\pi_i(H(X,E))\simeq H_i(X,E)$$

In [Wa1], section 15, Waldhausen associates to any ring $R$ and any group $G$ an assembly map: $H(BG,K(R))\fl K(R[G])$ which is a map of infinite loop spaces. This
assembly map induces assembly maps $H(BG,\Sigma^nK(R))\fl K(\Sigma^n R[G])$ and then an assembly map $h:H(BG,\underline K(R))\fl\underline K(R[G])$ which is a map of
spectra. So we get a fibration of spectra:
$$H(BG,\underline K(R))\build\fl_{}^h\underline K(R[G])\fl \underline Wh^R(G)$$
The spectrum $\underline Wh^R(G)$ is called the Whitehead spectrum of $G$ relative to $R$.

For every integer $i$, we set:
$$Wh_i(G)=\pi_i(\underline Wh^\Z(G))$$
For $i<0$, the group $Wh_i(G)$ is isomorphic to $K_i(Z[G])$ and we have exact sequences:
$$0\fl \Z\fl K_0(Z[G])\fl Wh_0(G)\fl 0$$
$$0\fl \Z/2\oplus H_1(G,\Z)\fl K_1(\Z[G])\fl Wh_1(G)\fl 0$$\vskip 12pt

More precisely, $Wh_0(G)$ is the reduced $K_0$-group of $\Z[G]$, $Wh_1(G)$ is the classical Whitehead group of $G$ and $Wh_2(G)$ is the second Whitehead group of $G$ as
defined in [HW].

Following Waldhausen, a group $G$ is said to be regular noetherian (resp. regular coherent) if, for every ring $R$ which is regular noetherian on the right,
$R[G]$ is regular noetherian (resp. regular coherent) on the right. Since $R^{op}[G]$ is isomorphic to $(R[G])^{op}$, this condition is equivalent to the condition
obtained by replacing right by left. We denote also by $\G$ the category of groups and monomorphisms of groups.

\vskip 12pt
\noi{\bf 4.1 Proposition:} {\sl We have the following properties:

  $\bullet$ If $G$ is the amalgamated free product of a diagram in $\G$:
  $$\diagram{H&\hfl{}{}&G_1\cr\vfl{}{}&&\cr G_2&&\cr}$$
  where $G_1$ and $G_2$ are regular coherent and $H$ regular noetherian, then $G$ is regular coherent.

  $\bullet$ If $G$ is the HNN extension of a diagram in $\G$:
  $$H\ \ \bhfl{\alpha}{\beta}\ \ G_1$$
  where $G_1$ is regular coherent and $H$ is regular noetherian, then $G$ is regular coherent.

  $\bullet$ If $G$ is the colimit of a filtered system $G_i$ in $\G$, where each $G_i$ is regular coherent, then $G$ is regular coherent.

  $\bullet$ A subgroup of a regular coherent group is regular coherent.}
\vskip 12pt
\noi{\bf Proof:} All these properties are proven in [Wa1] (in theorem 19.1) except the third one.

Let $G$ be the colimit of a filtered system $G_i$ in $\G$ and $R$ be a ring which is regular noetherian on the right. Suppose each $G_i$ is regular coherent.
Set: $A=R[G]$ and $A_i=R[G_i]$.

Each ring $A_i$ is regular coherent on the right and for each $i\in I$, the ring $A$ is free on the left over $A_i$.

Let $M$ be a finitely presented right $A$-module. We have an exact sequence of right $A$-modules:
$$F_1\build\fl_{}^f F_0\fl M\fl 0$$
where $F_0$ and $F_1$ are finitely generated free $A$-modules. The morphism $f$ is represented by a finite matrix with entries in $A$. Since $A$ is the colimit of the
$A_i$'s, there is an element $i\in I$ such that $A_i$ contains all the entries of $f$. Therefore $f$ comes from a finite matrix with entries in $A_i$ and there exist a
finitely presented right $A_i$-module $M'$ and an isomorphism $M'\otimes_{A_i}A\simeq M$.

Since $A_i$ is regular coherent on the right we have an exact sequence of right $A_i$-modules:
$$0\fl C_n\fl C_{n-1}\fl\dots\fl C_0\fl M'\fl 0$$
where each $C_k$ is finitely generated projective. Since $A$ is free on the left over $A_i$, we have an exact sequence of right $A$-modules:
$$0\fl C_n\build\otimes_{A_i}^{}A\fl C_{n-1}\build\otimes_{A_i}^{}A\fl\dots\fl C_0\build\otimes_{A_i}^{}A\fl M\fl 0$$
But each $C_k\build\otimes_{A_i}^{}A$ is finitely generated projective right $A$-module. Therefore every finitely presented $A$-module has a finite resolution by
finitely generated projective $A$-modules and $A$ is regular coherent. The result follows.\cqfd
\vskip 12pt

Let Cl be the class of groups defined by Waldhausen in [Wa1]. This class is the smallest class of groups satisfying the following:

$\bullet$ The trivial group belongs to Cl.

$\bullet$ If $G$ is the amalgamated free product of a diagram in $\G$:
$$\diagram{H&\hfl{}{}&G_1\cr\vfl{}{}&&\cr G_2&&\cr}$$
where $G_1$ and $G_2$ are in Cl and $H$ regular coherent, then $G$ belongs to Cl.

$\bullet$ If $G$ is the HNN extension of a diagram in $\G$:
  $$H\ \ \bhfl{\alpha}{\beta}\ \ G_1$$
where $G_1$ is in Cl and $H$ is regular coherent, then $G$ belongs to Cl.

$\bullet$ If $G$ is the colimit of a filtered system $G_i$ in $\G$, where each $G_i$ is in Cl, then $G$ belongs to Cl.
\vskip 12pt
This class contains free groups, torsion free abelian groups, poly-$\Z$-groups, torsion free one-relator groups and fundamental groups of many low-dimensional manifolds.
It is also closed under taking subgroups. See theorem 19.5 in [Wa1].
\vskip 12pt
\noi{\bf 4.2 Theorem:} {\sl For every group $G$ in Cl and every ring $R$ which is regular noetherian on the right, the Whitehead spectrum $\underline Wh^R(G)$ is
contractible.}
\vskip 12pt
\noi{\bf Proof:} This is essentially theorem 19.4 in [Wa1]. We just have to replace spaces $Nil(A,S)$ by spectra $\underline Nil(A,S)$. Since all these spectra are
contractible, the result follows.\cqfd
\vskip 12pt
We'll construct a class of groups Cl$_1$ obtained by replacing the condition ''$H$ is regular coherent'' (in the definition of Cl) by a weaker condition in such a way that
theorem 4.2 is still true for groups in Cl$_1$.
\vskip 12pt
Consider a diagram of groups:
$$\diagram{H&\hfl{\alpha}{}&G_1\cr\vfl{\beta}{}&&\cr G_2&&\cr}$$
We say that this diagram is regular coherent if the following holds:

$\bullet$ $\alpha$ and $\beta$ are monomorphisms

$\bullet$ for every $x\in G_1\setminus\alpha(H)$ and every $y\in G_2\setminus\beta(H)$, the intersection of the two groups $\alpha^{-1}(x\alpha(H)x^{-1})$ and
$\beta^{-1}(y\beta(H)y^{-1})$ is regular coherent.

Consider a diagram of groups:
  $$H\ \ \bhfl{\alpha}{\beta}\ \ G_1$$
We say that this diagram is regular coherent if the following holds:

$\bullet$ $\alpha$ and $\beta$ are monomorphisms

$\bullet$ for every $x\in G_1\setminus\alpha(H)$ and every $y\in G_1\setminus\beta(H)$, the intersection of the two groups $\alpha^{-1}(x\alpha(H)x^{-1})$ and
$\beta^{-1}(y\beta(H)y^{-1})$ is regular coherent.

$\bullet$ for every $x\in G_1$, the group $\beta^{-1}(x\alpha(H)x^{-1})$ is regular coherent.

Since the condition ''regular coherent'' is stable under taking subgroups, it is easy to see that diagrams above are regular coherent if the subgroup $H$ is regular
coherent.
\vskip 12pt
So we define the class Cl$_1$ as the smallest class of groups satisfying the following:

$\bullet$ The trivial group belongs to Cl$_1$.

$\bullet$ If $G$ is the amalgamated free product of a diagram $D$:
$$\diagram{H&\hfl{}{}&G_1\cr\vfl{}{}&&\cr G_2&&\cr}$$
where $G_1$ and $G_2$ are in Cl$_1$ and $D$ is regular coherent, then $G$ belongs to Cl$_1$.

$\bullet$ If $G$ is the HNN extension of a diagram $D'$:
  $$H\ \ \bhfl{\alpha}{\beta}\ \ G_1$$
where $G_1$ is in Cl$_1$ and $D'$ is regular coherent, then $G$ belongs to Cl$_1$.

$\bullet$ If $G$ is the colimit of a filtered system $G_i$ in $\G$, where each $G_i$ is in Cl$_1$, then $G$ belongs to Cl$_1$.
\vskip 12pt
\noi{\bf 4.3 Theorem:} {\sl Let $R$ be a ring which is regular noetherian on the right and $G$ be the amalgamated free product of a regular coherent diagram of groups:
  $$\diagram{H&\hfl{\alpha}{}&G_1\cr\vfl{\beta}{}&&\cr G_2&&\cr}$$
  Then this diagram induces a homotopically cartesian diagram of spectra:
  $$\diagram{\underline Wh^R(H)&\hfl{\alpha}{}&\underline Wh^R(G_1)\cr\vfl{\beta}{}&&\vfl{}{}\cr \underline Wh^R(G_2)&\hfl{}{}&\underline Wh^R(G)\cr}$$}
\vskip 12pt
\noi{\bf 4.4 Theorem:} {\sl Let $R$ be a ring which is regular noetherian on the right and $G$ be the HNN extension of a regular coherent diagram of groups:
  $$H\ \ \bhfl{\alpha}{\beta}\ \ G_1$$
  Then this diagram induces a homotopy fibration of spectra:
  $$\underline Wh^R(H)\build\fl_{}^f\underline Wh^R(G_1)\fl \underline Wh^R(G)$$
where $f$ is the difference (in $\Omega$sp) of maps induced by $\alpha$ and $\beta$. }
\vskip 12pt
\noi{\bf Proofs of theorems 4.3 and 4.4:} In the amalgamated case, we have a commutative diagram of spectra:
$$\diagram{H(BH,\underline K(R))&\hfl{\alpha\oplus -\beta}{}&H(BG_1,\underline K(R))\oplus H(BG_2,\underline K(R))&\hfl{}{}&H(BG,\underline K(R))\cr
  \vfl{}{}&&\vfl{}{}&&\vfl{}{}\cr
  \underline K(R[H])&\hfl{\alpha\oplus -\beta}{}&\underline K(R[G_1])\oplus \underline K(R[G_2])&\hfl{}{}&\underline K(R[G])'\cr}$$
where horizontal lines are fibrations and vertical maps are assembly maps. Moreover theorem 2 implies a homotopy equivalence:
$$\underline K(R[G])\simeq \underline K(R[G])'\oplus\Omega^{-1}\underline Nil(R[H]\times R[H],S)$$
for some $R[H]\times R[H]$ bimodule $S$. Therefore we have a fibration:
$$\underline Wh^R(H)\fl\underline Wh^R(G_1)\oplus\underline Wh^R(G_2)\fl\underline Wh^R(G)'$$
and a homotopy equivalence:
$$\underline Wh^R(G)\simeq\underline Wh^R(G)'\oplus\Omega^{-1}\underline Nil(R[H]\times R[H],S)$$

We can do the same for the HNN extension and we get a fibration:
$$\underline Wh^R(H)\build\fl_{}^f \underline Wh^R(G_1)\fl\underline Wh^R(G)'$$
and a homotopy equivalence:
$$\underline Wh^R(G)\simeq\underline Wh^R(G)'\oplus\Omega^{-1}\underline Nil(R[H]\times R[H],S)$$
for some $R[H]\times R[H]$ bimodule $S$.

Hence the only thing to do is to prove that $\underline Nil(R[H]\times R[H],S)$ is contractible.
\vskip 12pt
Let's denote by $C$ the ring $R[H]$. With the notations of 1.4, the bimodule $S$ is determined by four $C$-bimodules ${}_i S_j$ (with $i,j\in\{1,2\}$). In order to
describe these bimodules we'll introduce the following terminology:

Let $G$ be a group. A $G$-biset is a set $X$ equipped with two compatible actions of $G$, one on the left and the other one on the right. We say that a $G$-biset $X$ is
free if both actions on $X$ are free. If $R$ is a ring, $R[X]$ is naturally a $R[G]$-bimodule and, if $X$ is free, $R[X]$ is free on both sides. For any free $G$-biset $X$
and any ring $R$, we set:
$$\underline Nil^R(G,X)=\underline Nil(R[G],R[X])$$

If $G_1$ and $G_2$ are two groups, we can also define a $(G_1,G_2)$-biset as a set equipped with two compatible actions: a left action of $G_1$ and a right action of
$G_2$. Then, for every ring $R$ and any $(G_1,G_2)$-biset $X$, $R[X]$ is a $(R[G_1],R[G_2])$-bimodule. 
\vskip 12pt
Consider the amalgamated case. We have two monomorphisms $\alpha:H\fl G_1$ and $\beta:H\fl G_2$. Denote by $X$ the complement of $\alpha(H)$ in $G_1$ and by $Y$ the
complement of $\beta(H)$ in $G_2$. The group $H$ acts on both sides on $X$ and $Y$ and $X$ and $Y$ are free $H$-bisets. Moreover we have:
$${}_2 S_1=R[X]\hskip 48pt {}_1 S_2=R[Y]\hskip 48pt{}_1S_1={}_2 S_2=0$$
and then:
$$S={}^2 R[X]^1\oplus {}^1 R[Y]^2$$

In the HNN extension case we have two monomorphisms $\alpha:H\fl G_1$ and $\beta:H\fl G_1$. Denote by $X$ the complement of $\alpha(H)$ in $G_1$ and by $Y$ the
complement of $\beta(H)$ in $G_1$. Denote also by $U$ (resp. V) the set $G_1$ where $H$ acts on the left by $\beta$ and on the right by
$\alpha$ (resp. $H$ acts on the right by $\beta$ and on the left by $\alpha$). These sets $X$, $Y$, $U$ and $V$ are free $H$-bisets.
In this case, the bimodule $S$ is characterized by the conditions:
$${}_2 S_1=R[X]\hskip 48pt {}_1 S_2=R[Y]\hskip 48pt{}_1S_1=R[U]\hskip 48pt{}_2 S_2=R[V]$$
and then:
$$S={}^2 R[X]^1\oplus {}^1 R[Y]^2\oplus {}^1 R[U]^1\oplus {}^2 R[V]^2$$

Consider the HNN extension case. For every $x\in G_1$ we set:
$$\Gamma(x)=\beta^{-1}(x\alpha(H)x^{-1})\hskip 48pt \Gamma'(x)=\alpha^{-1}(x\beta(H)x^{-1})$$
For each $x\in G_1$, $\Gamma(x)$ and $\Gamma'(x)$ are subgroups of $H$ and $\Gamma(x)$ is regular coherent. Moreover, we have a group homomorphism $\lambda_x:\Gamma(x)
\fl H$ such that:
$$\forall \gamma\in\Gamma(x),\ \ \beta(\gamma)=x\alpha(\lambda_x(\gamma))x^{-1}$$
It is easy to see that $\lambda_x$ is an isomorphism from $\Gamma(x)$ to $\Gamma'(x^{-1})$. Thus groups $\Gamma'(x)$ are also regular coherent.

Let:
$$U=\ \build{\scriptstyle\coprod}_i^{} U_i$$
be the decomposition of $U$ by orbits. Then, for every $i$, there exists an element $x\in U$ such that:
$$U_i=\beta(H)x\alpha(H)$$
Let $H_1$ (resp. ${}_xH$) be the $(H,\Gamma(x))$-biset (resp. the $(\Gamma(x),H)$-biset) $H$, where $H$ acts in the standard way on the left (resp. on the right) and
$\Gamma(x)$ acts by the inclusion on the right (resp. by the morphism $\lambda_x$ on the left). Then the map: $(u,v)\mapsto\beta(u)x\alpha(v)$ from  $H\times H$ to 
$\beta(H)x\alpha(H)$ induces an isomorphism of $H$-bisets:
$$H_1\build\times_{\Gamma(x)}^{} {}_xH\simeq \beta(H)x\alpha(H)$$
and we have:
$$R[U_i]\simeq R[H_1\build\times_{\Gamma(x)}^{} {}_xH]$$
$$\Longrightarrow {}^1 R[U_i]^1\simeq {}^1R[H_1]\build\otimes_{R[\Gamma(x)]}^{} R[{}_xH]^1$$
Moreover the ring $R[\Gamma(x)]$ is regular coherent.

Hence, because of theorem 6, we have a homotopy equivalence of spectra:
$$\underline Nil(C\times C,{}^2 R[X]^1\oplus {}^1 R[Y]^2\oplus {}^2 R[V]^2)\build\fl_{}^\sim\hskip 48pt$$
$$\hskip 30pt\underline Nil(C\times C,{}^2 R[X]^1\oplus {}^1 R[Y]^2\oplus {}^1 R[U]^1\oplus {}^2 R[V]^2)$$

We proceed the same with the biset $V$ and we get a homotopy equivalence of spectra:
$$\underline Nil(C\times C,{}^2 R[X]^1\oplus {}^1 R[Y]^2)\build\fl_{}^\sim\underline Nil(C\times C,{}^2 R[X]^1\oplus {}^1 R[Y]^2\oplus {}^2 R[V]^2)$$

Hence, in both amalgamated case and HNN extension case, we have a homotopy equivalence:
$$\underline Nil(C\times C,{}^2 R[X]^1\oplus {}^1 R[Y]^2))\build\fl_{}^\sim\underline Nil(C\times C,S)$$
and, because of theorem 4, we have a homotopy equivalence:
$$\underline Nil(C\times C,S)\simeq \underline Nil(C,R[X\build\times_H^{}Y])$$
 
Denote by $Z_j$ the orbits of the biset $X\build\times_X^{}Y$. Then, for each $j$, there is an element $(x,y)\in X\times Y$ such that:
$$Z_j=\alpha(H)xy\beta(H)$$

For each $x\in X$ and each $y\in Y$ we have the groups:
$$\Gamma_1(x)=\alpha^{-1}(x\alpha(H)x^{-1})\hskip 30pt \Gamma_2(y)=\beta^{-1}(y\beta(H)y^{-1})\hskip 30pt H(x,y)=\Gamma_1(x)\cap\Gamma_2(y)$$
We have group homomorphisms $\lambda_x:\Gamma_1(x)\fl H$ and $\mu_y:\Gamma_2(y)\fl H$ defined by:
$$\forall \gamma\in\Gamma_1(x),\ \ \alpha(\lambda_x(\gamma))=x\alpha(\gamma)x^{-1}$$
$$\forall \gamma\in\Gamma_2(y),\ \ \beta(\gamma)=y\beta(\mu_y(\gamma))y^{-1}$$

Denote by $H_x$ the $(H,H(x,y))$-biset $H$ where $H$ acts in the standard way on the left and $H(x,y)$ acts via $\lambda_x$ on the right. Denote also by ${}_y H$ the
$(H(x,y),H)$-biset where $H$ acts in the standard way on the right and $H(x,y)$ acts via $\mu_y$ in the left. Then the map $(u,v)\mapsto \alpha(u)xy\beta(v)$
from $H\times H$ to $\alpha(H)xy\beta(H)$ induces an isomorphism:
$$H_x\build\times_{H(x,y)}^{}{}_y H\build\fl_{}^\sim\alpha(H)xy\beta(H)$$
where $H(x,y)$ is regular coherent. Hence, because of theorem 6, the spectrum $\underline Nil(R[H],R[X\build\times_H^{}Y])$ is contractible and so is
$\underline Nil(R[H]\times R[H],S)$.\cqfd
\vskip 12pt
\noi{\bf 4.6 Theorem:} {\sl For every group $G$ in Cl$_1$ and every ring $R$ which is regular noetherian on the right, the Whitehead spectrum $\underline Wh^R(G)$ is
contractible.}
\vskip 12pt
\noi{\bf Proof:} Denote by Cl$_2$ the class of groups $G$ such that $\underline Wh^R(G)$ is contractible for every ring $R$ which is regular noetherian on the right.
Because of theorem 4.2, Cl$_2$ contains the class Cl.

Since the functor $\underline Wh$ commutes with filtered colimits, the class Cl$_2$ is stable under filtered colimits. Therefore it's enough to prove that Cl$_2$ is
stable under taking amalgamated free products and HNN extensions of regular coherent diagrams. But that follows directly from theorems 4.3 and 4.4.
\cqfd
\vskip 12pt
\noi{\bf 4.7 Example:} Consider the group $H$ with two generators $x$ and $t$ and the following relations:
$$\forall n\in\Z,\ \ xt^nxt^{-n}=t^nxt^{-n}x$$

For every integer $p\not=0$, the correspondence $x\mapsto x$ and $t\mapsto t^p$ induces a monomorphism $f_p:H\fl H$. Consider a monomorphism of groups $\alpha:H\fl G$ and
denote by $\Gamma$ the amalgamated free product of the diagram:
$$\diagram{H&\hfl{f_p}{}&H\cr\vfl{\alpha}{}&&\cr G&&\cr}$$
\vskip 12pt
\noi{\bf 4.8 Proposition:} {\sl For every ring $R$ which is regular noetherian on the right, the morphism $G\fl \Gamma$ induces a homotopy equivalence of spectra:
  $\underline Wh^R(G)\build\fl_{}^\sim\underline Wh^R(\Gamma)$.

  Moreover, if $G$ belongs to Cl$_1$, then the group $\Gamma$ is also in Cl$_1$.} 
\vskip 12pt
\noi{\bf Proof:}  Denote by $H'$ the normal closure of $x$ in $H$. This group is commutative and freely generated by the elements: $x_n=t^n xt^{-n}$ for $n\in\Z$. The
correspondence. $x_n\mapsto x_{n+1}$ is an automorphism $\tau:H'\build\fl_{}^\sim H'$ and $H$ is the semidirect product of $H'$ and $\Z$, or equivalently, the HNN
extension of $H'$ with morphisms Id, $\tau:H'\fl H'$. On the other hand, $R[H']$ is regular coherent on the right (but not noetherian) and $H$ belongs to the class Cl.
Hence the Whitehead spectrum $\underline Wh^R(H)$ is contractible.

Denote by $H_p$ the image of $f_p:H\fl H$ and by $X$ its complement in $H$. For every $z\in H$, denote by $\Gamma(z)$ the subgroup $f_p^{-1}(zf_p(H)z^{-1})$ of
$H$. We have the following formula:
$$\Gamma(f_p(a)zf_p(b))=a\Gamma(z)a^{-1}$$
for every $a,b,z$ in $H$ and the conjugacy class of $\Gamma(z)$ depends only on the class of $z$ in the set $Y=H_p\backslash H/H_p$. Let $z$ be an element in $X$.
A direct computation shows the following:

$\bullet$ if $z$ is congruent in $Y$ to an element in $H'$ then $\Gamma(z)$ is the group $H'$

$\bullet$ if $z$ is congruent in $Y$ to a power of $t$ then $\Gamma(z)$ is conjugate to the subgroup of $H$ generated by $t$

$\bullet$ in the other cases $\Gamma(z)$ is the trivial group.

Therefore $\Gamma(z)$ is always a free abelian group. Hence, for every $y$ in $G\setminus\alpha(H)$, the group $\Gamma(z,y)=\Gamma(z)\cap \alpha^{-1}(y\alpha(H)y^{-1})$
is also a free abelian group and the ring $R[\Gamma(z,y)]$ is regular coherent on the right. Then theorem 4.3 applies and the result follows.\cqfd
\vskip 12pt
The class Cl$_1$ seems to be strictly bigger than the class Cl. For example the amalgamated free product $\Gamma$ of the diagram:
$$\diagram{H&\hfl{f_p}{}&H\cr\vfl{f_q}{}&&\cr H&&\cr}$$
with $p,q>1$, belongs to the class Cl$_1$. But in this case, Waldhausen's theorems cannot be used to prove that $\Gamma$ belongs to the class Cl because of the following
result:
\vskip 12pt
\noi{\bf 4.9 Proposition:} {\sl The ring $\Z[H]$ is not regular coherent.}
\vskip 12pt
\noi{\bf Proof:} Let $f:\Z[H]\oplus\Z[H]\fl \Z[H]$ be the following morphism:
$$(U,V)\mapsto f(U,V)=(1-t+tx)U-(1-t+t^2xt^{-1})V$$
and $K$ be its kernel. We'll prove that $K$ is not finitely generated and that will imply that $\Z[H]$ is not coherent and therefore not regular coherent.

Denote by $A$ the ring $\Z[H']$. This ring is the ring of Laurent polynomials in the $x_i$'s. Then $A$ is an integral domain and every element $u\in\Z[H]$ can be written
in a unique way on a finite sum:
$$u=\build\sum_{i\in\Z}^{} t^iu_i$$
with each $u_i$ in $A$. So $u$ may be considered as a Laurent polynomial in $t$ and has a valuation $\nu(u)$ and a degree $\partial^\circ u$ (at least if $u$ is
not zero). If $u=0$, we set: $\nu(u)=+\infty$ and $\partial^\circ u=-\infty$.

Define the elements $y_i$ and $z_i$ in $A$ by:
$$\forall i\in\Z,\ \ y_i=1-x_i=1-t^i x t^{-i}\ \ \ z_i=y_i-y_{i-1}=x_{i-1}-x_i$$
and, for every integer $n\geq0$, we have the following elements in $\Z[H]$:
$$U_n=z_{-n}-t^{n+1}z_1y_0y_{-1}\dots y_{-n}$$
$$V_n=z_{-n}+\build\sum_{0<i\leq n}^{}t^iz_{-n}z_1y_0y_{-1}\dots y_{2-i}-t^{n+1}z_1y_0y_{-1}\dots y_{1-n}y_{-1-n}$$

An explicit computation shows that, for each $n\geq0$, $W_n=(U_n,V_n)$ is killed by $f$ and each $W_n$ belongs to $K$. 
\vskip 12pt
\noi{\bf 4.10 Lemma:} {\sl The set $\{W_0,W_1,W_2,\dots\}$ is a generating set of $K$.}
\vskip 12pt
\noi{\bf Proof:} Denote by $E$ the $\Z[H]$-submodule of $\Z[H]\oplus\Z[H]$ generated by the set $\{W_0,W_1,W_2,\dots\}$. Since each $W_n$ is in $K$, we have an inclusion
$E\subset K$ and we have to prove that this inclusion is an equality.

For each integer $n\geq0$, denote by $K_n$ the set of the elements $(U,V)$ in $K$ satisfying the following:
$$\nu(U)\geq0\hskip 24pt\partial^\circ U\leq n\hskip 48pt \nu(V)\geq0\hskip 24pt\partial^\circ V\leq n$$
We denote also by $I_n$ the ideal $(z_0,z_{-1},z_{-2},\dots,z_{1-n})\subset A$ and by $J_n$ the right ideal of $\Z[H]$ generated by $I_n$.

The quotient $B_n=A/I_n$ is the quotient of $A$ by the relations $$x_0=x_{-1}=x_{-2}=\dots=x_{-n}$$ and $B_n$ is a Laurent polynomial ring where $z_1$, $z_{-n}$ and all
$y_i$ are not zero.

Suppose we have proven that $K_{n-1}$ is contained in $E$. Let $U=\build\sum_{0\leq i\leq n}^{}t^i u_i$ and $V=\build\sum_{0\leq i\leq n}^{}t^i v_i$ be two elements
in $\Z[H]$, with $u_i$ and $v_i$ in $A$. For $W=(U,V)$ we have the following equivalences:
$$W\in K_n\ \Longleftrightarrow\ (1-ty_0)U=(1-ty_1)V$$
$$\Longleftrightarrow\ \build\sum_i^{} t^i(u_i-v_i)=ty_0\build\sum_i^{} t^iu_i-ty_1\build\sum_i^{} t^iv_i$$
$$\Longleftrightarrow\ \build\sum_i^{} t^i(u_i-v_i)=\build\sum_i^{}t^{i+1}y_{-i}u_i-\build\sum_i^{}t^{i+1}y_{1-i}v_i$$
$$\Longleftrightarrow\ \forall i,\ u_i-v_i=y_{1-i}u_{i-1}-y_{2-i}v_{i-1}$$
And these conditions are equivalent to the following:
$$v_0=u_0$$
$$v_1=u_1+z_1u_0$$
$$v_2=u_2+z_0u_1+z_1y_0u_0$$
$$\dots$$
$$v_n=u_n+\build\sum_{0\leq i<n}^{} z_{1-i}y_{-i}y_{-1-i}\dots y_{2-n}u_i$$
$$0=\build\sum_{0\leq i\leq n}^{} z_{1-i}y_{-i}y_{-1-i}\dots y_{1-n}u_i$$

This last relation implies the following relation in $B_n$:
$$z_1y_0y_{-1}\dots y_{1-n}u_0\equiv 0\ \in\ B_n$$
and we have: $u_0\equiv 0$ in $B_n$ because $B_n$ is an integral domain.

Suppose $W$ is in $K_n$. Then $u_0$ is in $I_n$ and there are elements $a_0,a_1,\dots a_{n-1}$ in $A$ such that:
$$u_0=z_0 a_0+z_{-1}a_1+\dots+z_{1-n}a_{n-1}$$
Set:
$$W'=W-(W_0 a_0+W_1 a_1+W_2 a_2+\dots W_{n-1}a_{n-1})$$
Since $W_0,W_1,\dots, W_{n-1}$ are in $K_n$, $W'$ belongs to $K_n$ and there exist elements $u'_i$ and $v'_i$ in $A$ such that:
$$W'=(\build\sum_{0\leq i\leq n}^{}t^i u'_i,\build\sum_{0\leq i\leq n}^{}t^i v'_i)$$
Moreover, because $W'$ is in $K_n$, we have $u'_0=v'_0$ and:
$$u'_0=u_0-(z_0 a_0+z_{-1}a_1+\dots+z_{1-n}a_{n-1})=0$$
So we have: $u'_0=v'_0=0$ and $W't^{-1}$ belongs to $K_{n-1}\subset E$. Therefore $W'$ and then $W$ belong to $E$ and we have: $K_n\subset E$.

Thus $K_n$ is contained in $E$ for every $n\geq0$. On the other hand, for every $W\in K$, there is some integer $p$ such that $W t^p$ belongs to some $K_n$. Hence $K$ is
contained in $E$ and the lemma is proven.\cqfd
\vskip 12pt
\noi{\bf 4.11 Lemma:} {\sl The module $K$ is not finitely generated.}
\vskip 12pt
\noi{\bf Proof:} For each integer $n>0$, denote by $E_n$ the submodule of $E$ generated by $\{W_0,W_1,W_2,\dots ,W_{n-1}\}$.

Let $F_n:\Z[H]^n\fl \Z[H]\oplus\Z[H]$ be the morphism:
$$(c_0,c_1,c_2,\dots,c_{n-1})\mapsto \build\sum_{0\leq i<n}^{}W_ic_i$$

The image if $F_n$ is the module $E_n$. Denote by $R_n$ the kernel of $F_n$ and by $\pi_n$ the last projection $\Z[H]^n\fl \Z[H]$. We have an exact sequence of right
$\Z[H]$-modules:
$$0\fl R_n\fl \Z[H]^n\build\fl_{}^{F_n}E_n\fl 0$$
Hence we have a commutative diagram of right $\Z[H]$-modules with exact lines:
$$\diagram{0&\hfl{}{}&\Z[H]^{n-1}&\hfl{}{}&\Z[H]^n&\hfl{\pi_n}{}&\Z[H]&\hfl{}{}&0\cr &&\vfl{F_{n-1}}{}&&\vfl{F_n}{}&&\vfl{}{}&&\cr
  0&\hfl{}{}&E_{n-1}&\hfl{}{}&E_n&\hfl{}{}&\Z[H]/J'_n&\hfl{}{}&0\cr}$$
with: $J'_n=\pi_n(R_n)$. We'll prove that $\Z[H]/J'_n$ is not trivial.

An easy computation shows that, for every integers $p,q$ with $0\leq p<q$, the following element:
$$W_qz_{-p}-W_p z_{-q}-W_{q-p-1}t^{p+1}z_1y_0y_{-1}\dots y_{-p}$$
is zero in $K$ and induces a well defined element $X(p,q)\in R_n$ (for every $n>q$) and, for every $p$ with $0\leq p<n-1$, we have: $\pi_n(X(p,n-1))=z_{-p}$.

Let $n>0$ be an integer. Suppose $J'_n$ is not contained in $J_n$. Then there is an element $X_0\in R_n$ such that $\pi_n(X_0)$ doesn't belong to $J_n$. Set:
$d=\partial^\circ\pi_n(X_0)$ and denote by $Z$ the set of $X\in R_n$ such that: $\partial^\circ\pi_n(X)=d$ and $\pi_n(X)-\pi_n(X_0)\in J_n$.

For each $X=(c_0,c_1,\dots,c_{n-1})\in Z$ we can associate three integers $\alpha,\beta,\gamma$ defined this way:

$\bullet$ $\alpha=\nu(c_{n-1})$

$\bullet$ $\beta$ is the lowest $\nu(c_k)$, for $k=0,1,\dots,n-1$

$\bullet$ $\gamma$ is the highest integer $k$ such that $\nu(c_k)=\beta$.

The triple $\chi(X)=(\alpha,\beta,\gamma)$ will be called the complexity of $X$. This complexity belongs to the set $C$ of triple $(\alpha,\beta,\gamma)\in\Z^3$ satisfying
the following conditions:
$$\beta\leq\alpha\leq d\ \ \hbox{and}\ \ 0\leq\gamma<n$$

The lexicographical order of $(d-\alpha,\alpha-\beta,\gamma)$ induces a well order relation on $C$ and we have:
$$(\alpha,\beta,\gamma)<(\alpha',\beta',\gamma')\Longleftrightarrow \alpha>\alpha'\ \hbox{or}\ \alpha=\alpha'\ \hbox{and}\ \beta>\beta'\ \hbox{or}\
(\alpha,\beta)=(\alpha',\beta')\ \hbox{and}\ \gamma<\gamma'$$

Since $C$ is well ordered, there is an element in $Z$ with a minimal complexity. Let $X=(c_0,c_1,\dots,c_{n-1})$ be such an element.

For each integer $k\in\{0,1,\dots,n-1\}$, we have a decomposition:
$$c_k=\build\sum_{\beta\leq i}^{}c_{ki}t^i$$
with $c_{ki}\in A$.

The condition $X\in R_n$ implies the following:
$$\build\sum_{0\leq k\leq \gamma}^{} z_{-k}c_{k\beta}=0$$
and that implies the congruence $z_{-\gamma}c_{\gamma\beta}\equiv0$ in $B_{\gamma}=A/I_{\gamma}$. But $z_{-\gamma}$ is not a zero divisor in $B_{\gamma}$ and
$c_{\gamma\beta}$ belongs to $I_{\gamma}$. So we have a decomposition in $A$:
$$c_{\gamma\beta}=\build\sum_{0\leq j<\gamma}^{} z_{-j}a_j$$ 
and we get a new element in $R_n$:
$$X'=X-\build\sum_{0\leq j<\gamma}^{} X(j,\gamma) a_jt^\beta$$

It is easy to see that $\pi_n(X')\equiv\pi_n(X)\equiv\pi_n(C_0)$ mod $J_n$ and that $X'$ belongs to $Z$.
Moreover we have the following: $\chi(X')<\chi(X)$. But that's impossible because $X$ was chosen with a minimal complexity.

Hence we get a contradiction and the module $J'_n$ is contained in $J_n$. As a consequence, by killing all the $z_i$'s, we get epimorphisms:
$$\Z[H]/J'_n\fl\Z[H]/J_n\fl \Z[x^{\pm1},t^{\pm1}]$$
and $\Z[H]/J'_n$ is not trivial. Hence the sequence $E_0\subset E_1\subset E_2\subset\dots$ is strictly increasing and $E=$Ker$(f)$ is not finitely generated. Therefore
the category of finitely presented right $\Z[H]$-modules is not abelian and $\Z[H]$ is not coherent. The result follows.\cqfd

\vskip 48pt
\noi{\bf References: }

\begin{list}{}{\leftmargin 24pt \labelsep 10pt \labelwidth 40pt \itemsep 0pt}
\item[{[B]}] H. Bass -- {\sl Algebraic K-theory}, Benjamin (1968).
\item[{[C]}]  P. M. Cohn -- {\sl Free ideal rings}, J. Algebra {\bf 1} (1964) 47--69.
\item[{[HW]}] A. Hatcher and J. Wagoner -- {\sl Pseudo isotopies of compact manifolds}, Ast\'erisque {\bf 6}, (1973).
\item[{[K]}] B. Keller -- {\sl Chain complexes and stable categories}, Manuscripta Mathematica. {\bf 67} (1990), 379–-417. doi:10.1007/BF02568439.
\item[{[KV]}] M. Karoubi and O. Villamayor -- {\sl K-th\'eorie alg\'ebrique et K-th\'eorie topologique}, Math. Scand. {\bf 28} (1971), 265--307.
\item[{[Q]}] D. Quillen -- {\sl Higher algebraic K-theory I}, Proc. Conf. alg. K-theory, Lecture Notes in Math. {\bf 341} (1973), 85--147.
\item[{[TT]}] R. W. Thomason and T. Trobaugh -- {\sl Higher algebraic K-Theory of schemes and of derived categories}, The Grothendieck Festschrift III, Progress in Math.,
  {\bf 88}, Birkh\"auser Boston, Boston, MA (1990) 247–-435. MR1106918
\item[{[Wa1]}] F. Waldhausen -- {\sl Algebraic K-theory of generalized free products, part 1 \& 2}, Annals of Math. {\bf 108} (1978), 135--256.
\item[{[Wa2]}] F. Waldhausen -- {\sl Algebraic K-theory of spaces}, Algebraic and geometric topology (New Brunswick, N.J., 1983), 318–-419, 
  Lecture Notes in Math. {\bf 1126} Springer, Berlin, (1985).
\item[{[We]}] C. A. Weibel -- {\sl The K-book: An introduction in algebraic K-theory}, Graduate Studies in Mathematics, {\bf 145}, Amer. Math. Soc. (2013).
\end{list}
\end{document}